\documentclass{amsart}
\usepackage{amsmath,amsrefs}

\usepackage[colorlinks=true, pdfborder={0 0 0}]{hyperref}
\usepackage{cleveref}
\usepackage{upref}

\DeclareMathOperator{\Scal}{Scal}
\DeclareMathOperator{\Weyl}{Weyl}
\DeclareMathOperator{\Span}{Span}
\DeclareMathOperator{\Eucl}{Eucl}
\DeclareMathOperator{\loc}{loc}

\DeclareMathOperator{\bigO}{O}
\DeclareMathOperator{\smallo}{o}
\def \rr {\mathbb{R}}
\def \nn {\mathbb{N}}
\def \rn {\mathbb{R}^n}
\def \ss {\mathbb{S}}

\def \ue {u_\epsilon}
\def \he {h_\epsilon}

\def \hue {\hat{u}_\epsilon}
\def \hhe {\hat{h}_\epsilon}

\def \tu {\tilde{u}_0}

\def \eps {\epsilon} 

\def \crit {2^\star}

\newtheorem{theorem}{Theorem}[section]
\newtheorem{proposition}{Proposition}[section]
\newtheorem{lem}{Lemma}[section]
\newtheorem{propdefi}{Proposition-Definition}[section]

\begin{document}

\title[Blowing-up solutions]{Blowing-up solutions for second-order critical elliptic equations: the impact of the scalar curvature}

\author{Fr\'ed\'eric Robert}

\address{Fr\'ed\'eric Robert, Institut \'Elie Cartan, UMR 7502, Universit\'e de Lorraine, BP 70239, F-54506 Vand{\oe}uvre-l\`es-Nancy, France}
\email{frederic.robert@univ-lorraine.fr}

\author{J\'er\^ome V\'etois}

\address{J\'er\^ome V\'etois, McGill University, Department of Mathematics and Statistics, 805 Sherbrooke Street West, Montreal, Quebec H3A 0B9, Canada}
\email{jerome.vetois@mcgill.ca}

\date{December 16, 2019}

\begin{abstract} Given a closed manifold $(M^n,g)$, $n\geq 3$, Olivier Druet~\cite{DruetJDG} proved that a necessary condition for the existence of energy-bounded blowing-up solutions to perturbations of the equation
$$\Delta_gu+h_0u=u^{\frac{n+2}{n-2}},\ u>0\hbox{ in }M$$
is that $h_0\in C^1(M)$ touches the Scalar curvature somewhere when $n\geq 4$ (the condition is different for $n=6$). In this paper, we prove that Druet's condition is also sufficient provided we add its natural differentiable version.  For $n\geq 6$, our arguments are local. For the low dimensions $n\in\{4,5\}$, our proof requires the introduction of a suitable mass that is defined only where Druet's condition holds. This mass carries global information both on $h_0$ and $(M,g)$.

\end{abstract}

\maketitle
\section{Introduction and main results}\label{Sec1}
Let $(M,g)$ be a compact Riemannian manifold of dimension $n\geq 3$ without boundary and $h_0\in C^p(M)$, $1\le p\le\infty$. We consider the equation
\begin{equation}\label{eq}
\Delta_gu+h_0u=u^{\crit-1},\ u>0\hbox{ in }M,
\end{equation}
where $\Delta_g:=-\hbox{div}_g(\nabla)$ is the Laplace--Beltrami operator and $\crit:=\frac{2n}{n-2}$. We investigate the existence of families $(h_\eps)_{\eps>0}\in C^p(M)$ and $(u_\eps)_{\eps>0}\in C^2(M)$ satisfying 
\begin{equation}\label{eque}
\Delta_gu_\eps+h_\eps u_\eps=u_\eps^{\crit-1},\ u_\eps>0\hbox{ in }M\text{ for all }\eps>0,
\end{equation}
and such that $h_\eps\to h_0$ in $C^p(M)$ and $\max_Mu_\eps\to\infty$ as $\eps\to0$. We say that $(u_\eps)_{\eps}$ {\it blows up} at some point $\xi_0\in M$ as $\eps\to0$ if for all $r>0$, $\lim_{\eps\to 0}\max_{B_r(\xi_0)}\ue=+\infty$. Druet~\cites{DruetJDG,DruetHDR} obtained the following necessary condition for blow-up:

\begin{theorem}[Druet~\cites{DruetJDG,DruetHDR}]\label{ThDruet}
Let $(M,g)$ be a compact Riemannian manifold of dimension $n\geq 4$. Let $h_0\in C^1(M)$ be such that $\Delta_g+h_0$ is coercive. Assume that there exist families $(h_\eps)_{\eps>0}\in C^1(M)$ and $(u_\eps)_{\eps>0}\in C^2(M)$ satisfying \eqref{eque} and such that $h_\eps\to h_0$ strongly in $C^1(M)$ and $u_\eps\rightharpoonup u_0$ weakly in $L^{\crit}(M)$. Assume that $(u_\eps)_\eps$ blows-up. Then there exists $\xi_0\in M$ such that $(u_\eps)_{\eps}$ blows up at $\xi_0$ and 
\begin{equation}\label{phi0}
\left(h_0-c_n\Scal_g\right)(\xi_0)=0\hbox{ if }n\neq 6\hbox{ and }\left(h_0-c_n\Scal_g-2u_0\right)(\xi_0)=0\hbox{ if }n= 6.
\end{equation}

Furthermore, if $n\in\{4,5\}$, then $u_0\equiv0$. \end{theorem}
Here $c_n:=\frac{n-2}{4(n-1)}$ and $\Scal_g$ is the Scalar curvature of $(M,g)$. This result does not hold in dimension $n=3$. Indeed, Hebey--Wei~\cite{HW} constructed examples of blowing-up solutions to \eqref{eque} on the standard sphere $(\ss^3,g_0)$, which are bounded in $L^{\crit}(\ss^3)$ but do not satisfy \eqref{phi0}.

\smallskip
This paper is concerned with the converse of Theorem~\ref{ThDruet} in dimensions $n\ge4$. For the sake of clarity, we state separately our results in the cases  $u_0\equiv0$ in dimension $n\ge4$ (Theorem~\ref{th:2}) and  $u_0>0$ in dimension $n\ge6$ (Theorem~\ref{th:3}):

\begin{theorem}\label{th:2} 
Let $(M,g)$ be a compact Riemannian manifold of dimension $n\geq 4$. Let $h_0\in C^p(M)$, $1\le p\le\infty$, be such that $\Delta_g+h_0$ is coercive. Assume that there exists a point $\xi_0\in M$ such that 
\begin{equation}\label{phi0th2}
\left(h_0-c_n\Scal_g\right)(\xi_0)=|\nabla\left(h_0-c_n\Scal_g\right)(\xi_0)|=0.
\end{equation}
Then there exist families $(h_\eps)_{\eps>0}\in C^p(M)$ and $(u_\eps)_{\eps>0}\in C^2(M)$ satisfying \eqref{eque} and such that $h_\eps\to h_0$ strongly in $C^p(M)$, $u_\eps\rightharpoonup 0$ weakly in $L^{\crit}(M)$ and $(u_\eps)_{\eps}$ blows up at $\xi_0$.
\end{theorem}
For convenience, for every $h_0,u_0\in C^0(M)$, we define
\begin{equation}\label{defphi0th3}
\varphi_{h_0}:=h_0-c_n\Scal_g\hbox{ and }\varphi_{h_0,u_0}:=\left\{\begin{aligned}&h_0-c_n\Scal_g&&\text{if }n\ne6\\&h_0-2u_0-c_n\Scal_g&&\text{if }n=6.\end{aligned}\right.
\end{equation}

\begin{theorem}\label{th:3}
Let $(M,g)$ be a compact Riemannian manifold of dimension $n\geq 6$. Let $h_0\in C^p(M)$, $1\le p\le\infty$, be such that $\Delta_g+h_0$ is coercive. Assume that there exist a solution $u_0\in C^2(M)$ of \eqref{eq} and a point $\xi_0\in M$ such that
\begin{equation}\label{phi0th3}
\varphi_{h_0,u_0}(\xi_0)=|\nabla\varphi_{h_0,u_0}(\xi_0)|=0.
\end{equation}
Then there exist families $(h_\eps)_{\eps>0}\in C^p(M)$ and $(u_\eps)_{\eps>0}\in C^2(M)$ satisfying \eqref{eque} and such that $h_\eps\to h_0$ strongly in $C^p(M)$, $u_\eps\rightharpoonup u_0$ weakly in $L^{\crit}(M)$ and $(u_\eps)_{\eps}$ blows up at $\xi_0$.
\end{theorem}
Compared with Theorem \ref{ThDruet}, we have assumed here that condition \eqref{phi0} is also satisfied at order 1. However, this stronger condition is actually expected to be necessary for the existence of blowing-up solutions (see Theorem \ref{th:cns} in the last section of this paper and the discussion in Druet~\cite{DruetHDR}*{Section~2.5}). 

\smallskip
We refer to Section~\ref{Sec2} for examples of functions $h_0$ and $u_0$ satisfying the assumptions in Theorem~\ref{th:3}. Recently, Premoselli--Thizy~\cite{PT} obtained a beautiful example of blowing-up solutions showing that in dimension $n\in\left\{4,5\right\}$,  condition \eqref{phi0th2} may not be satisfied at all blow-up points. 

\smallskip
When $h_0\equiv c_n\Scal_g$, that is when \eqref{eq} is the Yamabe equation, several examples of blowing-up solutions have been obtained. In the perturbative case, that is when $h_\eps\not\equiv c_n\Scal_g$, examples of blowing-up solutions have been obtained by Druet--Hebey~\cite{DH}, Esposito--Pistoia--V\'etois~\cite{EPV}, Morabito--Pistoia--Vaira~\cite{MPV}, Pistoia--Vaira~\cite{PV} and Robert--V\'etois~\cite{RVJDG}. In the nonpertubative case $h_\eps\equiv c_n\Scal_g$, we refer to Brendle~\cite{Brendle} and Brendle--Marques~\cite{BMJDG} regarding the non-compactness of Yamabe metrics. When solutions blow-up not only pointwise but also in energy, the function $\varphi_{h_0}$ may not vanish (see Chen--Wei--Yan \cite{CWY} for $n\ge5$ and V\'etois--Wang~\cite{VW} for $n=4$).

\smallskip
When there does not exist any blowing-up solutions to the equations~\eqref{eque}, then equation \eqref{eq} is {\it stable}. We refer to the survey of Druet~\cites{DruetHDR} and the book of Hebey~\cite{HebeyETH} for exhaustive studies of the various concepts of stability. Stability also arises in the Lin--Ni--Takagi problem (see for instance del Pino--Musso--Roman--Wei~\cite{delPMRW} for a recent reference on this topic).  In Geometry, stability is linked to the problem of compactness of the Yamabe equation (see Schoen~\cites{Schoen2,Schoen3}, Li--Zhu~\cite{LZ3}, Druet~\cite{DruetIMRN}, Marques~\cite{Marques}, Li--Zhang~\cites{LZ1,LZ2}, Khuri--Marques--Schoen~\cite{KMS}). 

\smallskip
Let us give some general considerations about the proofs. Theorem~\ref{ThDruet} yields {\it local} information on blow-up points. It is essentially the consequence of the concentration of the $L^2$--norm of the solutions at one of the blow-up points when $n\geq 4$. However, in our construction, the problem may be {\it both local and global}. Indeed, we reduce the problem to finding critical points of a functional defined on a finite-dimensional space. The first term in the asymptotic expansion of the reduced functional is local. This is due to the $L^2$--concentration of the standard bubble in the definition of our ansatz. The second term in the expansion plays a decisive role for obtaining critical points. For the high dimensions $n\geq 6$, this term is also local (see e.g. \eqref{est:J:W:2}). However, for $n\in\left\{4,5\right\}$, the second term is global and we are then compelled to introduce a suitable notion of mass, which carries global information on $h_0$ and $(M,g)$, and to add a corrective term to the standard bubble (see \eqref{def:W:beta}) in order to obtain a reasonable expansion (see e.g. \eqref{est:J:1:2}). Unlike the case where $n=3$ or $h_0\equiv c_n \Scal_g$, the mass is not defined at all points in the manifold, but only at the points where the condition \eqref{phi0th3} is satisfied.

\smallskip
More precisely, Theorems~\ref{th:2} and~\ref{th:3} are consequences of Theorems~\ref{th:main} and~\ref{th:main:u0} below. The latter are the core results of our paper. In these theorems, we fix a linear perturbation $h_\eps=h_0+\eps f$ for some function $f\in C^p(M)$. Furthermore, we specify the behavior of the blowing-up solutions that we obtain. More precisely, we say that $(\ue)_\eps$ blows up with one bubble at some point $\xi_0\in M$ if $u_\eps=u_0+U_{\delta_\eps,\xi_\eps}+\smallo(1)$ as $\eps\to0$ in $H^2_1(M)$, where $u_0\in H_1^2(M)$ is such that $u_\eps\rightharpoonup u_0$ weakly in $H_1^2(M)$,  $U_{\delta_\eps,\xi_\eps}$ is as in \eqref{def:peak}, $(\delta_\eps,\xi_\eps)\to(0,\xi_0)$ and $o(1)\to 0$ strongly in $H_1^2(M)$ as $\eps\to0$. 

\smallskip
Our first result deals with the case where $u_0\equiv0$ in dimension $n\ge4$:

\begin{theorem}\label{th:main} 
Let $(M,g)$ be a compact Riemannian manifold of dimension $n\geq 4$. Let $h_0\in C^{p}(M)$, $p\geq 2$, be such that $\Delta_g+h_0$ is coercive. Assume that there exists a point $\xi_0\in M$ satisfying \eqref{phi0th2}. Assume in addition that $\xi_0$ is a {\it nondegenerate} critical point of $h_0-c_n\Scal_g$ and
\begin{equation}\label{Kh}
K_{h_0}(\xi_0):=\left\{\begin{aligned}&m_{h_0}(\xi_0)&&\text{if }n=4,5\\&\Delta_g\left(h_0-c_n\Scal_g\right)(\xi_0)+\frac{c_n}{6}|\Weyl_g(\xi_0)|_g^2&&\text{if }n\ge6\end{aligned}\right\}\ne0,
\end{equation}
where $m_{h_0}(\xi_0)$ is the mass of $\Delta_g+h_0$ at the point $\xi_0$ (see Proposition-Definition~\ref{propdefi:1}), and $\Weyl_g$ is the Weyl curvature tensor of the manifold. We fix a function $f\in C^{p}(M)$ such that $f(\xi_0)\times K_{h_0}(\xi_0)>0$. Then for small $\epsilon>0$, there exists $\ue\in C^2(M)$ satisfying
\begin{equation}\label{eq:ue}
\Delta_g \ue+(h_0+\eps f)\ue=\ue^{\crit-1}\hbox{ in }M,\ \ue>0,
\end{equation}
and such that $\ue\rightharpoonup 0$ weakly in $L^{\crit}(M)$ and $(\ue)_\eps$ blows up with one bubble at $\xi_0$.
\end{theorem}
The definition of $K_{h_0}(\xi_0)$ outlines the major difference between high- and low-dimensions that was mentioned above: for $n\geq 6$, it is a local quantity, but for $n\in \{4,5\}$, it carries global information (see Section~\ref{sec:mass} for more discussions). 

\smallskip
Next we deal  with the case where $u_0>0$ in dimension $n\ge6$:

\begin{theorem}\label{th:main:u0} 
Let $(M,g)$ be a compact Riemannian manifold of dimension $n\geq 6$. Let $h_0\in C^{p}(M)$, $p\geq 2$, be such that $\Delta_g+h_0$ is coercive. Assume that there exist a nondegenerate solution $u_0\in C^2(M)$ to equation \eqref{eq} and  $\xi_0\in M$ satisfying \eqref{phi0th3}. Assume in addition that $\xi_0$ is a {\it nondegenerate} critical point of $\varphi_{h_0,u_0}$ and
\begin{equation}\label{Khu}
K_{h_0,u_0}(\xi_0):=\left\{\begin{aligned}&\Delta_g\varphi_{h_0,u_0}(\xi_0)+\frac{c_6}{6}|\Weyl_g(\xi_0)|_g^2&&\hbox{if }n=6\\
&u_0(\xi_0)&&\hbox{if }7\leq n\leq 9\\
&672u_0(\xi_0)+\Delta_g\varphi_{h_0,u_0}(\xi_0)+\frac{c_{10}}{6}|\Weyl_g(\xi_0)|_g^2&&\hbox{if }n=10\\
&\Delta_g\varphi_{h_0,u_0}(\xi_0)+\frac{c_n}{6}|\Weyl_g(\xi_0)|_g^2&&\hbox{if }n\geq 11\end{aligned}\right\}\ne0.
\end{equation}
We fix a function $f\in C^{p}(M)$ such that 
\begin{equation}\label{condition6}
K_{h_0,u_0}(\xi_0)\times\left\{\begin{aligned}&\left[f+2(\Delta_g+h_0-2u_0)^{-1}(fu_0)\right](\xi_0) &&\text{if }n=6\\&f(\xi_0)&&\text{if }n>6\end{aligned}\right\}>0.
\end{equation}
Then for small $\epsilon>0$, there exists $\ue\in C^2(M)$ satisfying \eqref{eq:ue} and such that $\ue\rightharpoonup u_0$ weakly in $L^{\crit}(M)$ and $(\ue)_\eps$ blows up with one bubble at $\xi_0$.
\end{theorem}

The paper is organized as follows. In Section~\ref{Sec2}, we discuss the question of existence of functions $h_0$ and $u_0$ satisfying the assumptions in Theorem~\ref{th:3}. In Section~\ref{Sec3}, we introduce our notations and discuss the general setting of the problem. In Section~\ref{Sec4}, we establish a general $C^1$-estimate on the energy functional, which holds in all dimensions. In Sections~\ref{Sec5},~\ref{Sec6} and~\ref{Sec7}, we then compute a $C^1$-asymptotic expansion of the energy functional in the case where $n\ge6$, which we divide in the following subcases: [$n\ge6$ and $u_0\equiv0$] in Section~\ref{Sec5}, [$n\ge7$ and $u_0>0$] in Section~\ref{Sec6} and [$n=6$ and $u_0>0$] in Section~\ref{Sec7}. In Section~\ref{sec:mass}, we discuss the specific setting of dimensions $n\in\{4,5\}$ and we define the mass of $\Delta_g+h_0$ in this case. In Section~\ref{sec:small:dim}, we then deal with the $C^1$-asymptotic expansion of the energy functional when $n\in\{4,5\}$. In Sections~\ref{sec:pf:1},~\ref{sec:pf:2},~\ref{sec:pf:3} and~\ref{sec:pf:4}, we complete the proofs of Theorems~\ref{th:main},~\ref{th:main:u0},~\ref{th:2} and \ref{th:3}, respectively. Finally, in Section~\ref{sec:nec}, we deal with the necessity of condition \eqref{phi0th2} on the gradient 

\section{Existence results for $h_0$ and $u_0$}\label{Sec2}

This short section deals with two results which provide conditions for the existence of functions $h_0$ and $u_0$ satisfying the assumptions in Theorem~\ref{th:3} with prescribed $\varphi_{h_0,u_0}$ and $\xi_0$. The first result is a straightforward consequence of classical works on the Yamabe equation:

\begin{theorem}(Aubin~\cite{Aubin}, Schoen~\cite{Schoen1}, Trudinger~\cite{Trudinger})
Assume that $n\ge3$. Then there exists $\eps_0\ge0$ depending only on $n$ and $(M,g)$ such that $\eps_0>0$ if $(M,g)$ is not conformally diffeomorphic to the standard sphere, $\eps_0=0$ otherwise, and for every $\varphi_0\in C^1(M)$ such that 
$$\varphi_0\le\eps_0\text{ and }\lambda_1(\Delta_g+h_0)>0,\text{ where }h_0:=\varphi_0+c_n\Scal_g,$$
there exists a solution $u_0\in C^2(M)$ of the equation \eqref{eq}. In particular, if $n\ne6$ and $\varphi_0(\xi_0)=|\nabla\varphi_0(\xi_0)|=0$ at some point $\xi_0\in M$, then $h_0$ satisfies \eqref{phi0th3}.
\end{theorem}

It remains to deal with the case where $n=6$. In this case, we obtain the following:

\begin{proposition}\label{Pr6}
Assume that $n=6$. Let $\varphi_0\in C^p(M)$, $1\le p\le\infty$, be such that 
\begin{equation}\label{condition6bis}
\lambda_1(\Delta_g+\varphi_0+c_n\Scal_g)<0.
\end{equation}
Then there exists $h_0\in C^p(M)$ such that the equation \eqref{eq} admits a solution $u_0\in C^2(M)$ satisfying $h_0-c_n\Scal_g-2u_0\equiv\varphi_0$. In particular, if $\varphi_0(\xi_0)=|\nabla\varphi_0(\xi_0)|=0$ at some point $\xi_0\in M$, then $(h_0,u_0)$ satisfy \eqref{phi0th3}.
\end{proposition}

\proof[Proof of Proposition~\ref{Pr6}]
Note that $\crit-1=2$ when $n=6$. In this case, we can rewrite the equation \eqref{eq} as
\begin{equation}\label{eq:u0bis}
\Delta_g u+(h_0-2u) u=-u^2,\ u>0\hbox{ in }M.
\end{equation}
Using \eqref{condition6bis} together with a standard variational method, we obtain that there exists a solution $u_0\in C^{p+1}(M)\subset C^2(M)$ of the equation \eqref{eq:u0bis} with $h_0:=\varphi_0+c_n\Scal_g+2u_0\in C^p(M)$. This ends the proof of Proposition~\ref{Pr6}.
\endproof

\section{Notations and general setting}\label{Sec3}

We follow the notations and definitions of Robert--V\'etois~\cite{RVProc}. 

\subsection{Euclidean setting} 

We define 
\begin{equation}\label{def:U}
U_{1,0}(x):=\Bigg(\frac{\sqrt{n(n-2)}}{1+|x|^2}\Bigg)^{\frac{n-2}{2}}\hbox{ for all }x\in\rn,
\end{equation}
so that $U_{1,0}$ is a positive solution to the equation 
$$\Delta_{\Eucl} U=U^{\crit-1}\hbox{ in }\rn,$$
where $\Eucl$ stands for the Euclidean metric. For every $\delta>0$ and $\xi\in \rn$, we define
\begin{equation}\label{def:U:delta}
U_{\delta,\xi}(x):=\delta^{-\frac{n-2}{2}}U\left(\delta^{-1}(x-\xi)\right)=\Bigg(\frac{\sqrt{n(n-2)}\delta}{\delta^2+|x-\xi|^2}\Bigg)^{\frac{n-2}{2}}\hbox{ for all }x\in\rn.
\end{equation}
We define
\begin{equation}\label{def:Z}
Z_0:=\left(\partial_\delta U_{\delta,\xi}\right)_{|(1,0)} \hbox{ and }Z_i:=\left(\partial_{\xi_i} U_{\delta,\xi}\right)_{|(1,0)}\hbox{ for all }i=1,\dotsc,n.
\end{equation}
As one checks, 
\begin{equation}\label{def:Z:0}
Z_0=-\frac{n-2}{2}U-(x,\nabla U)=\sqrt{n(n-2)}^{\frac{n-2}{2}}\frac{n-2}{2}\frac{|x|^2-1}{(1+|x|^2)^{\frac{n}{2}}}
\end{equation} and 
\begin{equation}\label{def:Z:i}
Z_i=-\partial_i U=\sqrt{n(n-2)}^{\frac{n-2}{2}}(n-2)\frac{x_i}{(1+|x|^2)^{\frac{n}{2}}}\hbox{ for all }i=1,\dotsc,n.
\end{equation} 
We denote $p=(p_0,p_1,\dotsc,p_n):=(\delta,\xi)\in (0,\infty)\times \rn$. Straightforward computations yield
\begin{equation}\label{der:U}
\partial_{p_i} U_{\delta,\xi}=\delta^{-1}(Z_i)_{\delta,\xi}:=\delta^{-1}\delta^{-\frac{n-2}{2}}Z_i\left(\delta^{-1}(x-\xi)\right)\hbox{ for all }i=0,\dotsc,n,
\end{equation}
\begin{equation}\label{est:der:U:delta}
\partial_\delta U_{\delta,\xi}=\sqrt{n(n-2)}^{\frac{n-2}{2}}\frac{n-2}{2}\delta^{\frac{n-2}{2}-1}\frac{|x-\xi|^2-\delta^2}{(\delta^2+|x-\xi|^2)^{n/2}}
\end{equation}
and
\begin{equation}\label{est:der:U:i}
\partial_{\xi_i}U_{\delta,\xi}=\sqrt{n(n-2)}^{\frac{n-2}{2}}(n-2)\delta^{\frac{n-2}{2}}\frac{(x-\xi)_i}{(\delta^2+|x-\xi|^2)^{n/2}}\hbox{ for all }i=1,\dotsc,n.
\end{equation} 
It follows from Rey~\cite{Rey} (see also Bianchi--Egnell~\cite{BE}) that 
$$\{\phi\in D_1^2(\rn):\ \Delta_{\Eucl} \phi=(\crit-1)U^{\crit-2}\phi\hbox{ in }\rn\}=\Span\{Z_i\}_{i=0,\dotsc,n}.$$

\subsection{Riemannian setting}\label{Riem} 

We fix $N>n-2$ to be chosen large later. It follows from Lee--Parker~\cite{LP} that there exists a function $\Lambda\in C^\infty(M\times M)$ such that, defining $\Lambda_\xi:=\Lambda(\xi,\cdot)$, we have
\begin{equation}\label{Lambda}
\Lambda_\xi>0,\ \Lambda_\xi(\xi)=1\hbox{ and }\nabla\Lambda_\xi(\xi)=0\hbox{ for all }\xi\in M
\end{equation}
and, defining the metric $g_\xi:=\Lambda_\xi^{\crit-2}g$ conformal to $g$, we have
\begin{equation}\label{prop:ge}
\Scal_{g_\xi}(\xi)=0,\ \nabla \Scal_{g_\xi}(\xi)=0,\ \Delta_g \Scal_{g_\xi}(\xi)=\frac{1}{6}|\Weyl_g(\xi)|_g^2
\end{equation}
and 
\begin{equation}\label{eq:elt:vol}
dv_{g_\xi}(x)=(1+\bigO(|x|^N))\,dx\hbox{ via the chart }\exp_{\xi}^{g_\xi} \hbox{ around }0,
\end{equation}
where $dx$ is the Euclidean volume element, $dv_{g_\xi}$ is the Riemannian volume element of $(M,g_\xi)$ and $\exp_{\xi}^{g_\xi}$ is the exponential chart at $\xi$ with respect to the metric $g_\xi$. The compactness of $M$ yields the existence of $r_0>0$ such that the injectivity radius of the metric $g_\xi$ satisfies $i_{g_\xi}(M)\geq3r_0$ for all $\xi\in M$. We let $\chi\in C^\infty(\rr)$ be such that $\chi(t)=1$ for all $t\leq r_0$, $\chi(t)=0$ for all $t\geq2r_0$ and $0\leq \chi\leq 1$. For every $\delta>0$ and $\xi\in M$, we then define the bubble as
\begin{align}\label{def:peak}
U_{\delta,\xi}(x):&=\chi(d_{g_\xi}(x,\xi))\Lambda_\xi(x)\delta^{-\frac{n-2}{2}} U_{1,0}(\delta^{-1}(\exp_{\xi}^{g_\xi})^{-1}(x))\\
&=\chi(d_{g_\xi}(x,\xi))\Lambda_\xi(x)\Bigg(\frac{\delta\sqrt{n(n-2)}}{\delta^2+d_{g_\xi}(x,\xi)^2}\Bigg)^{\frac{n-2}{2}},\nonumber
\end{align}
where $d_{g_\xi}(x,\xi)$ is the geodesic distance between $x$ and $\xi$ with respect to the metric $g_\xi$. Since there will never be ambiguity, to avoid unnecessary heavy notations, we will keep the notations $U_{\delta,\xi}$ as \eqref{def:U:delta} when $p=(\delta,\xi)\in (0,\infty)\times\rn$, and as \eqref{def:peak} when $p=(\delta,\xi)\in (0,\infty)\times M$. Finally, for every $p=(\delta,\xi)\in (0,\infty\times M$, we define
$$K_{\delta,\xi}:=\Span\{(Z_i)_{\delta,\xi}\}_{i=0,\dotsc,n},$$
where 
$$(Z_i)_{\delta,\xi}(x):=\chi(d_{g_\xi}(x,\xi))\Lambda_\xi(x)\delta^{-\frac{n-2}{2}} Z_i(\delta^{-1}(\exp_{\xi}^{g_\xi})^{-1}(x))$$
for all $x\in M$ and $i=0,\dotsc,n$.

\subsection{General reduction theorem} 

For every $1\le q\le\infty$, we let $\left\|\cdot\right\|_q$ be the usual norm of $L^q(M)$. We let $H_1^2(M)$ be the completion of $C^\infty(M)$ for the norm 
$$\Vert u\Vert_{H^2_1}:=\Vert\nabla u\Vert_2+\Vert u\Vert_2.$$
For every $h\in C^0(M)$, we define
$$J_h(u):=\frac{1}{2}\int_M\left(|\nabla u|_g^2+h u^2\right)dv_g-\frac{1}{\crit}\int_Mu_+^{\crit} dv_g\hbox{ for all }u\in H_1^2(M),$$
where $u_+:=\max(u,0)$. The space $H_1^2(M)$ is endowed with the bilinear form $\langle\cdot,\cdot\rangle_h$, where 
$$\langle u,v\rangle_h:=\int_M(\nabla u\nabla v+huv)\,dv_g\hbox{ for all }u,v\in H_1^2(M).$$ 
If $\Delta_g+h_0$ is coercive and $\Vert h-h_0\Vert_\infty$ is small enough, then $\langle\cdot,\cdot\rangle_h$ is positive definite and $(H_1^2(M),\langle\cdot,\rangle_{h})$ is a Hilbert space. We then have that $J_h\in C^1(H_1^2(M))$ and
$$J_h'(u)[\phi]=\int_M(\nabla u\nabla\phi+hu\phi)\,dv_g-\int_Mu_+^{\crit-1}\phi\, dv_g=\langle u,\phi\rangle_h-\int_Mu_+^{\crit-1}\phi\, dv_g$$
for all $u,\phi\in H_1^2(M)$. We let $(\delta,\xi)\to B_{h,\delta,\xi}=B_h(\delta,\xi)$ be a function in $C^1((0,\infty)\times M,H_1^2(M))$ such that for every $\delta>0$, there exists $\eps(\delta)>0$ independent of $h$ and $\xi$ such that 
\begin{equation}\label{hyp:B}
\Vert B_{h,\delta,\xi}\Vert_{H_1^2}+\delta\Vert \partial_p B_{h,\delta,\xi}\Vert_{H_1^2}<\epsilon(\delta)\hbox{ for all }p=(\delta,\xi)\in (0,\infty)\times M
\end{equation}
and $\eps(\delta)\to0$ as $\delta\to0$. The function $B_{h,\delta,\xi}$ will be fixed later. We also let $\tu\in C^2(M)$. We define
$$W_{h,\tu,\delta,\xi}:=\tu+U_{\delta,\xi}+B_{h,\delta,\xi}.$$
We fix a point $\xi_0\in M$ and a function $h_0\in C^0(M)$ such that $\Delta_g+h_0$ is coercive. We let $u_0\in C^2(M)$ be a solution of the equation
$$\Delta_gu_0+h_0u_0=u_0^{\crit-1},\ u_0\ge0\hbox{ in }M.$$
It follows from the strong maximum principle that either $u_0\equiv 0$ or $u_0>0$. We assume that $u_0$ is {\it nondegenerate}, that is, for every $\phi\in H_1^2(M)$, 
$$\Delta_g\phi+h_0\phi=(\crit-1)u_0^{\crit-2}\phi\ \Longleftrightarrow\  \phi\equiv 0.$$
It then follows from Robert--V\'etois~\cite{RVProc} that there exist $\epsilon_0>0$, $U_0\subset M$ a small open neighborhood of $\xi_0$ and $\Phi_{h,\tu}\in C^1((0,\epsilon_0)\times U_0,H_1^2(M))$ such that, when $\Vert h-h_0\Vert_{\infty}<\eps_0$ and $\Vert \tu-u_0\Vert_{C^2}<\eps_0$, we have
\begin{equation}\label{eq:ortho}
\Pi_{K_{\delta,\xi}^\perp}(W_{h,\tu,\delta,\xi}+\Phi_{h,\tu,\delta,\xi}-(\Delta_g+h)^{-1}((W_{h,\tu,\delta,\xi}+\Phi_{h,\tu,\delta,\xi})_+^{\crit-1}))=0
\end{equation}
and
\begin{equation}\label{control:R:1}
\Vert\Phi_{h,\tu,\delta,\xi}\Vert_{H_1^2}\leq C\Vert W_{h,\tu,\delta,\xi}-(\Delta_g+h)^{-1}((W_{h,\tu,\delta,\xi})_+^{\crit-1})\Vert_{H_1^2}\leq C\left\Vert R_{\delta_,\xi} \right\Vert_{\frac{2n}{n+2}}
\end{equation}
for all $(\delta,\xi)\in (0,\eps_0)\times U_0$, where $C>0$ does not depend on $(h,\tu,\delta,\xi)$, $\Phi_{h,\tu,\delta,\xi}:=\Phi_{h,\tu}(\delta,\xi)$,  $\Pi_{K_{\delta,\xi}^\perp}$ is the orthogonal projection of $H^2_1(M)$ onto $K_{\delta,\xi}^\perp$ (here, the orthogonality is taken with respect to $\langle\cdot,\cdot\rangle_h$) and
\begin{equation}\label{def:R}
R_{\delta_,\xi}:=(\Delta_g+h)W_{h,\tu,\delta,\xi}-(W_{h,\tu,\delta,\xi})_+^{\crit-1}.
\end{equation}
Furthermore, for every $(\delta_0,\xi_0)\in (0,\eps_0)\times U_0$, we have
\begin{multline}\label{cns}
J_h^\prime (W_{h,\tu,\delta_0,\xi_0}+\Phi_{h,\tu,\delta_0,\xi_0})=0\\
\Longleftrightarrow\  (\delta_0,\xi_0)\hbox{ is a critical point of }(\delta,\xi)\mapsto J_h(W_{h,\tu,\delta,\xi}+\Phi_{h,\tu,\delta,\xi}).
\end{multline}
It follows from Robert--V\'etois~\cite{RVProc} that
\begin{equation}
J_{h}(W_{h,\tu,\delta,\xi}+\Phi_{h,\tu,\delta,\xi})=J_h\left(W_{h,\tu,\delta,\xi}\right)+\bigO(\Vert\Phi_{h,\tu,\delta,\xi}\Vert_{H_1^2}^2)\label{est:fin:1}
\end{equation}
uniformly with respect to $(\delta,\xi)\in (0,\eps_0)\times U_0$ and $(h,\tu)$ such that $\Vert h-h_0\Vert_\infty<\eps_0$ and $\Vert\tu-u_0\Vert_{C^2}<\eps_0$. 

\smallskip\noindent
{\bf Conventions:}
\begin{itemize}\samepage
\item  To avoid unnecessarily heavy notations, we will often drop the indices $(h,,\tu,\delta,\xi)$, so that $U:=U_{\delta,\xi}$, $B:=B_{h,\delta,\xi}$, $W:=W_{h,\tu,\delta,\xi}$, $\Phi:=\Phi_{h,\tu,\delta,\xi}$, etc. The differentiation with respect to the variable $(\delta,\xi)$ will always be denoted by $\partial_p$, and the differentiation with respect to $x\in M$ (or $\rn$) by $\partial_x$. For example,
$$\partial_{x_i}\partial_{p_j} W=\left\{\begin{aligned}&\frac{\partial^2 W_{h,\tu,\delta,\xi}(x)}{\partial x_i\partial\delta}&&\hbox{if }j=0\\
&\frac{\partial^2 W_{h,\tu,\delta,\xi}(x)}{\partial x_i\partial\xi_j}&&\hbox{if }j\ge1.
\end{aligned}\right.$$
\smallskip
\item For every $\xi\in U_0$, we identify the tangent space $T_\xi M$ with $\rn$. Indeed, assuming that the neighborhood $U_0$ is small enough, it follows from the Gram--Schmidt orthonormalization procedure that there exists an orthonormal frame with respect to the metric $g_\xi$, which is smooth with respect to the point $\xi$. Such a frame provides a smooth family of linear isometries $(\psi_\xi)_{\xi\in U_0}$, $\psi_\xi:\rn\to T_\xi M$, which allow to identify $T_\xi M$ with $\rn$. In particular, in this paper, the chart $\exp_\xi^{g_\xi}$ will denote the composition of the usual exponential chart with the isometry $\psi_\xi$. 
\smallskip
\item Throughout the paper, $C$ will denote a positive constant such that
\begin{itemize} 
\item $C$ depends on $n$, $(M,g)$, $\xi_0\in M$, the functions $h_0,u_0\in C^2(M)$ and a constant $A>0$ such that $\left\|h_0\right\|_{C^2}<A$ and $\lambda_1(\Delta_g+h_0)>1/A$. In the case where $u_0>0$, we also assume that $\left\|u_0\right\|_{C^2}<A$ and $u_0>1/A$. 
\item $C$ does not depend on $x\in M$ (or $x\in\rn$, depending on the context), $\xi$ in the neighborhood $U_0$, $\delta>0$ small and $h\in C^2(M)$ such that $\left\|h\right\|_{C^2}<A$ and $\lambda_1(\Delta_g+h)>1/A$. In the case where $u_0>0$, $C$ is also independent of $\tilde{u}_0\in C^2(M)$ such that $\left\|\tilde{u}_0\right\|_{C^2}<A$ and $\tilde{u}_0>1/A$. 
\end{itemize}
The value of $C$ might change from line to line, and even in the same line. 
\smallskip
\item For every $f,g\in\rr$, the notations $f=\bigO(g)$ and $f=\smallo(g)$ will stand for $\left|f\right|\le C\left|g\right|$ and $\left|f\right|\le C\epsilon(h,\delta,\xi)\left|g\right|$, respectively, where $\epsilon(h,\delta,\xi)\to0$ as $h\to h_0$ in $C^2(M)$, $\delta\to0$ and $\xi\to\xi_0$. 
\end{itemize}

\section{$C^1$-estimates for the energy functional}\label{Sec4}

For every $\delta>0$ and $\xi\in U_0$, we define
\begin{equation}\label{def:tildeU}
\tilde{U}_{\delta,\xi}(x):=\Bigg(\frac{\delta\sqrt{n(n-2)}}{\delta^2+d_{g_\xi}(x,\xi)^2}\Bigg)^{\frac{n-2}{2}}\hbox{ for all }x\in M.
\end{equation}
Our first result is the differentiable version of \eqref{est:fin:1}.

\begin{proposition}\label{prop:c1} 
In addition to the assumptions of Section~\ref{Sec3}, we assume that
\begin{equation}\label{hyp:B:p}
|B_{h,\delta,\xi}(x)|+\delta |\partial_pB_{h,\delta,\xi}(x)|\leq C(U_{\delta,\xi}(x)+\delta\tilde{U}_{\delta,\xi}(x)).
\end{equation}
We then have
\begin{multline}
\partial_pJ_h(W+\Phi)=\partial_pJ_h(W)+\bigO(\delta^{-1}\Vert \Phi\Vert_{H_1^2}(\Vert R\Vert_{\frac{2n}{n+2}}+\delta\Vert\partial_p R\Vert_{\frac{2n}{n+2}}+\Vert\Phi\Vert_{H_1^2}))\\
+\bigO({\bf 1}_{n\geq 7}\delta^{-1}\Vert \Phi\Vert_{H_1^2}^{\crit-1}),\label{est:c1}
\end{multline}
where $R=R_{\delta,\xi}$ is as in \eqref{def:R}.
\end{proposition}

\proof[Proof of Proposition~\ref{prop:c1}]
It follows from \eqref{eq:ortho} that there exist real numbers $\lambda_j:=\lambda_j(\delta,\xi)$ such that
$$W+\Phi-(\Delta_g+h)^{-1}(W+\Phi)_+^{\crit-1}=\sum_{j=0}^n\lambda_j Z_j.$$
This can be written as
\begin{equation}\label{eq:Jh}
J_h'(W+\Phi)=\sum_{j=0}^n\lambda_j \langle Z_j,\cdot\rangle_h.
\end{equation}
We fix $i\in\{0,\dotsc,n\}$.  We obtain
\begin{align}\label{est:c1:1}
&\partial_{p_i}J_h(W+\Phi)= J_h'(W+\Phi)[\partial_{p_i}W+\partial_{p_i}\Phi]\\
&\qquad= J_h'(W)[\partial_{p_i}W]+\left(J_h'(W+\Phi)-J_h'(W)\right)[\partial_{p_i}W] +J_h'(W+\Phi)[\partial_{p_i}\Phi]\nonumber\\
&\qquad= J_h'(W)[\partial_{p_i}W]+\left(J_h'(W+\Phi)-J_h'(W)\right)[\partial_{p_i}W]+\sum_{j=0}^n\lambda_j \langle Z_j,\partial_{p_i}\Phi\rangle_h\nonumber\\
&\qquad= \partial_{p_i}J_h(W)+\left(J_h'(W+\Phi)-J_h'(W)\right)[\partial_{p_i}W]-\sum_{j=0}^n\lambda_j \langle\partial_{p_i} Z_j,\Phi\rangle_h,\nonumber
\end{align}
where, for the last line, we have used that $\langle (Z_i)_{\delta,\xi}, \Phi_{h,\tu,\delta,\xi}\rangle_h=0$ for all $(\delta,\xi)$ since $\Phi_{h,\tu,\delta,\xi}\in K_{\delta,\xi}^\perp$. We estimate separately the two last terms in the right-hand side of \eqref{est:c1:1}. As regards the first of these two term, we have 
\begin{align}\label{est:c1:2}
&\left(J_h'(W+\Phi)-J_h'(W)\right)[\partial_{p_i}W]\\
&\quad=\int_M (\nabla \Phi\nabla\partial_{p_i}W+h\Phi\partial_{p_i}W)-\int_M((W+\Phi)_+^{\crit-1}-W_+^{\crit-1})\partial_{p_i}W\, dv_g\nonumber\\
&\quad=\int_M \Phi ((\Delta_g+h)\partial_{p_i}W-(\crit-1)W_+^{\crit-1}\partial_{p_i}W)\,dv_g\nonumber\\
&\qquad-\int_M((W+\Phi)_+^{\crit-1}-W_+^{\crit-1}-(\crit-1)W_+^{\crit-1}\Phi)\partial_{p_i}W\, dv_g.\nonumber
\end{align}
With the definition \eqref{def:R}, H\"older's and Sobolev's inequalities, we obtain
\begin{multline}
\int_M \Phi ((\Delta_g+h)\partial_{p_i}W-(\crit-1)W_+^{\crit-1}\partial_{p_i}W)\,dv_g\\
=\int_M \Phi \partial_{p_i}R\, dv_g=\bigO(\Vert\Phi\Vert_{\crit}\Vert\partial_{p_i}R\Vert_{\frac{2n}{n+2}})=\bigO(\Vert\Phi\Vert_{H_1^2}\Vert\partial_{p_i}R\Vert_{\frac{2n}{n+2}}).\label{est:c1:3}
\end{multline}
In the sequel, we will need the following lemma:

\begin{lem}\label{lem:upper:U} 
We have
\begin{equation}\label{est:lem:1}
U_{\delta,\xi}(x)+\delta|\partial_pU_{\delta,\xi}(x)|\leq C\tilde{U}_{\delta,\xi}(x)
\end{equation}
for all $(\delta,\xi)\in (0,\eps_0)\times U_0$ and $x\in M$.
\end{lem}

\proof[Proof of Lemma~\ref{est:lem:1}] 
Most of the proof is easy computations. The only delicate point is to prove that $|\partial_{\xi}d_{g_\xi}(x,\xi)^2|\leq C d_{g_\xi}(x,\xi)$ for all $x\in M$ and $\xi\in U_0$. We define $F(x,\xi):=d_{g_\xi}(x,\xi)^2$ and $G(\xi,y):=\exp_{\xi}^{g_\xi}(y)$. Proving the desired inequality amounts to proving that $(\partial_{\xi}F(x,\xi))_{|\xi=x}=0$ for all $x\in M$. Note that $F(G(\xi,y),\xi)=|y|^2$ for small $y\in \rn$. Differentiating this equality with respect to $\xi$ yields a relation between $\partial_xF$ and $\partial_\xi F$, and the requested inequality follows.
\endproof

\proof[End of proof of Proposition~\ref{prop:c1}] 
Using Lemma~\ref{lem:upper:U}, the assumption \eqref{hyp:B:p} on $B_{h,\delta,\xi}$, and that $\partial_{p_i}\tu=0$, we obtain
\begin{multline*}
\left|\int_M((W+\Phi)_+^{\crit-1}-W_+^{\crit-1}-(\crit-1)W_+^{\crit-2}\Phi)\partial_{p_i}W\, dv_g\right| \\
\leq C\delta^{-1}\int_M|(W+\Phi)_+^{\crit-1}-W_+^{\crit-1}-(\crit-1)W_+^{\crit-2}\Phi|\tilde{U}\, dv_g.
\end{multline*}
We split the integral in two. First
\begin{multline*}
\int_{|W|\leq 2|\Phi|}|(W+\Phi)_+^{\crit-1}-W_+^{\crit-1}-(\crit-1)W_+^{\crit-2}\Phi|\tilde{U}\, dv_g\\
\leq C \int_M|\Phi|^{\crit-1}\tilde{U}\, dv_g\leq C\Vert \Phi\Vert_{\crit}^{\crit-1}\Vert\tilde{U}\Vert_{\crit}\leq C\Vert \Phi\Vert_{H_1^2}^{\crit-1}.
\end{multline*}
As regards the other part, looking carefully at the signs of the different terms, we obtain
\begin{align*}
&\int_{|\Phi|\leq |W|/2}|(W+\Phi)_+^{\crit-1}-W_+^{\crit-1}-(\crit-1)W_+^{\crit-2}\Phi|\tilde{U}\, dv_g\\
&\qquad=\int_{|\Phi|\leq |W|/2}|W|^{\crit-1}\left|\left(1+\frac{\Phi}{W}\right)^{\crit-1}-1-(\crit-1)\frac{\Phi}{W}\right|\tilde{U}\, dv_g\\
&\qquad\leq C \int_{|\Phi|\leq |W|/2}|W|^{\crit-1}\left(\frac{\Phi}{W}\right)^2 \tilde{U}\, dv_g=C\int_{|\Phi|\leq |W|/2}|W|^{\crit-3}\Phi^2\tilde{U}\, dv_g.
\end{align*}
In case $n\le6$, that is $\crit\geq 3$, we obtain
$$\int_{|\Phi|\leq |W|/2}|W|^{\crit-3}|\Phi|^2\tilde{U}\, dv_g\leq \int_{M}\tilde{U}^{\crit-2}|\Phi|^2 dv_g\leq C\Vert \tilde{U}\Vert_{\crit}^{\crit-2}\Vert\Phi\Vert_{\crit}^2\leq C\Vert \Phi\Vert_{H_1^2}^{2}.$$
In case $n\ge7$, that is $\crit< 3$, arguing as above, we obtain
$$\int_{|\Phi|\leq |W|/2}|W|^{\crit-3}\Phi^2 \tilde{U}\, dv_g\leq C\int_{M}|\Phi|^{\crit-1}\tilde{U}\, dv_g\leq C \Vert \Phi\Vert_{H_1^2}^{\crit-1}.$$
Plugging these estimates together yields
\begin{multline}
\left|\int_M((W+\Phi)_+^{\crit-1}-W_+^{\crit-1}-(\crit-1)W_+^{\crit-2}\Phi)\partial_{p_i}W\, dv_g\right|\\
\leq C\delta^{-1}(\Vert\Phi\Vert_{H_1^2}^2+{\bf 1}_{n\geq 7}\Vert\Phi\Vert_{H_1^2}^{\crit-1}).\label{est:c1:4}
\end{multline}
As regards the last term in the right-hand side of \eqref{est:c1:1}, arguing as in the proof of Lemma~\ref{lem:upper:U}, we obtain $\Vert \partial_{p_i}Z_j\Vert_{H_1^2}\leq C/\delta$ for all $i,j=0,\dotsc,n$. Therefore, we obtain
\begin{equation}\label{est:c1:6}
\left|\sum_{j=0}^n\lambda_j \langle\partial_{p_i} Z_j,\Phi\rangle_h\right|\leq C\delta^{-1}\Lambda\Vert\Phi\Vert_{H_1^2},\hbox{ where }\Lambda:=\sum_{j=0}^n|\lambda_j|.
\end{equation}
It follows from \eqref{eq:Jh} that
$$J_h'(W+\Phi)[Z_i]=\sum_{j=0}^n\lambda_j\langle Z_i,Z_j\rangle_h$$
for all $i=0,\dotsc,n$. Since $\langle Z_i,Z_j\rangle_h\to 0$ if $i\neq j$ and $\to1$ if $i=j$ as $\delta\to 0$ and uniformly with respect to $\xi\in U_0$, we obtain
$$\Lambda\leq C\sum_{i=0}^n\left|J_h'(W+\Phi)[Z_i]\right|.$$
For every $i=0,\dotsc,n$, using that $\langle \Phi,Z_i\rangle_h=0$ and $\Vert W\Vert_{\crit}+\Vert Z_i\Vert_{\crit}\leq C$, we obtain 
\begin{align*}
\left|J_h'(W+\Phi)[Z_i]\right|&\leq \left|J_h'(W)[Z_i]\right|+\left|\langle \Phi,Z_i\rangle_h-\int_M((W+\Phi)_+^{\crit-1}-W_+^{\crit-1})Z_i\, dv_g\right|\\
&\leq \left|\int_M R Z_i\, dv_g\right|+C\int_M (|W|^{\crit-2}|\Phi|+|\Phi|^{\crit-1})|Z_i|\, dv_g\\
&\leq C\Vert R\Vert_{\frac{2n}{n+2}}+C(\Vert\Phi\Vert_{\crit}+\Vert\Phi\Vert_{\crit}^{\crit-1})\leq C\Vert R\Vert_{\frac{2n}{n+2}}+C\Vert\Phi\Vert_{\crit}.
\end{align*}
Therefore,
\begin{equation}\label{est:c1:5}
\Lambda\leq C\Vert R\Vert_{\frac{2n}{n+2}}+C\Vert\Phi\Vert_{\crit}.
\end{equation}
Plugging \eqref{est:c1:2}, \eqref{est:c1:3}, \eqref{est:c1:4},  \eqref{est:c1:6} and \eqref{est:c1:5} into \eqref{est:c1:1} yields \eqref{est:c1}. This proves Proposition~\ref{prop:c1}.
\endproof

\section{Energy and remainder estimates: the case $n\geq 6$ and $u_0\equiv\tilde{u}_0\equiv 0$}\label{Sec5}

In this section, we consider the case where $n\ge6$ and $u_0\equiv\tilde{u}_0\equiv 0$. In this case, we set $B_{h,\delta,\xi}\equiv 0$. Then $W_{h,\tu,\delta,\xi}=W_{\delta,\xi}\equiv U_{\delta,\xi}$ and the assumptions of Proposition~\ref{prop:c1} are satisfied. We prove the following estimates for $R=R_{\delta,\xi}$:

\begin{proposition}\label{prop:est:R:1}
Assume that $n\geq 6$ and $u_0\equiv\tilde{u}_0\equiv 0$. Then
\begin{equation}\label{control:R:2}
\Vert R\Vert_{\frac{2n}{n+2}}+\delta\Vert \partial_p R\Vert_{\frac{2n}{n+2}}\leq C \left\{\begin{aligned}
&\delta^{2}+D_{h,\xi}\delta^2\left(\ln(1/\delta)\right)^{2/3}&&\hbox{if }n=6\\
&\delta^{\frac{n-2}{2}} +D_{h,\xi}\delta^2&&\hbox{if }7\le n\le9\\
&\delta^{4}\left(\ln(1/\delta)\right)^{3/5} +D_{h,\xi}\delta^2&&\hbox{if }n=10\\
&\delta^4+D_{h,\xi}\delta^2&&\hbox{if }n\ge11,
\end{aligned}\right.
\end{equation}
where 
\begin{equation}\label{Dhxi}
D_{h,\xi}:=\Vert h-h_0\Vert_\infty+d_g(\xi,\xi_0)^2.
\end{equation}
\end{proposition}

\proof[Proof of Proposition~\ref{prop:est:R:1}]
Define the conformal Laplacian as the operator $L_g:=\Delta_g+c_n \Scal_g$. For a metric $g'=w^{4/(n-2)}g$ conformal to $g$ ($w\in C^\infty(M)$ is positive), the conformal invariance law gives that 
\begin{equation}\label{conf:law}
L_{g'}\phi=w^{-(\crit-1)}L_g(w\phi)\hbox{ for all }\phi\in C^\infty(M).
\end{equation}
Therefore, we have 
\begin{align*}
R=(\Delta_g+h)U-U^{\crit-1}&= L_gU-U^{\crit-1}+\varphi_hU\\
 &= \Lambda_\xi^{\crit-1}(L_{g_\xi}(\Lambda_\xi^{-1}U)-(\Lambda_\xi^{-1}U)^{\crit-1}) +\varphi_hU\\
 &= \Lambda_\xi^{\crit-1}(\Delta_{g_\xi}(\Lambda_\xi^{-1}U)-(\Lambda_\xi^{-1}U)^{\crit-1})+\hat{h}_{\xi}U,
\end{align*}
where $\varphi_h$ is as in \eqref{defphi0th3} and
\begin{equation}\label{def:hat:h}
\hat{h}_{\xi}:=\varphi_h+ c_n \Lambda_{\xi}^{\crit-2}\Scal_{g_{\xi}}.
\end{equation}
Via the exponential chart, using the radial symmetry of $U_{\delta,0}:\rn\to \rr$, we obtain that around $0$,
\begin{equation}
\Delta_{g_\xi}(\Lambda_\xi^{-1}U)-(\Lambda_\xi^{-1}U)^{\crit-1}=\Delta_{\Eucl}U_{\delta,0}+\frac{\partial_r\sqrt{|g_{\xi}|}}{\sqrt{|g_{\xi}|}}\partial_r U_{\delta,0} -U_{\delta,0} ^{\crit-1}=\frac{\partial_r\sqrt{|g_{\xi}|}}{\sqrt{|g_{\xi}|}}\partial_r U_{\delta,0}.\label{eq:rest:U} 
\end{equation}
It then follows from \eqref{eq:elt:vol} that
\begin{equation}\label{eq:W:lap}
R(x)=\hat{h}_{\xi}(x)U(x) +\delta^{\frac{n-2}{2}}\Theta_{\delta,\xi}(x),\hbox{ where }|\Theta_{\delta,\xi}(x)|+|\partial_p\Theta_{\delta,\xi}(x)|\leq C
\end{equation}
for all $(\delta,\xi)\in (0,\infty)\times U_0$ and $x\in M$. Note that these estimates are a consequence of \eqref{eq:rest:U} when $x$ is close to $\xi$, and they are straightforward when $x$ is far from $\xi$. Using Lemma~\ref{lem:upper:U}, we then obtain
\begin{equation}\label{rest:1}
|R(x)|+\delta|\partial_\delta R(x)|\leq C\delta^{\frac{n-2}{2}}+C|\hat{h}_{\xi}(x)|\tilde{U}_{\delta,\xi}(x)
\end{equation}
and
\begin{equation}\label{rest:2}
\delta|\partial_{\xi}R(x)|\leq C\delta^{\frac{n-2}{2}}+C|\tilde{h}_{\xi}(x)|\tilde{U}_{\delta,\xi}(x)+C\delta|\partial_p\tilde{h}_{\xi}(x)|\tilde{U}_{\delta,\xi}(x).
\end{equation}
Since \eqref{phi0th3} and \eqref{prop:ge} hold, we have
\begin{equation}\label{rest:3}
|\hat{h}_{\xi}(x)|\leq CD_{h,\xi}+Cd_{g_\xi}(x,\xi)^2\hbox{ and }|\partial_{\xi}\hat{h}_{\xi}(x)|\leq C d_{g_\xi}(x,\xi).
\end{equation}
It is a straightforward computation that for every $\alpha>0$ and $p\geq 1$, we have
\begin{equation}\label{step:comput}
\Vert d_{g_\xi}(\cdot,\xi)^\alpha \tilde{U}_{\delta,\xi}\Vert_p\leq C\left\{\begin{aligned}
&\delta^{\frac{n-2}{2}}&&\hbox{if }n>(n-2-\alpha)p\\
&\delta^{\frac{n-2}{2}}\left(\ln(1/\delta)\right)^{1/p}&&\hbox{if }n=(n-2-\alpha)p\\
&\delta^{\frac{n}{p}+\alpha-\frac{n-2}{2}}&&\hbox{if }n<(n-2-\alpha)p.
\end{aligned}\right.
\end{equation}
Plugging together \eqref{rest:1}, \eqref{rest:2}, \eqref{rest:3} and \eqref{step:comput}, long but painless computations yield \eqref{control:R:2}. This ends the proof of Proposition~\ref{prop:est:R:1}.
\endproof

Since $n\geq 6$, that is $\crit-1\leq 2$, we have $\Vert\Phi\Vert_{H_1^2}^2=\bigO(\Vert\Phi\Vert_{H_1^2}^{\crit-1})$. Plugging together \eqref{est:fin:1}, \eqref{control:R:1}, \eqref{est:c1} and \eqref{control:R:2}, we obtain
\begin{equation}\label{est:J:21}
J_h(W+\Phi)=J_h(W)+\bigO\left(\begin{aligned}
&\delta^{4}+D_{h,\xi}^2\delta^{4}\left(\ln(1/\delta)\right)^{4/3}&&\hbox{if }n=6\\
&\delta^{n-2} +D_{h,\xi}^2\delta^4&&\hbox{if }7\le n\le9\\
&\delta^8\left(\ln(1/\delta)\right)^{6/5} +D_{h,\xi}^2\delta^4&&\hbox{if }n=10\\
&\delta^8+D_{h,\xi}^2\delta^4&&\hbox{if }n\ge11
\end{aligned}\right)
\end{equation}
and 
\begin{multline}\label{est:der:J:1}
\partial_{p_i}J_h(W+\Phi)=\partial_{p_i}J_h(W)\\
+\bigO\left(\delta^{-1}\right)\left\{\begin{aligned}
&\delta^{4}+D_{h,\xi}^2\delta^{4}\left(\ln(1/\delta)\right)^{4/3}&&\hbox{if }n=6\\
&(\delta^{\frac{n-2}{2}} +D_{h,\xi}\delta^2)^{\crit-1}&&\hbox{if }7\le n\le9\\
&(\delta^{4}\left(\ln(1/\delta)\right)^{3/5} +D_{h,\xi}\delta^2)^{\crit-1}&&\hbox{if }n=10\\
&(\delta^4+D_{h,\xi}\delta^2)^{\crit-1}&&\hbox{if }n\ge11
\end{aligned}\right.
\end{multline}
for all $i=0,\dotsc,n$. We now estimate $J_h(W+\Phi)$:

\begin{proposition}\label{prop:est:J:W:1} 
Assume that $n\geq 6$ and $u_0\equiv\tilde{u}_0\equiv 0$. Then
\begin{multline}\label{est:J:W:2}
J_h(W+\Phi)= \frac{1}{n}\int_{\rn}U_{1,0}^{\crit}\,dx+\frac{1}{2}\varphi_h(\xi)\delta^2\int_{\rn}U_{1,0}^2\, dx\\
-\frac{1}{4n}
\left\{\begin{aligned}
&24^2\omega_5K_{h_0}(\xi_0)\delta^4\ln(1/\delta)+\bigO(\delta^4(\smallo(\ln(1/\delta)+D_{h,\xi}^2(\ln(1/\delta))^{4/3}))&&\hbox{if }n=6\\
&K_{h_0}(\xi_0)\delta^4\int_{\rn}|x|^2U_{1,0}^2\, dx+\smallo(\delta^{4})&&\hbox{if }n\geq 7
\end{aligned}\right.
\end{multline}
as $\delta\to0$, $\xi\to\xi_0$ and $h\to h_0$ in $C^2(M)$, where $K_{h_0}(\xi_0)$ is as in \eqref{Kh}.
\end{proposition}

\proof[Proof of Proposition~\ref{prop:est:J:W:1}]
Integrating by parts, we obtain 
\begin{align}\label{est:J:1:0}
J_{h}(U)&= \frac{1}{2}\int_M \left[(\Delta_g+h)U\right]U\, dv_g-\frac{1}{\crit}\int_M U^{\crit} dv_g\\
&= \frac{1}{2}\int_M [(\Delta_g+h)U-U^{\crit-1}]U\, dv_g+\frac{1}{n}\int_M U^{\crit} dv_g.\nonumber
\end{align}
It follows from \eqref{eq:W:lap} that
\begin{equation}\label{est:J:2}
\int_M (\Delta_g U+h U-U^{\crit-1})U\, dv_g=\int_M \hat{h}_\xi U^2 dv_g+\bigO(\delta^{n-2}).
\end{equation}
Using the volume estimate \eqref{eq:elt:vol}, we obtain
\begin{align}\label{eq:3}
\int_M U^{\crit} dv_g=\int_M (\Lambda_{\xi}^{-1}U)^{\crit} dv_{g_\xi}&=\int_{B_{r_0}(0)}U_{\delta,0}^{\crit}(1+\bigO(|x|^N)\,dx+ \bigO(\delta^n)\\
&=\int_{\rn}U_{1,0}^{\crit}\,dx+\bigO(\delta^n).\nonumber
\end{align}
Plugging \eqref{est:J:2} and \eqref{eq:3} into \eqref{est:J:1:0}, we obtain
$$J_{h}(U)= \frac{1}{2}\int_M \hat{h}_\xi U^2 dv_g+\frac{1}{n}\int_{\rn}U_{1,0}^{\crit}\,dx+\bigO(\delta^{n-2}).$$
With the change of metric, the definition of the bubble \eqref{def:peak} and the property of the volume \eqref{eq:elt:vol}, we obtain
\begin{equation}\label{inthU2}
\int_M \hat{h}_\xi U^2 dv_g= \int_{B_{r_0}(\xi)}\hat{h}_\xi U^2 dv_g+\bigO(\delta^{n-2})=\int_{B_{r_0}(0)} A_{h,\xi} U_{\delta,0}^2\, dx+\bigO(\delta^{n-2}),
\end{equation}
where 
\begin{equation}\label{def:A}
A_{h,\xi}(x):=(\hat{h}_{\xi}\Lambda_{\xi}^{2-\crit})(\exp_{\xi}^{g_{\xi}}(x)).
\end{equation}
Using the radial symmetry of $U_{\delta,0}$ and since $h_0\in C^2(M)$, we obtain  
\begin{align}\label{est:A:U:L2}
&\int_{B_{r_0}(0)} A_{h,\xi} U_{\delta,0}^2\, dx= \int_{B_{r_0}(0)}(A_{h,\xi}(0)+\partial_{x_\alpha} A_{h,\xi}(0)x^\alpha \\
&\qquad\quad+\frac{1}{2}\partial_{x_\alpha}\partial_{x_\beta}A_{h,\xi}(0)x^\alpha x^\beta+\bigO(\left\|h-h_0\right\|_{C^2}|x|^2+\eps_{h_0,\xi}(x)|x|^2))U_{\delta,0}^2\, dx\nonumber\\
&\qquad=A_{h,\xi}(0)\int_{B_{r_0}(0)}U_{\delta,0}^2\, dx-\frac{1}{2n}\Delta_{\Eucl}A_{h,\xi}(0)\int_{B_{r_0}(0)} |x|^2 U_{\delta,0}^2\, dx\nonumber\\
&\qquad\quad+\bigO\left(\int_{B_{r_0}(0)}(\left\|h-h_0\right\|_{C^2}+\eps_{h_0,\xi}(x))|x|^2U_{\delta,0}^2\, dx \right)+\bigO(\delta^{n-2}),\nonumber
\end{align}
where $\eps_{h_0,\xi}(x)\to0$ as $x\to0$, uniformly in $\xi\in U_0$. With a change of variable and Lebesgue convergence theorem,  we obtain
\begin{equation}
\int_{B_{r_0}(0)} U_{\delta,0}^2\, dx= \delta^2\int_{\rn}U_{1,0}^2\, dx+\bigO(\delta^{n-2}),\label{est:U:L2}
\end{equation}
\begin{equation}
\int_{B_{r_0}(0)} |x|^2 U_{\delta,0}^2\, dx=\left\{\begin{aligned}
&24^2\omega_5\delta^4\ln(1/\delta)+\bigO(\delta^4)&&\hbox{if }n=6\\
&\delta^4\int_{\rn}|x|^2U_{1,0}^2\, dx+\bigO(\delta^5)&&\hbox{if }n\geq 7,
\end{aligned}\right.\label{est:U:L2:2}
\end{equation}
and
\begin{equation}
\int_{B_{r_0}(0)}\eps_{h_0,\xi}(x)|x|^2 U_{\delta,0}^2\, dx=\smallo\left(\begin{aligned}
&\delta^4\ln(1/\delta)&&\hbox{if }n=6\\
&\delta^4&&\hbox{if }n\geq 7
\end{aligned}\right).\label{est:U:L2:eps}
\end{equation}
Furthermore, we have $A_{h,\xi}(0)=\varphi_h(\xi)$ and
\begin{align}\label{est:lap}
\Delta_{\Eucl}A_{h,\xi}(0)&= \Delta_{g_{\xi}}(\hat{h}_\xi\Lambda_{\xi}^{2-\crit})(\xi)= L_{g_{\xi}}(\varphi_h\Lambda_{\xi}^{2-\crit})(\xi)+ c_n\Delta_{g_{\xi}}\Scal_{g_{\xi}}(\xi)\\
&= L_g(\varphi_h\Lambda_{\xi}^{3-\crit})(\xi)+\frac{c_n}{6}|\Weyl_g(\xi)|^2_g\nonumber\\
&= L_g(\varphi_{h_0}\Lambda_{\xi}^{3-\crit})(\xi)+\frac{c_n}{6}|\Weyl_g(\xi)|^2_g+\bigO(\Vert h-h_0\Vert_{C^2})\nonumber\\
&= K_{h_0}(\xi_0)+\bigO(\epsilon_{h_0}(\xi)+\Vert h-h_0\Vert_{C^2}),\nonumber
\end{align}
where $\epsilon_{h_0}(\xi)\to0$ as $\xi\to\xi_0$. Therefore, plugging together these identities yields
\begin{multline}\label{est:J:W}
J_h(U)= \frac{1}{n}\int_{\rn}U_{1,0}^{\crit}\,dx+\frac{1}{2}\varphi_h(\xi)\delta^2\int_{\rn}U_{1,0}^2\, dx\\
-\frac{1}{4n}
\left\{\begin{aligned}
&24^2\omega_5K_{h_0}(\xi_0)\delta^4\ln(1/\delta)+\smallo\left(\delta^{4}\ln(1/\delta)\right)&&\hbox{if }n=6\\
&K_{h_0}(\xi_0)\delta^4\int_{\rn}|x|^2U_{1,0}^2\, dx+\smallo(\delta^{4})&&\hbox{if }n\geq 7.
\end{aligned}\right.
\end{multline}
Plugging together \eqref{est:J:21} and \eqref{est:J:W}, we obtain \eqref{est:J:W:2}. This ends the proof of Proposition~\ref{prop:est:J:W:1}.
\endproof

We now estimate the derivatives of $J_h(W+\Phi)$:

\begin{proposition}\label{prop:est:der:J:1}
Assume that $n\geq 6$ and $u_0\equiv\tilde{u}_0\equiv 0$. Then
\begin{multline}\label{der:delta:2}
\partial_\delta J_h(W+\Phi)= \varphi_h(\xi)\delta\int_{\rn}U_{1,0}^2 dx\\
-\left\{\begin{aligned}
&24^2\omega_5K_{h_0}(\xi_0)\delta^3\ln(1/\delta)+\smallo\left(\delta^3\ln(1/\delta)\right)+\bigO(D_{h,\xi}^2\delta^3(\ln(1/\delta))^{4/3})&&\hbox{if }n=6\\
&K_{h_0}(\xi_0)\delta^3\int_{\rn}|x|^2U_{1,0}^2\, dx+\smallo(\delta^3)+\bigO(D_{h,\xi}^{\crit-1}\delta^{\frac{n+6}{n-2}})&&\hbox{if }n\geq 7
\end{aligned}\right.
\end{multline}
and
\begin{multline}\label{der:xi:2}
\partial_{\xi_i} J_h(W+\Phi)=\frac{1}{2}\partial_{\xi_i}\varphi_{h}(\xi)\delta^2\int_{\rn}U_{1,0}^2\, dx\\
+\bigO\left(\begin{aligned}
&\smallo(\delta^3\ln(1/\delta))+\bigO(D_{h,\xi}^2\delta^3\left(\ln(1/\delta)\right)^{4/3})&&\hbox{if }n=6\\
&\smallo(\delta^3) +\bigO(D_{h,\xi}^{\crit-1}\delta^{\frac{n+6}{n-2}})&&\hbox{if }n\geq 7
\end{aligned}\right)
\end{multline}
for all $i=1,\dotsc,n$, as $\delta\to0$, $\xi\to\xi_0$ and $h\to h_0$ in $C^2(M)$.
\end{proposition}

\proof[Proof of Proposition~\ref{prop:est:der:J:1}]
We fix $i\in\{0,\dotsc,n\}$. Using \eqref{eq:W:lap} and \eqref{est:lem:1} and arguing as in \eqref{inthU2}, we obtain
\begin{align}
\partial_{p_i}J_h(U)&=J_h^\prime(U)[\partial_{p_i}U]=\int_M(\Delta_gU+hU-U^{\crit-1})\partial_{p_i}U\, dv_g\label{est:26}\\
&=\int_M R\partial_{p_i}U\, dv_g=\int_M \hat{h}_{\xi}U\partial_{p_i}U\, dv_g+\bigO\left(\delta^{\frac{n-2}{2}}\int_M |\partial_{p_i}U|\, dv_g\right)\nonumber\\
&=  \int_M \hat{h}_{\xi}U\partial_{p_i}U\, dv_g+\bigO\left(\delta^{-1}\delta^{n-2}\right)\nonumber\\
&= \int_{B_{r_0}(0)}A_{h,\xi}U_{\delta,\xi}(\Lambda_\xi^{-1}\partial_{p_i}U)(\exp_{\xi}^{g_\xi}(x))\,dx+\bigO\left(\delta^{-1}\delta^{n-2}\right)\nonumber
\end{align}
As in \eqref{est:A:U:L2}, we write 
\begin{multline}
A_{h,\xi}(x)=A_{h,\xi}(0)+\partial_{x_\alpha}A_{h,\xi}(0)x^\alpha+\frac{1}{2}\partial_{x_j}\partial_{x_k}A_{h,\xi}(0)x^jx^k\\
+\bigO(\epsilon_{h_0,\xi}(x)|x|^2+\Vert h-h_0\Vert_{C^2}|x|^2)\label{taylor:A}
\end{multline}
for all $x\in B_{r_0}(0)$, where $\eps_{h_0,\xi}(x)\to0$ as $x\to0$, uniformly in $\xi\in U_0$. With \eqref{est:lem:1}, \eqref{est:U:L2:2} and \eqref{est:U:L2:eps}, we obtain 
\begin{multline}\label{control:23}
\left|\int_{B_{r_0}(0)}(\epsilon_{h_0,\xi}(x)+\Vert h-h_0\Vert_{C^2})|x|^2U_{\delta,0}(\Lambda_\xi^{-1}\partial_{p_i}U)(\exp_{\xi}^{g_\xi}(x))\,dx\right|\\
\leq C\delta^{-1}\int_{B_{r_0}(0)}(\epsilon_{h_0,\xi}(x)+\Vert h-h_0\Vert_{C^2})|x|^2\tilde{U}_{\delta,0}^2\, dx=\smallo(\delta^{-1})\left\{\begin{aligned}
&\delta^4\ln(1/\delta) &&\hbox{if }n=6\\
&\delta^4&&\hbox{if }n\geq 7.
\end{aligned}\right.
\end{multline}
We write
\begin{align*}
&\int_{B_{r_0}(0)}A_{h,\xi}U_{\delta,0}(\Lambda_\xi^{-1}\partial_{p_i}U)(\exp_{\xi}^{g_\xi}(x))\,dx\\
&\qquad= \int_{B_{r_0}(0)}A_{h,\xi}U_{\delta,0}\partial_{p_i}(\Lambda_\xi^{-1}U)(\exp_{\xi}^{g_\xi}(x))\,dx\\
&\qquad\quad-\int_{B_{r_0}(0)}A_{h,\xi}U_{\delta,0}^2(\Lambda_\xi^{-1}\partial_{p_i}\Lambda_\xi^{-1})(\exp_{\xi}^{g_\xi}(x))\,dx.
\end{align*}
Since $\nabla\Lambda_\xi(\xi)=0$, we obtain
\begin{multline*}
\int_{B_{r_0}(0)}A_{h,\xi}U_{\delta,0}^2(\Lambda_\xi^{-1}\partial_{p_i}\Lambda_\xi^{-1})(\exp_{\xi}^{g_\xi}(x))\,dx\\
=\bigO\left(A_{\delta,\xi}(0)\int_{B_{r_0}(0)} |x|U_{\delta,0}^2\, dx\right)+\bigO\left(\int_{B_{r_0}(0)}|x|^2U_{\delta,0}^2\,dx\right).
\end{multline*}
With the definition \eqref{def:A} of $A_{h,\xi}$ and the assumption \eqref{phi0th3} on $h_0$, it follows that
\begin{multline*}
\int_{B_{r_0}(0)}A_{h,\xi}U_{\delta,0}^2(\Lambda_\xi^{-1}\partial_{p_i}\Lambda_\xi^{-1})(\exp_{\xi}^{g_\xi}(x))\,dx\\
=\bigO\left(\delta^{-1}\delta^4\left(D_{h,\xi}+\left\{\begin{aligned}
&\delta\ln(1/\delta) &&\hbox{if }n=6\\
&\delta &&\hbox{if }n\geq 7
\end{aligned}\right\}\right)\right).
\end{multline*}
This estimate, the Taylor expansion \eqref{taylor:A} and the estimate   \eqref{control:23} yield
\begin{multline}
\int_{B_{r_0}(0)}A_{h,\xi}\Lambda_\xi^{-1}U_{\delta,0}(\Lambda_\xi^{-1}\partial_{p_i}U)(\exp_{\xi}^{g_\xi}(x))\,dx\label{est:24}\\
= A_{h,\xi}(0)\int_{B_{r_0}(0)}U_{\delta,0}\partial_{p_i}(\Lambda_\xi^{-1}U)(\exp_{\xi}^{g_\xi}(x))\,dx\\
+\partial_{x_\alpha} A_{h,\xi}(0)\int_{B_{r_0}(0)}x^\alpha U_{\delta,0}\partial_{p_i}(\Lambda_\xi^{-1}U)(\exp_{\xi}^{g_\xi}(x))\,dx\\
+\frac{1}{2}\partial_{x_j}\partial_{x_k}A_{h,\xi}(0)\int_{B_{r_0}(0)}x^jx^kU_{\delta,0}\partial_{p_i}(\Lambda_\xi^{-1}U)(\exp_{\xi}^{g_\xi}(x))\,dx\\
+\smallo(\delta^{-1})\left\{\begin{aligned}
&\delta^4\ln(1/\delta)&&\hbox{if }n=6\\
&\delta^4&&\hbox{if }n\geq 7.
\end{aligned}\right.
\end{multline}
We first deal with the case $i=0$, that is $\partial_{p_i}=\partial_{p_0}=\partial_\delta$. For every homogeneous polynomial $Q$ on $\rn$, it follows from \eqref{def:U:delta} and \eqref{der:U} that
\begin{align*}
&\int_{B_{r_0}(0)}QU_{\delta,0}\partial_{\delta}(\Lambda_\xi^{-1}U)(\exp_{\xi}^{g_\xi}(x))\,dx\\
&\qquad=\int_{B_{r_0}(0)}Q\delta^{-1}\delta^{-\frac{n-2}{2}}U_{1,0}\left(x/\delta\right)\delta^{-\frac{n-2}{2}}Z_0\left(x/\delta\right)\,dx.
\end{align*}
The explicit expressions \eqref{def:U} and \eqref{def:Z} of $U$ and $Z_0$ and their radial symmetry then yield
$$\int_{B_{r_0}(0)}U_{\delta,0}\partial_{\delta}(\Lambda_\xi^{-1}U)(\exp_{\xi}^{g_\xi}(x))\,dx=\delta^{-1}\delta^2\int_{\rn}U_{1,0} Z_0\, dx+\bigO(\delta^{-1}\delta^{n-2})\hbox{ for }n\geq 6,$$
$$\int_{B_{r_0}(0)}x^jU_{\delta,0}\partial_{\delta}(\Lambda_\xi^{-1}U)(\exp_{\xi}^{g_\xi}(x))\,dx=0\hbox{ for }n\geq 6,$$
and
\begin{multline*}
\int_{B_{r_0}(0)}x^jx^kU_{\delta,0}\partial_{\delta}(\Lambda_\xi^{-1}U)(\exp_{\xi}^{g_\xi}(x))\,dx\\
=\frac{\epsilon_{jk}}{n}\delta^{-1}\delta^4\left\{\begin{aligned}
&c'_6\ln(1/\delta)+\bigO(\delta^{-1}\delta^{4})&&\hbox{if }n=6\\
&\int_{\rn}|x|^2U_{1,0} Z_0\, dx+\bigO(\delta^{-1}\delta^{n-2})&&\hbox{if }n\geq 7,\end{aligned}\right.
\end{multline*}
where $\epsilon_{jk}$ is the Kronecker symbol and $c'_6>0$ is a constant that will be discussed later. Putting these estimates in \eqref{est:26}, and \eqref{est:24}, we obtain
\begin{multline*}
\partial_\delta J_h(U)= A_{h,\xi}(0)\delta^{-1}\delta^2\int_{\rn}U_{1,0} Z_0\, dx\\
-\frac{1}{2n}\delta^{-1}\delta^4\left\{\begin{aligned}
&c'_6\Delta_{Eucl}A_{h,\xi}(0)\ln(1/\delta)+\smallo(\ln(1/\delta))&&\hbox{if }n=6\\
&\Delta_{Eucl}A_{h,\xi}(0)\int_{\rn}|x|^2U_{1,0} Z_0\, dx+\smallo(1)&&\hbox{if }n\geq 7.
\end{aligned}\right.
\end{multline*}
For every $\delta>0$, we have
$$\int_{\rn}U_{\delta,0}^2\, dx=\delta^2\int_{\rn}U_{1,0}^2\, dx\hbox{ for }n\geq 5$$
and
$$\int_{\rn}|x|^2U_{\delta,0}^2\, dx=\delta^4\int_{\rn}|x|^2U_{1,0}^2\, dx\hbox{ for }n\geq 7.$$
Differentiating these equalities with respect to $\delta$ at $\delta=1$, we obtain
$$\int_{\rn}U_{1,0} Z_0\, dx=\int_{\rn}U_{1,0}^2\hbox{ for }n\geq 5$$
and
$$\int_{\rn}|x|^2U_{1,0} Z_0\, dx=2\int_{\rn}|x|^2U_{1,0}^2\hbox{ for }n\geq 7.$$
Therefore, with the computation \eqref{est:lap} and the definition \eqref{Kh}, we obtain
\begin{multline}
\partial_\delta J_h(U)= \varphi_h(\xi)\delta^{-1}\delta^2\int_{\rn}U_{1,0}^2\, dx\label{der:xi:0}\\
-\frac{1}{n}\delta^{-1}\delta^4\left\{\begin{aligned}
&c'_6K_{h_0}(\xi_0)\ln(1/\delta)+\smallo(\ln(1/\delta))&&\hbox{if }n=6\\
&K_{h_0}(\xi_0)\int_{\rn}|x|^2U_{1,0}^2\, dx+\smallo(1)&&\hbox{if }n\geq 7.
\end{aligned}\right.
\end{multline}
Differentiating \eqref{est:J:W}, we obtain $c^{\prime}_6/2=24^2\omega_5$. Therefore, with \eqref{est:der:J:1}, we obtain \eqref{der:delta:2}. We now deal with the case where $i\geq 1$, that is $\partial_{p_i}=\partial_{\xi_i}$. We first claim that
\begin{equation}\label{eq:30}
[\partial_{\xi_i}(\Lambda_\xi^{-1}U_{\delta,\xi})](\xi, \exp_{\xi}^{g_\xi}(x))+[\partial_{x_i}(\Lambda_\xi^{-1}U_{\delta,\xi})](\xi, \exp_{\xi}^{g_\xi}(x))=\bigO\left(\frac{\delta^{\frac{n-2}{2}}|x|^3}{\left(\delta^2+|x|^2\right)^{n/2}}\right),
\end{equation}
where the differential for $\xi$ is taken via the exponential chart. Before proving this claim, let us remark that it is trivial in the Euclidean context. Indeed, for every $\xi,x\in\rn$ and $\delta>0$, with the notation \eqref{def:U:delta}, we have
$$\partial_{\xi_i}U_{\delta,\xi}(x)=\partial_{\xi_i}(\delta^{-{\frac{n-2}{2}}}U(\delta^{-1}(x-\xi)))=-\partial_{x_i}U_{\delta,\xi}(x).$$
We now prove the claim \eqref{eq:30}. We fix $\xi\in U_0$. We define the path $\xi(t):=\exp_{\xi}^{g_{\xi}}(t \vec{e}_i)$ for small $t\in\rr$, where $\vec{e}_i$ is the $i$-th vector in the canonical basis of $\rn$. With \eqref{def:tildeU}, we obtain
\begin{align}\label{eq:31}
[\partial_{x_i}(\Lambda_\xi^{-1}U_{\delta,\xi})](\xi, \exp_{\xi}^{g_\xi}(x))&=\frac{d}{dt}\tilde{U}_{\delta, \xi}(\exp_{\xi}^{g_{\xi}}(x+t\vec{e}_i))_{|t=0}\\
&=-\frac{n-2}{2}\frac{\delta^{\frac{n-2}{2}}}{(\delta^2+|x|^2)^{n/2}}\cdot2x_i\nonumber
\end{align}
and
\begin{align}\label{eq:31:bis}
[\partial_{\xi_i}(\Lambda_\xi^{-1}U_{\delta,\xi})](\xi, \exp_{\xi}^{g_\xi}(x))&=\frac{d}{dt}\tilde{U}_{\delta, \xi(t)}(\exp_{\xi}^{g_{\xi}}(x))_{|t=0}\\
&= -\frac{n-2}{2}\frac{\delta^{\frac{n-2}{2}}}{(\delta^2+|x|^2)^{n/2}}\cdot \frac{d}{dt} d_{g_{\xi(t)}}^2(\xi(t), \exp_{\xi}^{g_\xi}(x)).\nonumber
\end{align}
It follows from Esposito--Pistoia--V\'etois~\cite{EPV}*{Lemma~A.2} that
\begin{equation}\label{eq:d}
\frac{d}{dt} d_{g_{\xi(t)}}^2(\xi(t), \exp_{\xi}^{g_\xi}(x))+2x_i=\bigO(|x|^3)\hbox{ as }x\to0.
\end{equation}
Putting together all these estimates yields \eqref{eq:30}. This proves the claim. With the definition \eqref{def:U:delta}, we obtain
$$\int_{B_{r_0}(0)}U_{\delta,0}\frac{\delta^{\frac{n-2}{2}}|x|^3}{\left(\delta^2+|x|^2\right)^{n/2}} dx=\bigO(\delta^3)\hbox{ for }n\geq 6,$$
$$\int_{B_{r_0}(0)}|x|U_{\delta,0}\frac{\delta^{\frac{n-2}{2}}|x|^3}{\left(\delta^2+|x|^2\right)^{n/2}} dx=\bigO\left(\begin{aligned}
&\delta^4\ln(1/\delta)&&\hbox{if }n=6\\
&\delta^4&&\hbox{if }n\geq 7
\end{aligned}\right)$$
and
$$\int_{B_{r_0}(0)}|x|^2U_{\delta,0} \frac{\delta^{\frac{n-2}{2}}|x|^3}{\left(\delta^2+|x|^2\right)^{n/2}} dx=\bigO\left(\begin{aligned}
&\delta^4 &&\hbox{if }n=6\\
&\delta^5\ln(1/\delta) &&\hbox{if }n=7\\
&\delta^5&&\hbox{if }n\geq 8.
\end{aligned}\right).$$
Noting that $[\partial_{x_i}(\Lambda_\xi^{-1}U_{\delta,\xi})](\xi, \exp_{\xi}^{g_\xi}(x))=\partial_{x_i}U_{\delta,0}$, we obtain by symmetry that
$$\int_{B_{r_0}(0)}\Lambda_\xi^{-1}U_{\delta,0}\partial_{x_i}(\Lambda_\xi^{-1}U)(\exp_{\xi}^{g_\xi}(x))\,dx=\int_{B_{r_0}(0)}U_{\delta,0}\partial_{x_i}U_{\delta,0}\, dx=0$$
and similarly,
$$\int_{B_{r_0}(0)}x^jx^kU_{\delta,0}\partial_{x_i}(\Lambda_\xi^{-1}U)(\exp_{\xi}^{g_\xi}(x))\,dx=0.$$
Integrating by parts, straightforward estimates yield
\begin{align*}
&\int_{B_{r_0}(0)}x^jU_{\delta,0}\partial_{x_i}(\Lambda_\xi^{-1}U)(\exp_{\xi}^{g_\xi}(x))\,dx=\int_{B_{r_0}(0)}x^jU_{\delta,0}\partial_{x_i}U_{\delta,0}\, dx\\
&\qquad=\frac{1}{2}\int_{B_{r_0}(0)}x^j\partial_{x_i}(U_{\delta,0}^2)\,dx= -\frac{\epsilon_{ij}}{2}\int_{B_{r_0}(0)}U_{\delta,0}^2\, dx+\frac{1}{2}\int_{\partial B_{r_0}(0)}x^j\vec{\nu}_iU_{\delta,0}^2\, d\sigma\\
&\qquad= -\frac{\epsilon_{ij}}{2}\delta^2\int_{\rn}U_{1,0}^2\, dx+\bigO(\delta^{n-2})\hbox{ for }n\geq 6,
\end{align*}
where $\vec{\nu}:=(\vec{\nu}_1,\dotsc,\vec{\nu}_n)$ is the outward unit normal vector and $d\sigma$ is the volume element of $\partial B_{r_0}(0)$. Since $A_{h,\xi}(0)=\bigO(D_{h,\xi})$, plugging these estimates together with \eqref{est:26} and \eqref{est:24}, we obtain
\begin{equation}\label{der:xi:1}
\partial_{\xi_i} J_h(U)=\frac{1}{2}\partial_{\xi_i} \varphi_{h}(\xi)\delta^2\int_{\rn}U_{1,0}^2\, dx+\smallo\left(\begin{aligned}
&\delta^{3}\ln(1/\delta)&&\hbox{if }n=6\\
&\delta^{3}&&\hbox{if }n\geq 7
\end{aligned}\right).
\end{equation}
With \eqref{est:der:J:1}, we then obtain \eqref{der:xi:2}. This ends the proof of Proposition~\ref{prop:est:der:J:1}. 
\endproof
 
Theorem~\ref{th:main} for $n\geq 6$ will be proved in Section~\ref{sec:pf:1}.

\section{Energy and remainder estimates: the case $n\geq 7$ and $u_0,\tilde{u}_0>0$}\label{Sec6}

In this section, we assume that $u_0,\tilde{u}_0>0$ and $n\geq 7$, that is $\crit-1<2$. As in the previous case, we set $B_{h,\delta,\xi}\equiv 0$, so that $W_{h,\tu,\delta,\xi}=W_{\tu,\delta,\xi}\equiv\tu+U_{\delta,\xi}$ and the assumptions of Proposition~\ref{prop:c1} are satisfied. We prove the following estimates for $R=R_{\delta,\xi}$:

\begin{proposition}\label{prop:est:R:2}
Assume that $n\geq 7$ and $u_0,\tilde{u}_0>0$. Then
\begin{equation}\label{est:R:2}
\Vert R\Vert_{\frac{2n}{n+2}}\leq C\Vert \Delta_g \tu+h\tu-\tu^{\crit-1}\Vert_{\infty}+C(D_{h,\xi}+\delta^2+\delta^{\frac{n-6}{2}})\delta^2\text{ and }\Vert \partial_p R\Vert_{\frac{2n}{n+2}}\leq C\delta,
\end{equation}
where $D_{h,\xi}$ is as in \eqref{Dhxi}.
\end{proposition}

\proof[Proof of Proposition~\ref{prop:est:R:2}]
We have
\begin{equation}
R=(\Delta_g \tu+h\tu-\tu^{\crit-1})+R^0-((\tu+U)^{\crit-1}-\tu^{\crit-1}-U^{\crit-1}),\label{est:R:u}
\end{equation}
where 
$$R^0:=\Delta_gU+hU-U^{\crit-1}.$$
Concerning the derivatives, given $i\in\{0,\dotsc,n\}$, we have
\begin{align}
\partial_{p_i}R&= \Delta_g\partial_{p_i}U+h\partial_{p_i}U-(\crit-1)(\tu+U)^{\crit-2}\partial_{p_i}U\label{est:der:R:u}\\
&=\partial_{p_i}R^0-(\crit-1)((\tu+U)^{\crit-2}-U^{\crit-2})\partial_{p_i}U.\nonumber
\end{align}
A straightforward estimate yields
$$|(\tu+U)^{\crit-1}-\tu^{\crit-1}-U^{\crit-1}|\leq C{\bf 1}_{U\leq \tu}\tu^{\crit-2}U+C{\bf 1}_{\tu\leq U}\tu U^{\crit-2}.$$
With the expression \eqref{def:peak}, we obtain
$$\{U(x)\leq \tu(x) \, \Rightarrow\, d_{g_\xi}(x,\xi)\geq c_1\sqrt{\delta}\}\hbox{ and }\{U(x)\geq \tu(x) \, \Rightarrow\, d_{g_\xi}(x,\xi)\leq c_2\sqrt{\delta}\}$$
for all $x\in M$, where $c_1,c_2>0$ depend only on $n$, $(M,g)$ and $A>0$ such that $1/A<\tilde{u}_0<A$
Therefore, with $r:=d_{g_\xi}(x,\xi)$,
$$|(\tu+U)^{\crit-1}-\tu^{\crit-1}-U^{\crit-1}|\leq C{\bf 1}_{r\geq c_1\sqrt{\delta}}U+C{\bf 1}_{r\leq c_2\sqrt{\delta}}U^{\crit-2}.$$
Since $U\leq C\delta^{\frac{n-2}{2}}(\delta^2+r^2)^{1-n/2}$, we then obtain
\begin{equation}
\Vert (\tu+U)^{\crit-1}-\tu^{\crit-1}-U^{\crit-1}\Vert_{\frac{2n}{n+2}}\leq C\delta^{\frac{n+2}{4}}\hbox{ for }n\geq 7.\label{est:u:r}
\end{equation}
Since $0<\crit-2< 1$, we have 
$$|(\tu+U)^{\crit-2}-U^{\crit-2}|\leq C.$$
Therefore, with \eqref{def:tildeU} and \eqref{est:lem:1}, we obtain
\begin{equation}\label{est:u:r:2}
\Vert ((\tu+U)^{\crit-2}-U^{\crit-2})\partial_{p_i}U\Vert_{\frac{2n}{n+2}}\leq C\delta^{-1}\Vert U\Vert_{\frac{2n}{n+2}}\leq C\delta^{-1}\delta^2\hbox{ for }n\geq 7.
\end{equation}
Merging the estimates \eqref{control:R:2}, \eqref{est:R:u}, \eqref{est:der:R:u}, \eqref{est:u:r} and \eqref{est:u:r:2}, we obtain \eqref{est:R:2}. This ends the proof of Proposition~\ref{prop:est:R:2}.
\endproof

Plugging \eqref{est:R:2} and \eqref{est:R:2} together with \eqref{est:fin:1}, \eqref{control:R:1} and \eqref{est:c1}, we obtain that
\begin{equation}\label{est:J:21bis}
J_h(W+\Phi)=J_h(W)+\bigO(\Vert \Delta_g \tu+h\tu-\tu^{\crit-1}\Vert_{\infty}^2+ D_{h,\xi}^2\delta^4+\delta^8+\delta^{n-2})
\end{equation}
and 
\begin{multline}\label{est:der:J:2bis}
\partial_{p_i}J_h(W+\Phi)=\partial_{p_i}J_h(W)+\bigO(\Vert \Delta_g \tu+h\tu-\tu^{\crit-1}\Vert_{\frac{2n}{n+2}}\delta\\
+\Vert \Delta_g \tu+h\tu-\tu^{\crit-1}\Vert_{\infty}^{\crit-1}\delta^{-1}+(D_{h,\xi}+\delta^2+\delta^{\frac{n-6}{2}})^{\crit-1}\delta^{\frac{n+6}{n-2}}+D_{h,\xi}\delta^3+\delta^5+\delta^{n/2})
\end{multline}
for all $i=0,\dotsc,n$. We now estimate $J_h(W+\Phi)$:

\begin{proposition}\label{prop:est:J:2}
Assume that $n\geq 7$ and $u_0,\tilde{u}_0>0$. Then
\begin{multline}
J_h(W+\Phi)=J_h(\tu)+\frac{1}{n}\int_{\rn}U_{1,0}^{\crit}\,dx+\frac{1}{2}\varphi_h(\xi)\delta^2\int_{\rn}U_{1,0}^2\, dx\label{est:c0phi}\\
-\frac{1}{4n}K_{h_0}(\xi_0)\delta^4\int_{\rn}|x|^2U_{1,0}^2\, dx+\smallo(\delta^4)-u_0(\xi_0)\delta^{\frac{n-2}{2}}\int_{\rn}U_{1,0}^{\crit-1} dx\\
+\bigO(\Vert \Delta_g \tu+h\tu-\tu^{\crit-1}\Vert_{\infty}^2+\delta^{\frac{n-2}{2}}(\Vert\Delta_g\tu+h\tu-\tu^{\crit-1}\Vert_\infty+\Vert \tu-u_0\Vert_\infty+\smallo(1)))
\end{multline}
as $\delta\to0$, $\xi\to\xi_0$ and $h\to h_0$ in $C^2(M)$.
\end{proposition}

\proof[Proof of Proposition~\ref{prop:est:J:2}]
We first write that
\begin{multline*}
J_h(\tu+U)= J_h(\tu)+J_h(U)-\int_M \tu U^{\crit-1} dv_g+\int_M (\Delta_g\tu+h\tu-\tu^{\crit-1})U\, dv_g\\-\frac{1}{\crit}\int_M ((\tu+U)^{\crit}-\tu^{\crit}-U^{\crit}-\crit\tu^{\crit-1}U-\crit \tu U^{\crit-1})\,dv_g.
\end{multline*}
We fix $0<\theta<\frac{2}{n-2}<\crit-2$. There exists $C>0$ such that
\begin{multline*}
|(\tu+U)^{\crit}-\tu^{\crit}-U^{\crit}-\crit\tu^{\crit-1}U-\crit \tu U^{\crit-1}|\\
\leq C{\bf 1}_{\tu\leq U}\tu^{1+\theta}U^{\crit-1-\theta}+C{\bf 1}_{U\leq \tu}\tu^{\crit-1-\theta}U^{1+\theta}.
\end{multline*}
Using the definition \eqref{def:peak} and arguing as in the proof of \eqref{est:u:r}, we obtain
$$\left|\int_M ((\tu+U)^{\crit}-\tu^{\crit}-U^{\crit}-\crit\tu^{\crit-1}U-\crit \tu U^{\crit-1})\,dv_g\right|\leq C \delta^{\frac{n-2}{2}+\frac{n-2}{2}\theta}.$$
Furthermore, we obtain 
\begin{align*}
\left|\int_M (\Delta_g\tu+h\tu-\tu^{\crit-1})U\, dv_g\right|&\le C\Vert\Delta_g\tu+h\tu-\tu^{\crit-1})\Vert_\infty\int_MUdv_g\\
&\le C\Vert\Delta_g\tu+h\tu-\tu^{\crit-1})\Vert_\infty\delta^{\frac{n-2}{2}}.
\end{align*}
Using \eqref{def:peak}, that $\Lambda_\xi(x)=1+\bigO(d_g(x,\xi)^2)$ for all $x\in M$ and that $U_{\delta,0}$ is radially symmetrical, we obtain
\begin{align*}
&\int_M \tu U^{\crit-1} dv_g= \int_{B_{r_0}(0)}\tu(\exp_\xi^{g_\xi}(x))(1+\bigO(|x|^2))U_{\delta,0}^{\crit-1} dx+\bigO(\delta^{\frac{n-2}{2}(\crit-1)})\\
&\ = \int_{B_{r_0}(0)}(\tu(\xi)+x^\alpha\partial_{x_\alpha}\tu(\exp_\xi^{g_\xi}(\xi))+\bigO(|x|^2))U_{\delta,0}^{\crit-1} dx+\bigO(\delta^{\frac{n+2}{2}})\\
&\ = \tu(\xi)\int_{B_{r_0}(0)}U_{\delta,0}^{\crit-1} dx+\bigO\left(\int_{B_{r_0}(0)}|x|^2U_{\delta,0}^{\crit-1} dx\right)+\bigO(\delta^{\frac{n+2}{2}})\\
&\ = \tu(\xi)\delta^{\frac{n-2}{2}}\int_{B_{r_0/\delta}(0)}U_{1,0}^{\crit-1}dx+\bigO\left(\delta^{\frac{n+2}{2}}\int_{B_{r_0/\delta}(0)}|x|^2U_{1,0}^{\crit-1} dx\right)+\bigO(\delta^{\frac{n+2}{2}}).
\end{align*}
Since $U_{1,0}\leq C(1+|x|^2)^{1-n/2}$, we obtain
$$\int_{B_{r_0/\delta}(0)}U_{1,0}^{\crit-1} dx=\int_{\rn}U_{1,0}^{\crit-1} dx+\bigO(\delta^2)$$
and
$$\int_{B_{r_0/\delta}(0)}|x|^2U_{1,0}^{\crit-1} dx=\bigO\left(\ln(1/\delta)\right)\hbox{ for }n\ge7.$$
Therefore, plugging all these estimates together yields
$$\int_M \tu U^{\crit-1} dv_g= \tu(\xi)\delta^{\frac{n-2}{2}}\int_{\rn}U_{1,0}^{\crit-1} dx+\bigO(\delta^{\frac{n+2}{2}}\ln(1/\delta)).$$
Consequently, we obtain that for every $0<\theta<\frac{2}{n-2}$,
\begin{multline*}
J_h(\tu+U)= J_h(\tu)+J_h(U)-\tu(\xi)\delta^{\frac{n-2}{2}}\int_{\rn}U_{1,0}^{\crit-1} dx\\
+\int_M (\Delta_g\tu+h\tu-\tu^{\crit-1})U\, dv_g+\bigO(\delta^{\frac{n-2}{2}+\frac{n-2}{2}\theta}).
\end{multline*}
Now, with the expansion \eqref{est:J:W}, we obtain that for $n\geq 7$,
\begin{multline}
J_h(\tu+U)= J_h(\tu)+\frac{1}{n}\int_{\rn}U_{1,0}^{\crit}\,dx+\frac{1}{2}\varphi_h(\xi)\delta^2\int_{\rn}U_{1,0}^2\, dx\label{est:c0}\\
-\frac{1}{4n}K_{h_0}(\xi_0)\delta^4\int_{\rn}|x|^2U_{1,0}^2\, dx+\smallo(\delta^4)-u_0(\xi_0)\delta^{\frac{n-2}{2}}\int_{\rn}U_{1,0}^{\crit-1} dx\\
+\bigO(\delta^{\frac{n-2}{2}}(\Vert\Delta_g\tu+h\tu-\tu^{\crit-1}\Vert_\infty+\Vert \tu-u_0\Vert_\infty+d_g(\xi,\xi_0)+\delta^{\frac{n-2}{2}\theta})).
\end{multline}
Plugging together \eqref{est:J:21bis} and \eqref{est:c0}, we then obtain \eqref{est:c0phi}. This ends the proof of Proposition~\ref{prop:est:J:2}.
\endproof

We now estimate the derivatives of $J_h(W+\Phi)$:

\begin{proposition}\label{prop:est:der:J:2}
Assume that $n\geq 7$ and $u_0,\tilde{u}_0>0$. Then
\begin{multline}
\partial_\delta J_h(W+\Phi)= \varphi_h(\xi)\delta\int_{\rn}U_{1,0}^2\, dx- \frac{1}{n}K_{h_0}(\xi_0)\delta^3\int_{\rn}|x|^2U_{1,0}^2\, dx+\smallo(\delta^3)\label{der:J:u:0}\\
-\frac{n-2}{2}u_0(\xi_0)\delta^{\frac{n-4}{2}}\int_{\rn}U_{1,0}^{\crit-1} dx+\bigO(\delta^{\frac{n-4}{2}}(\Vert \tu-u_0\Vert_\infty+\smallo(1))\\
+\Vert \Delta_g \tu+h\tu-\tu^{\crit-1}\Vert_\infty\delta+\Vert \Delta_g \tu+h\tu-\tu^{\crit-1}\Vert_\infty^{\crit-1}\delta^{-1}+D_{h,\xi}^{\crit-1}\delta^{\frac{n+6}{n-2}})
\end{multline}
and
\begin{multline}
\partial_{\xi_i} J_h(W+\Phi)=\frac{1}{2}\partial_{\xi_i}\varphi_h(\xi)\delta^2\int_{\rn}U_{1,0}^2\, dx+\smallo(\delta^3)\label{der:J:u:i}\\
+\bigO(\delta^{\frac{n-4}{2}}(\Vert \tu-u_0\Vert_\infty+\smallo(1))
+\Vert \Delta_g \tu+h\tu-\tu^{\crit-1}\Vert_\infty\delta\\
+\Vert \Delta_g \tu+h\tu-\tu^{\crit-1}\Vert_\infty^{\crit-1}\delta^{-1}+D_{h,\xi}^{\crit-1}\delta^{\frac{n+6}{n-2}})
\end{multline}
for all $i=1,\dotsc,n$, as $\delta\to0$, $\xi\to\xi_0$ and $h\to h_0$ in $C^2(M)$.
\end{proposition}

\proof[Proof of Proposition~\ref{prop:est:der:J:2}]
We fix $i\in\{0,\dotsc,n\}$. We have
\begin{multline*}
\partial_{p_i}J_h(\tu+U)=\int_M (\Delta_g\tu+h\tu-\tu^{\crit-1})\partial_{p_i}U\, dv_g-(\crit-1)\int_M \tu U^{\crit-2}\partial_{p_i}U\, dv_g\\
+\partial_{p_i}J_h(U)-\int_M ((\tu+U)^{\crit-1}-U^{\crit-1}-(\crit-1)\tu U^{\crit-2})\partial_{p_i}U\, dv_g.
\end{multline*}
There exists $C>0$ such that
\begin{align*}
&|(\tu+U)^{\crit-1}-\tu^{\crit-1}-U^{\crit-1}-(\crit-1) \tu U^{\crit-2}|\\
&\qquad\qquad\qquad\qquad\qquad\qquad\qquad\quad\leq C{\bf 1}_{\tu\leq U}\tu^{\crit-1}+C{\bf 1}_{U\leq \tu}U^{\crit-1}.
\end{align*}
Since $|\partial_{p_i}U|\leq C \tilde{U}/\delta$ (see \eqref{est:lem:1}), arguing as in the proof of \eqref{est:u:r}, we obtain
$$\left|\int_M ((\tu+U)^{\crit-1}-U^{\crit-1}-(\crit-1)\tu U^{\crit-2})\partial_{p_i}U\, dv_g\right|\leq C\delta^{\frac{n-2}{2}}.$$
Furthermore, we obtain 
\begin{align*}
\left|\int_M (\Delta_g\tu+h\tu-\tu^{\crit-1})\partial_{p_i}U\, dv_g\right|&\le C\Vert\Delta_g\tu+h\tu-\tu^{\crit-1})\Vert_\infty\delta^{-1}\int_M\tilde{U}dv_g\\
&\le C\Vert\Delta_g\tu+h\tu-\tu^{\crit-1})\Vert_\infty\delta^{-1}\delta^{\frac{n-2}{2}}.
\end{align*}
Independently, using again \eqref{est:lem:1}, straightforward computations yield
\begin{align*}
&\int_M \tu U^{\crit-2}\partial_{p_i}U\, dv_g= \int_M \left(u_0(\xi_0)+\bigO(\Vert \tu-u_0\Vert_\infty+d_g(.,\xi_0)\right) U^{\crit-2}\partial_{p_i}U\, dv_g\\
&\qquad= u_0(\xi_0) \int_M U^{\crit-2}\partial_{p_i}U\, dv_g\\
&\qquad\quad+\bigO\left(\delta^{-1}\int_M(\Vert \tu-u_0\Vert_\infty+d_g(\xi,\xi_0)+d_g(.,\xi)) \tilde{U}^{\crit-1}dv_g\right)\\
&\qquad= u_0(\xi_0) \int_M U^{\crit-2}\partial_{p_i}U\, dv_g+\bigO(\delta^{-1}\delta^{\frac{n-2}{2}}(\Vert \tu-u_0\Vert_\infty+d_g(\xi,\xi_0)+\delta)).
\end{align*}
Arguing as in the proof of \eqref{est:24}, we obtain
\begin{align*}
\int_M U^{\crit-2}\partial_{p_i}U\, dv_g&=\int_{B_{r_0}(0)} (\Lambda_\xi U)^{\crit-2}\partial_{p_i}(\Lambda_\xi^{-1} U)(\exp_\xi^{g_\xi}(x))\,dx\\
&\qquad\qquad+\bigO\left(\delta^{-1}\int_{B_{r_0}(0)}|x|\tilde{U}^{\crit-1} dx\right)\\
&=\int_{B_{r_0}(0)} (\Lambda_\xi U)^{\crit-2}\partial_{p_i}(\Lambda_\xi^{-1} U)(\exp_\xi^{g_\xi}(x))\,dx+\bigO(\delta^{\frac{n-2}{2}}).
\end{align*}
We first deal with the case where $i=0$, that is $\partial_{p_i}=\partial_{p_0}=\partial_\delta$. With \eqref{der:U}, we obtain
\begin{multline*}
\int_{B_{r_0}(0)} (\Lambda_\xi U)^{\crit-2}\partial_{\delta}(\Lambda_\xi^{-1} U)(\exp_\xi^{g_\xi}(x))\,dx=\int_{B_{r_0}(0)} U_{\delta,0}^{\crit-2}\partial_{\delta}U_{\delta,0}\, dx\\
= \int_{B_{r_0}(0)} U_{\delta,0}^{\crit-2}\partial_{\delta}U_{\delta,0} dx=\delta^{-1}\int_{B_{r_0}(0)} (\delta^{-\frac{n-2}{2}}U_{1,0}(\delta^{-1}x))^{\crit-2}\delta^{-\frac{n-2}{2}}Z_{0}(\delta^{-1}x)\,dx\\
= \delta^{-1}\delta^{\frac{n-2}{2}}\int_{B_{r_0/\delta}(0)} U_{1,0}^{\crit-2}Z_0\, dx.
\end{multline*}
Since $Z_0\leq C U_{1,0}$, an asymptotic estimate yields
$$\int_{B_{r_0}(0)} (\Lambda_\xi U)^{\crit-2}\partial_\delta(\Lambda_\xi^{-1} U)(\exp_\xi^{g_\xi}(x))\,dx=\delta^{-1}\delta^{\frac{n-2}{2}}\int_{\rn} U_{1,0}^{\crit-2}Z_0\, dx+\bigO(\delta^{\frac{n}{2}}).$$
Note that for every $\delta>0$, we have 
$$\int_{\rn}U_{\delta,0}^{\crit-1} dx=\delta^{\frac{n-2}{2}}\int_{\rn}U_{1,0}^{\crit-1} dx.$$
Differentiating this equality with respect to $\delta$ at $1$, we obtain
$$(\crit-1)\int_{\rn}U_{1,0}^{\crit-2}Z_0\, dx=\frac{n-2}{2}\int_{\rn}U_{1,0}^{\crit-1} dx.$$
Therefore, we obtain
$$(\crit-1)\int_M U^{\crit-2}\partial_{\delta}U\, dv_g=\frac{n-2}{2}\delta^{-1}\delta^{\frac{n-2}{2}}\int_{\rn}U_{1,0}^{\crit-1} dx+\bigO(\delta^{\frac{n-2}{2}}).$$
We now deal with the case $i\ge1$, that is $\partial_{p_i}=\partial_{\xi_i}$. It follows from \eqref{eq:31:bis} and \eqref{eq:d} that
\begin{align*}
&\int_{B_{r_0}(0)} (\Lambda_\xi U)^{\crit-2}\partial_{\xi_i}(\Lambda_\xi^{-1} U)(\exp_\xi^{g_\xi}(x))\,dx\\
&=\int_{B_{r_0}(0)} U_{\delta,0}^{\crit-2}\left(-\frac{n-2}{2}\right)\frac{\delta^{\frac{n-2}{2}}}{(\delta^2+|x|^2)^{n/2}}\left(-2x_i+\bigO(|x|^3)\right) dx\\
&=\bigO\left(\int_{B_{r_0}(0)} U_{\delta,0}^{\crit-2}\frac{U_{\delta,0}}{\delta^2+|x|^2}|x|^3\,dx\right)=\bigO\left(\int_{B_{r_0}(0)} |x|U_{\delta,0}^{\crit-1}dx\right)=\bigO(\delta^{\frac{n-2}{2}}).
\end{align*}
Putting these results together yields
\begin{align*}
&\partial_{\xi_i}J_h(\tu+U)=\partial_{\xi_i}J_h(U)-\frac{n-2}{2}\eps_{i,0}u_0(\xi_0)\delta^{-1}\delta^{\frac{n-2}{2}}\int_{\rn}U_{1,0}^{\crit-1} dx\\
& +\bigO(\delta^{-1}\delta^{\frac{n-2}{2}}(\Vert \Delta_g\tu+h\tu-\tu^{\crit-1}\Vert_{\infty}+\Vert \tu-u_0\Vert_\infty+d_g(\xi,\xi_0)+\delta))
\end{align*}
for all $i=0,\dotsc,n$. Using the estimates \eqref{der:xi:0} and \eqref{der:xi:1} for the derivatives of $J_h(U_{\delta,\xi})$, we obtain
\begin{multline*}
\partial_\delta J_h(\tu+U)= \varphi_h(\xi)\delta^{-1}\delta^2\int_{\rn}U_{1,0}^2\, dx- 4K_{h_0}(\xi_0)\delta^3\int_{\rn}|x|^2U_{1,0}^2\, dx\\
-\frac{n-2}{2}u_0(\xi_0)\delta^{\frac{n-4}{2}}\int_{\rn}U_{1,0}^{\crit-1} dx+\smallo(\delta^3)\\
+\bigO(\delta^{\frac{n-4}{2}}(\Vert \Delta_g\tu+h\tu-\tu^{\crit-1}\Vert_{\infty}+\Vert \tu-u_0\Vert_\infty+d_g(\xi,\xi_0)+\delta))
\end{multline*}
and
\begin{multline*}
\partial_{\xi_i} J_h(\tu+U)=\frac{1}{2}\partial_{\xi_i}\varphi_h(\xi)\delta^2\int_{\rn}U_{1,0}^2\, dx+\smallo(\delta^4)\\
+\bigO(\delta^{\frac{n-4}{2}}(\Vert \Delta_g\tu+h\tu-\tu^{\crit-1}\Vert_{\infty}+\Vert \tu-u_0\Vert_\infty+d_g(\xi,\xi_0)+\delta)).
\end{multline*}
With \eqref{est:der:J:2bis}, we then obtain \eqref{der:J:u:0} and \eqref{der:J:u:i}. This ends the proof of Proposition~\ref{prop:est:der:J:2}.
\endproof

Theorem~\ref{th:main:u0} for $n\geq 7$ will be proved in Section~\ref{sec:pf:2}.

\section{Energy and remainder estimates: the case $n=6$ and $u_0,\tilde{u}_0>0$}\label{Sec7} 

In this section, we assume that $u_0,\tilde{u}_0>0$ and $n=6$, that is $\crit-1=2$. Here again, we set $B_{h,\delta,\xi}\equiv 0$, so that $W_{h,\tu,\delta,\xi}=W_{\tu,\delta,\xi}\equiv\tu+U_{\delta,\xi}$ and the assumptions of Proposition~\ref{prop:c1} are satisfied. The remark underlying this section is that
$$\Delta_g(u_0+U)+h(u_0+U)-(u_0+U)^2=\Delta_gU+(h-2u_0)U-U^2.$$
Therefore, to obtain a good approximation of the blowing-up solution, we will subtract a perturbation of $2u_0$ to the potential. We first estimate $R=R_{\delta,\xi}$:

\begin{proposition}\label{prop:est:R:3}
Assume that $n=6$ and $u_0,\tilde{u}_0>0$. Then
\begin{equation}\label{est:R:3}
\Vert R\Vert_{3/2}+\delta\Vert \partial_pR\Vert_{3/2}\leq C\Vert  \Delta_g\tu+h\tu-\tu^2\Vert_\infty+ C \delta^{2} (1+\overline{D}_{h,\xi}\left(\ln(1/\delta)\right)^{2/3}),
\end{equation}
where
\begin{equation}\label{Dhxi2}
\overline{D}_{h,\xi}:=\Vert \bar{h}-\bar{h}_0\Vert_\infty+d_g(\xi,\xi_0)^2.
\end{equation}
\end{proposition}

\proof[Proof of Proposition~\ref{prop:est:R:3}]
Since $\crit-1=2$, we have
\begin{align*}
R&= \Delta_g(\tu+U)+h(\tu+U)-(\tu+U)^2\\
&= \Delta_g\tu+h\tu-\tu^2+ \Delta_gU+(h-2\tu)U-U^2
\end{align*}
and
$$\partial_{p_i}R= \partial_{p_i}\left(\Delta_gU+(h-2\tu)U-U^2\right)$$
for all $i=0,\dotsc,n$. For convenience, we write 
$$\bar{h}:=h-2\tu\hbox{ and }\bar{h}_0:=h_0-2u_0.$$ 
The estimate \eqref{est:R:3} then follows from \eqref{control:R:2}. This ends the proof of Proposition~\ref{prop:est:R:3}.
\endproof

We now estimate the derivatives of $J_h(W+\Phi)$:

\begin{proposition}\label{prop:est:der:J:3}
Assume that $n=6$ and $u_0,\tilde{u}_0>0$. Then
\begin{multline}
J_h(W+\Phi)= J_h(\tu)+\frac{1}{n}\int_{\rn}U_{1,0}^{\crit}\,dx+\frac{1}{2}\varphi_{h,\tilde{u}_0}(\xi)\delta^2\int_{\rn}U_{1,0}^2\, dx\label{co:6}\\
-24\omega_5K_{h_0,u_0}(\xi_0)  \delta^4\ln(1/\delta)+\bigO(\Vert\Delta_g\tu+h\tu-\tu^2\Vert_\infty^2+\Vert\Delta_g\tu+h\tu-\tu^2\Vert_\infty\delta^2)\\
+\bigO(\delta^4\ln(1/\delta)(\smallo(1)+\overline{D}_{h,\xi}^2\left(\ln(1/\delta)\right)^{1/3})),
\end{multline}
\begin{multline}
\partial_\delta J_h(W+\Phi)= \varphi_{h,\tilde{u}_0}(\xi)\delta\int_{\rn}U_{1,0}^2\, dx-96\omega_5K_{h_0,u_0}(\xi_0)\delta^3\ln(1/\delta)\label{der:0:6}\\
+\bigO(\Vert\Delta_g\tu+h\tu-\tu^2\Vert_\infty\delta+\Vert\Delta_g\tu+h\tu-\tu^2\Vert_\infty^2\delta^{-1})\\
+\bigO(\delta^3\ln(1/\delta)(\smallo(1)+ \overline{D}_{h,\xi}^2\left(\ln(1/\delta)\right)^{1/3}))
\end{multline}
and
\begin{multline}
\partial_{\xi_i} J_h(W+\Phi)= \frac{1}{2}\partial_{\xi_i}\varphi_{h,\tilde{u}_0}(\xi)\delta^2\int_{\rn}U_{1,0}^2\, dx\label{der:i:6}\\
+\bigO(\Vert\Delta_g\tu+h\tu-\tu^2\Vert_\infty\delta+\Vert\Delta_g\tu+h\tu-\tu^2\Vert_\infty^2\delta^{-1})\\
+\bigO(\delta^3\ln(1/\delta)(\smallo(1)+\overline{D}_{h,\xi}^2\left(\ln(1/\delta)\right)^{1/3}))
\end{multline}
for all $i=1,\dotsc,n$, as $\delta\to0$, $\xi\to\xi_0$ and $h\to h_0$ in $C^2(M)$, where $\varphi_{h,\tilde{u}_0}$, $K_{h_0,u_0}(\xi_0)$ and $\overline{D}_{h,\xi}$ are as in \eqref{defphi0th3}, \eqref{Khu} and \eqref{Dhxi2}.
\end{proposition}

\proof[Proof of Proposition~\ref{prop:est:der:J:3}]
As one checks, since $n=6$ and $\crit=3$, we have
$$J_h(\tu+U)=J_h(\tu)+J_{\bar{h}}(U)+\int_M\left(\Delta_g\tu+h\tu-\tu^2\right) U\, dv_g$$
and
$$\partial_{p_i}J_h(\tu+U)=\partial_{p_i}J_{\bar{h}}(U)+\int_M\left(\Delta_g\tu+h\tu-\tu^2\right)\partial_{p_i} U\, dv_g$$
for all $i=0,\dotsc,n$. Using the definition \eqref{def:peak} and since $|\partial_{p_i}U|\le C\tilde{U}/\delta$, we obtain
\begin{align*}
\left|\int_M (\Delta_g\tu+h\tu-\tu^2)U\, dv_g\right|&\le C\Vert\Delta_g\tu+h\tu-\tu^2)\Vert_\infty\int_MUdv_g\\
&\le C\Vert\Delta_g\tu+h\tu-\tu^2)\Vert_\infty\delta^2.
\end{align*}
and
\begin{align*}
\left|\int_M (\Delta_g\tu+h\tu-\tu^2)\partial_{p_i}U\, dv_g\right|&\le C\Vert\Delta_g\tu+h\tu-\tu^2)\Vert_\infty\delta^{-1}\int_M\tilde{U}dv_g\\
&\le C\Vert\Delta_g\tu+h\tu-\tu^2)\Vert_\infty\delta^{-1}\delta^2.
\end{align*}
Putting these estimates together with \eqref{phi0th3}, \eqref{est:fin:1}, \eqref{est:c1}, \eqref{est:R:3}, \eqref{est:J:W}, \eqref{der:xi:0} and \eqref{der:xi:1}, we obtain \eqref{co:6}, \eqref{der:0:6} and \eqref{der:i:6}. This ends the proof of Proposition~\ref{prop:est:der:J:3}.
\endproof

Theorem~\ref{th:main:u0} for $n=6$ will be proved in Section~\ref{sec:pf:2}.

\section{Setting and definition of the mass in dimensions $n=3,4,5$}\label{sec:mass}

In this section, we assume that $n\leq 5$. Our first lemma is a simple computation:

\begin{lem}\label{lem:lap}
There exist two functions $(\xi,x)\mapsto f_i(\xi,x)$, $i=1,2$, defined and smooth on $M\times M$ such that for every function $f:\rn\to\rr$ that is radially symmetrical, we have
\begin{multline*}
(\Delta_g+h)(\chi(r)\Lambda_\xi(x)f(r))= \Lambda_\xi(x)^{\crit-1}\chi\Delta_{\Eucl}(f(r))+f_1(\xi,x)f'(r)+f_2(\xi,x)f(r)\\
+\hat{h}_\xi\chi(x)\Lambda_\xi(x)f(r)
\end{multline*}
for all $x\in M\backslash\left\{\xi\right\}$, where $r:=d_{g_\xi}(x,\xi)$ and $\hat{h}_\xi$ is as in \eqref{def:hat:h}. Furthermore, $f_i(\xi,x)=0$ when $d_{g}(x,\xi)\geq r_0$ and there exists $C_N>0$ such that
$$|f_1(\xi,x)(x)|\leq C_N d_g(x,\xi)^{N-1}\text{ and }|f_2(\xi,x)|\leq C_N d_g(x,\xi)^{N-2}\hbox{ for all }x,\xi\in M.$$
\end{lem}

The proof of Lemma~\ref{lem:lap} follows the computations in \eqref{eq:W:lap}. We leave the details to the reader.

\smallskip
We define
$$\Gamma_{\xi}(x):=\frac{\chi(d_{g_\xi}(x,\xi))\Lambda_\xi(x)}{(n-2)\omega_{n-1}d_{g_\xi}(x,\xi)^{n-2}}$$
for all $x\in M\backslash\left\{\xi\right\}$. It follows from Lemma~\ref{lem:lap} and the definition \eqref{def:U:delta} that
\begin{equation}\label{exp:U:delta}
\Delta_gU_{\delta,\xi}+hU_{\delta,\xi}=U_{\delta,\xi}^{\crit-1}+F_\delta(\xi,x)\delta^{\frac{n-2}{2}}+\hat{h}_\xi U_{\delta,\xi}
\end{equation}
and
$$(\Delta_g+h )\Gamma_{\xi}=\frac{F_0(\xi,x)}{k_n}+\hat{h}_\xi \Gamma_{\xi},$$
where
$$k_n:=(n-2)\omega_{n-1}\sqrt{n(n-2)}^{\frac{n-2}{2}}$$
and $(t,\xi,x)\to F_t(\xi,x)$ is of class $C^p$ on $[0,\infty)\times M\times M$, with $p$ being as large as we want provided we choose $N$ large enough. This includes $t=0$ and, therefore,
\begin{equation}\label{lim:F}
\lim_{t\to 0}F_t=F_0\hbox{ in }C^p(M\times M).
\end{equation}
For every $t\geq 0$, we define $\beta_{h,t,\xi}\in H_1^2(M)$ as the unique solution to
\begin{align}\label{def:beta:h}
(\Delta_g+h)\beta_{h,t,\xi}&=-\left(\frac{F_t(\xi,x)}{k_n}+\hat{h}_\xi \frac{\chi(d_{g_\xi}(\xi,x))\Lambda_\xi(x)}{(n-2)\omega_{n-1}(t^2+d_{g_\xi}(\xi,x)^2)^{\frac{n-2}{2}}}\right)\\
&=-\frac{F_t(\xi,x)}{k_n}-\hat{h}_\xi\left\{\begin{aligned}&\frac{U_{t,\xi}}{k_nt^{\frac{n-2}{2}}}&&\text{if }t>0\\&\Gamma_\xi&&\text{if }t=0.\end{aligned}\right.\nonumber
\end{align}
Since $N>n-2$ and $n\leq 5$, the right-hand-side is uniformly bounded in $L^q(M)$ for some $q>\frac{2n}{n+2}$, independently of $t\geq 0$, $\xi\in U_0$ and $h\in C^2(M)$ satisfying $\left\|h\right\|_\infty<A$. Furthermore, we have $\lambda_1(\Delta_g+h)>1/A$. Therefore, $\beta_{h,t,\xi}$ is well defined and we have 
\begin{equation}\label{cv:beta:h}
\Vert\beta_{h,t,\xi}-\beta_{h,0,\xi}\Vert_{H_1^2}=\smallo(1)\hbox{ as }t\to0
\end{equation}
uniformly with respect to $\xi$ and $h$. Furthermore, we have $\beta_{h,t,\xi}\in C^2(M)$ when $t>0$. As one checks, with these definitions, we obtain that
$$G_{h,\xi}:=\Gamma_\xi+\beta_{h,0,\xi}$$
is the Green's function of the operator $\Delta_g+h$ at the point $\xi$. We now define the {\it mass} of $\Delta_g+h$ at the point $\xi$:

\begin{propdefi}\label{propdefi:1}
Assume that $3\leq n\leq 5$ and $N>n-2$. Let $h\in C^2(M)$ be such that $\Delta_g+h$ is coercive. In the case where $n\in\left\{4,5\right\}$, assume in addition that there exists $\xi\in M$ such that $\varphi_h(\xi)=|\nabla \varphi_h(\xi)|=0$, where $\varphi_h$ is as in \eqref{defphi0th3}. Then $\beta_{h,0,\xi}\in C^0(M)$. Furthermore, the number $\beta_{h,0,\xi}(\xi)$ does not depend on the choice of $N>n-2$ and $g_\xi$ satisfying \eqref{Lambda} and \eqref{eq:elt:vol}. We then define the mass of $\Delta_g+h$ at the point $\xi$ as $m_{h}(\xi):=\beta_{h,0,\xi}(\xi)$.
\end{propdefi}

\proof[Proof of Proposition-Definition~\ref{propdefi:1}] 
As one checks, when $n=3$, we have 
$$\hat{h}_\xi(x)\Gamma_\xi(x)=\bigO(d_g(x,\xi)^{-1})$$
and when $n\in\left\{4,5\right\}$ and $\varphi_h(\xi)=|\nabla \varphi_h(\xi) |=0$, we have 
$$\hat{h}_\xi(x)\Gamma_\xi(x)=\bigO(d_g(x,\xi)^{4-n}).$$
Furthermore, we have
$$F_0(\xi,x)=\bigO(d_g(x,\xi)^{N-n}).$$
When $N>n-2$, this implies that $\beta_{h,0,\xi}\in C^0(M)$. The fact that the number $\beta_{h,0,\xi}(\xi)$ does not depend on the choice of $N$ and $g_\xi$ then follows from the uniqueness of conformal normal coordinates up to the action of $O(n)$ and the choice of the metric's one-jet at the point $\xi$ (see Lee--Parker~\cite{LP}). This ends the proof of Proposition-Definition~\ref{propdefi:1}.
\endproof

We now prove a differentiation result that will allow us to obtain Theorem~\ref{th:2}: 

\begin{proposition}\label{prop:diff:mass} 
Assume that $3\leq n\leq 5$. Let $h\in C^2(M)$ be such that $\Delta_g+h$ is coercive. In the case where $n\in\left\{4,5\right\}$, assume that there exists $\xi\in M$ such that $\varphi_h(\xi)=|\nabla \varphi_h(\xi) |=0$. Let  $H\in C^2(M)$ be such that $H(\xi)=|\nabla H(\xi)|=0$. Then $m_{h+\eps H}(\xi)$ is well defined for small $\eps\in \rr$ and differentiable with respect to $\eps$. Furthermore,
$$\partial_\eps(m_{h+\eps H}(\xi))_{|0}=-\int_M HG_{h,\xi}^2\,dv_g.$$
\end{proposition}

\proof[Proof of Proposition~\ref{prop:diff:mass}] 
In order to differentiate the mass with respect to the potential function $h$, it is convenient to write 
$$G_{h,\xi}=G_{c_n\Scal_g,\xi}+\hat{\beta}_{h,\xi},$$
where $\hat{\beta}_{h,\xi}\in H_1^2(M)$ is the solution to
\begin{equation}\label{def:hat}
(\Delta_g+h)\hat{\beta}_{h,\xi}=-\varphi_hG_{c_n\Scal_g,\xi}.
\end{equation}
Under the assumptions of the proposition, we have $\hat{\beta}_{h,\xi}\in C^0(M)$ and
$$\hat{\beta}_{h,\xi}(\xi)=-\int_M \varphi_hG_{c_n\Scal_g,\xi}G_{h,\xi}\, dv_g.$$
Furthermore, as one checks, we have
\begin{equation}\label{eq:m:hat}
m_h(\xi)=m_{c_n \Scal_g}(\xi)-\hat{\beta}_{h,\xi}(\xi).
\end{equation}
It follows from standard elliptic theory that $\hat{\beta}_{h+\epsilon H,\xi}$ is differentiable with respect to $\epsilon$. Differentiating \eqref{def:hat} then yields
$$(\Delta_g+h)\partial_\eps(\hat{\beta}_{h+\eps H,\xi})_{|0}+H \hat{\beta}_{h,\xi}=- HG_{c_n\Scal_g,\xi},$$
which gives
$$(\Delta_g+h)\partial_\eps(\hat{\beta}_{h+\eps H,\xi})_{|0}=- H G_{h,\xi}.$$
Therefore,
$$\partial_\eps(\hat{\beta}_{h+\eps H,\xi}(x))_{|0}=-\int_M G_{h,x} H G_{h,\xi}\, dv_g.$$
It then follows from \eqref{eq:m:hat} that
$$\partial_\eps(m_{h+\eps H}(\xi))_{|0}=-\int_M HG_{h,\xi}^2\,dv_g.$$
This proves Proposition~\ref{prop:diff:mass}.
\endproof

\section{Energy and remainder estimates in dimensions $n=3,4,5$}\label{sec:small:dim}

In this section, we assume that $n\le5$ and $u_0\equiv\tilde{u}_0\equiv0$. When $n\in\left\{4,5\right\}$, we assume in addition that \eqref{phi0th2} is satisfied. We define 
\begin{equation}\label{def:W:beta}
W_{h,\tu,\delta,\xi}=W_{h,\delta,\xi}:=U_{\delta,\xi}+B_{h,\delta,\xi}
\end{equation}
with
\begin{equation}\label{def:B:beta}
B_{h,\delta,\xi}:=k_n\delta^{\frac{n-2}{2}}\beta_{h,\delta,\xi}.
\end{equation}
In order to use the $C^1-$estimates of Proposition~\ref{prop:c1}, our first step is to obtain estimates for $\beta_{h,\delta,\xi}$ and its derivatives in $H_1^2(M)$:

\begin{proposition}\label{prop:est:beta:1} 
For $3\leq n\leq 5$, let $B_{h,\delta,\xi}$ be as in \eqref{def:B:beta}. Then \eqref{hyp:B} holds.
\end{proposition}

\proof[Proof of Proposition~\ref{prop:est:beta:1}]
It follows from \eqref{cv:beta:h} that
$$\Vert \beta_{h,\delta,\xi}\Vert_{H_1^2}\leq C.$$
Differentiating \eqref{def:beta:h} with respect to $\xi_i$, $i=1,\dotsc,n$, we obtain
$$(\Delta_g+h)(\partial_{\xi_i}\beta_{h,\delta,\xi})=-\frac{1}{k_n}\left(\partial_{\xi_i}F_\delta(\xi,\cdot)+\partial_{\xi_i}\hat{h}_\xi \frac{U_{\delta,\xi}}{\delta^{\frac{n-2}{2}}}+\hat{h}_\xi \frac{\partial_{\xi_i}U_{\delta,\xi}}{\delta^{\frac{n-2}{2}}}\right).$$
It follows from \eqref{lim:F} that 
$$\Vert\partial_{\xi_i}F_\delta(\xi,\cdot)\Vert_\infty\leq C$$
With the definition \eqref{def:hat:h} of $\hat{h}_\xi$, we obtain 
$$\partial_{\xi_i}\hat{h}_{\xi}=\partial_{\xi_i}(c_n \Scal_{g_\xi}\Lambda_\xi^{2-\crit})=\bigO(d_g(\cdot,\xi)).$$ 
Therefore, with \eqref{def:U:delta}, we obtain
$$\left|\partial_{\xi_i}\hat{h}_\xi \frac{U_{\delta,\xi}}{\delta^{\frac{n-2}{2}}}\right|\leq C\frac{d_g(x,\xi)}{(\delta^2+d_g(x,\xi)^2)^{\frac{n-2}{2}}}\,.$$
With \eqref{eq:30} and \eqref{eq:31}, we obtain
$$|\delta^{-\frac{n-2}{2}}\partial_{\xi_i}U_{\delta,\xi}|\leq C \frac{1}{(\delta^2+d_g(x,\xi)^2)^{\frac{n-2}{2}}}+C\frac{d_g(x,\xi)}{(\delta^2+d_g(x,\xi)^2)^{n/2}}\,.$$
The definition \eqref{def:hat:h} of $\hat{h}_\xi$ and the assumption $\varphi_{h_0}(\xi_0)=|\nabla\varphi_{h_0}(\xi_0)|=0$ yield 
\begin{equation}\label{esthath}
\hat{h}_\xi(x)=\bigO(d_g(x,\xi)^2+D_{h,\xi}),
\end{equation}
where $D_{h,\xi}$ is as in \eqref{Dhxi}. Putting together these inequalities yields
\begin{equation}
|(\Delta_g+h)(\partial_{\xi_i}\beta_{h,\delta,\xi})|\leq C+C \frac{d_g(x,\xi)}{(\delta^2+d_g(x,\xi)^2)^{\frac{n-2}{2}}}+CD_{h,\xi}\frac{\delta^2+d_{g}(x,\xi)}{(\delta^2+d_g(x,\xi)^2)^{n/2}}\,.\label{ineq:56}
\end{equation}
It then follows from elliptic theory and straightforward computations that
$$\Vert \partial_{\xi_i}\beta_{h,\delta,\xi}\Vert_{H_1^2}
\leq  C\left\{\begin{aligned}
&1 &&\hbox{if }n=3\\
&\left(\ln(1/\delta)\right)^{4/3}&&\hbox{if }n=4\\
&\delta^{-1/2} &&\hbox{if }n=5.
\end{aligned}\right.$$
Similarly, differentiating with respect to $\delta$, we obtain
\begin{align}\label{ineq:57}
\left|(\Delta_g+h)(\partial_{\delta}\beta_{h,\delta,\xi})\right|&=\left|-\frac{1}{k_n}\left(\partial_{\delta}F_\delta(\xi,\cdot)+\hat{h}_\xi \partial_{\delta}(\delta^{-\frac{n-2}{2}}U_{\delta,\xi})\right)\right|\\
&\leq  C+C\frac{\delta(d_g(x,\xi)^2+D_{h,\xi})}{(\delta^2+d_g(x,\xi)^2)^{n/2}}\nonumber
\end{align}
and, therefore, elliptic estimates and straightforward computations yield
$$\Vert \partial_{\delta}\beta_{h,\delta,\xi}\Vert_{\frac{2n}{n+2}}\leq C+C\left\Vert \frac{\delta}{(\delta^2+d_g(x,\xi)^2)^{n/2}}\right\Vert_{H_1^2}\leq C\left\{\begin{aligned}&1&&\hbox{if }n=3\\&\delta^{2-n/2}&&\hbox{if }n=4,5.\end{aligned}\right.$$
With the definition \eqref{def:B:beta}, all these estimates yield \eqref{hyp:B}. This ends the proof of Proposition~\ref{prop:est:beta:1}.
\endproof
 
The sequel of the analysis requires a pointwise control for $\beta_{h,\delta,\xi}$ and its derivatives. This is the objective of the following proposition:

\begin{proposition}\label{prop:est:beta:2} 
We have
\begin{equation}\label{bnd:beta}
|\beta_{h,\delta,\xi}(x)|\leq  C\left\{\begin{aligned}
&1&&\hbox{if }n=3\\
&1+|\ln\left(\delta^2+d_g(x,\xi)^2\right)|&&\hbox{if }n=4\\
&\left(\delta^2+d_g(x,\xi)^2\right)^{-1/2}&&\hbox{if }n=5,
\end{aligned}\right.
\end{equation}
\begin{equation}\label{bnd:beta:d}
|\partial_\delta\beta_{h,\delta,\xi}(x)|\leq  C+CD_{h,\xi}\delta\ln(1/\delta)\left(\delta^2+d_g(x,\xi)^2\right)^{-\frac{n-2}{2}}
\end{equation}
and
\begin{equation}\label{bnd:beta:xi}
|\partial_{\xi_i}\beta_{h,\delta,\xi}(x)|\leq  C+C\left\{\begin{aligned}
&D_{h,\xi}\left|\ln(\delta^2+d_g(x,\xi)^2)\right|&&\hbox{if }n=3\\
&D_{h,\xi}\left(\delta^2+d_g(x,\xi)^2\right)^{-1/2}&&\hbox{if }n=4\\
&\left|\ln(\delta^2+d_g(x,\xi)^2)\right|+D_{h,\xi}(\delta^2+d_g(x,\xi)^2)^{-1}&&\hbox{if }n=5
\end{aligned}\right.
\end{equation}
for all $i=1,\dotsc,n$, where $D_{h,\xi}$ is as in \eqref{Dhxi}.
\end{proposition}

\proof[Proof of Proposition~\ref{prop:est:beta:2}] 
These estimates will be consequences of Green's representation formula and Giraud's Lemma. More precisely, it follows from \eqref{def:beta:h} that
\begin{equation}\label{pointwise:beta}
\beta_{h,\delta,\xi}(x)=-\int_M G_{h,x}(y)\left(\frac{F_\delta(\xi,y)}{k_n}+\hat{h}_\xi \frac{\chi(d_{g_\xi}(y,\xi))\Lambda_\xi(y)}{(n-2)\omega_{n-1}(\delta^2+d_{g_\xi}(y,\xi)^2)^{\frac{n-2}{2}}}\right)dv_g(y)
\end{equation}
for all $x\in M$. With \eqref{lim:F} and the standard estimates of the Green's function $0<G_{h,x}(y)\leq C d_g(x,y)^{2-n}$ for all $x,y\in M$, $x\neq y$, we obtain
\begin{equation}\label{pointwise:beta:2}
|\beta_{h,\delta,\xi}(x)|\leq C +C\int_M  \frac{d_g(x,y)^{2-n}}{(\delta^2+d_{g}(y,\xi)^2)^{\frac{n-2}{2}}}\,dv_g(y).
\end{equation}
Recall Giraud's Lemma (see~\cite{DHR} for the present statement): For every $\alpha, \beta$ such that $0<\alpha,\beta<n$ and $x,z\in M$, $x\neq z$, we have
$$\int_M d_g(x,y)^{\alpha-n}d_g(y,z)^{\beta-n} dv_g(z)\leq C \left\{\begin{aligned}
&d_g(x,z)^{\alpha+\beta-n} &&\hbox{if }\alpha+\beta<n\\
&1+|\ln d_g(x,z)| &&\hbox{if }\alpha+\beta=n\\
&1 &&\hbox{if }\alpha+\beta>n.
\end{aligned}\right.$$
Therefore, \eqref{pointwise:beta:2} yields \eqref{bnd:beta} when $d_g(x,\xi)\geq \delta$. When $d_g(x,\xi)\leq\delta$, \eqref{pointwise:beta:2} yields
$$|\beta_{h,\delta,\xi}(x)|\leq C +C\int_M  \frac{d_g(x,y)^{2-n}}{(\delta^2+d_{g}(y,x)^2)^{\frac{n-2}{2}}}\,dv_g(y),$$
which in this case also yields \eqref{bnd:beta}. To prove \eqref{bnd:beta:xi}, we use \eqref{ineq:56} and the same method as for \eqref{bnd:beta}. The inequality \eqref{bnd:beta:d} is a little more delicate. With \eqref{ineq:57} and Green's identity, we obtain
\begin{align*}
|\partial_\delta\beta_{h,\delta,\xi}(x)|&=\left|\int_M G_{h,x}(y)\left(\Delta_g+h\right)\partial_\delta\beta_{h,\delta,\xi}(y)\,dv_g(y)\right|\\
&\leq  C+C\int_M d_g(x,y)^{2-n}\frac{\delta(d_g(y,\xi)^2+D_{h,\xi})}{(\delta^2+d_g(y,\xi)^2)^{n/2}}\,dv_g(y).
\end{align*}
We then obtain
\begin{multline*}
|\partial_\delta\beta_{h,\delta,\xi}(x)|\leq  C+C\delta\int_M d_g(x,y)^{2-n}d_g(y,\xi)^{2-n} dv_g(y)\\
+C\delta D_{h,\xi}\int_M\frac{d_g(x,y)^{2-n}}{(\delta^2+d_g(y,\xi)^2)^{n/2} }\,dv_g(y).
\end{multline*}
We estimate the first two terms in the right-hand side by using Giraud's lemma as in the proof of \eqref{bnd:beta}. We split the integral of the third term as
$$\int_M\frac{d_g(x,y)^{2-n}}{(\delta^2+d_g(y,\xi)^2)^{n/2} }\,dv_g(y)=\int_{\{d_g(x,y)<d_g(x,\xi)/2\}}+\int_{\{d_g(x,y)\geq d_g(x,\xi)/2\}}.$$
Since $d_g(y,\xi)>d_g(x,\xi)/2$ when $d_g(x,y)<d_g(x,\xi)/2$, we have
\begin{multline*}
\int_{\{d_g(x,y)<d_g(x,\xi)/2\}}\frac{d_g(x,y)^{2-n}}{(\delta^2+d_g(y,\xi)^2)^{n/2} }\,dv_g(y)\\
\leq Cd_g(x,\xi)^{-n}\int_{\{d_g(x,y)<d_g(x,\xi)/2\}} d_g(x,y)^{2-n} dv_g(y)\leq C d_g(x,\xi)^{2-n}.
\end{multline*}
As regards the second part of the integral, we have
\begin{multline*}
\int_{\{d_g(x,y)\geq d_g(x,\xi)/2\}}\frac{d_g(x,y)^{2-n}}{(\delta^2+d_g(y,\xi)^2)^{n/2} }\,dv_g(y)\\
\leq C d_g(x,\xi)^{2-n}\int_M (\delta^2+d_g(y,\xi)^2)^{-n/2} dv_g(y)\leq C d_g(x,\xi)^{2-n}\ln(1/\delta).
\end{multline*}
This yields \eqref{bnd:beta:d} when $d_g(x,\xi)>\delta$. Finally, we treat the case $d_g(x,\xi)\leq \delta$ in the same way as \eqref{bnd:beta:xi}. This ends the proof of Proposition~\ref{prop:est:beta:2}.
\endproof

It is a direct consequence of Proposition~\ref{prop:est:beta:2} that \eqref{hyp:B} is satisfied. Therefore Proposition~\ref{prop:c1} applies. It follows from \eqref{control:R:1}, \eqref{est:fin:1} and \eqref{est:c1} that
\begin{equation}\label{est:J:1}
J_{h}(W+\Phi)=J_h\left(W\right)+\bigO(\Vert R\Vert_{\frac{2n}{n+2}}^2)
\end{equation}
and, since $n\leq 5$, 
\begin{equation}\label{est:derJ:1}
\partial_{p}J_h(W+\Phi)=\partial_{p}J_h(W)+\bigO(\delta^{-1}\Vert R\Vert_{\frac{2n}{n+2}}(\Vert R\Vert_{\frac{2n}{n+2}}+\delta\Vert\partial_p R\Vert_{\frac{2n}{n+2}})),
\end{equation}
where $R=R_{\delta,\xi}$ is as in \eqref{def:R}. We prove the following estimates for $R$:

\begin{proposition}\label{prop:est:R}
We have 
\begin{equation}\label{est:R}
\Vert R\Vert_{\frac{2n}{n+2}}+\delta\Vert \partial_pR\Vert_{\frac{2n}{n+2}}\leq C\left\{\begin{aligned}
&\delta&&\hbox{if }n=3\\
&\delta^2\ln(1/\delta)&&\hbox{if }n=4\\
&D_{h,\xi}\delta^2 \ln(1/\delta)+\delta^2&&\hbox{if }n=5.
\end{aligned}\right.
\end{equation}
\end{proposition}

\proof[Proof of Proposition~\ref{prop:est:R}]
Note that since $n<6$, we have $\crit>3$. The definitions \eqref{def:beta:h}, \eqref{def:W:beta} and \eqref{def:B:beta} combined with \eqref{exp:U:delta} yield
\begin{align}\label{est:R:1}
R&=(\Delta_g+h)U +(\Delta_g+h)B-(U+B)_+^{\crit-1}=U^{\crit-1}-(U+B)_+^{\crit-1}\\
&=-(\crit-1)U^{\crit-2} B +\bigO(U^{\crit-3}B^2+|B|^{\crit-1}),\nonumber
\end{align}
where we have used that $U\geq 0$. Therefore,
$$\Vert R\Vert_{\frac{2n}{n+2}}\leq C \Vert U^{\crit-2} B\Vert_{\frac{2n}{n+2}} +\Vert |B|^{\crit-1}\Vert_{\frac{2n}{n+2}}.$$
Since $B=k_n\delta^{\frac{n-2}{2}}\beta$, the pointwise estimate \eqref{bnd:beta}, the estimate $U\leq C\tilde{U}$ and the estimates \eqref{step:comput} yield
$$\Vert R\Vert_{\frac{2n}{n+2}}\leq  C\left\{\begin{aligned}
&\delta&&\hbox{if }n=3\\
&\delta^2\ln(1/\delta)&&\hbox{if }n=4\\
&\delta^2&&\hbox{if }n=5.
\end{aligned}\right.$$
We now deal with the gradient term. We fix $i\in \{0,\dotsc,n\}$. We have
\begin{align*}
\partial_{p_i}R&=\partial_{p_i}(U^{\crit-1}-(U+B)_+^{\crit-1})\\
&=-(\crit-1)((U+B)_+^{\crit-2}(\partial_{p_i}U+\partial_{p_i}B)-U^{\crit-2}\partial_{p_i}U)\\
&=-(\crit-1)(((U+B)_+^{\crit-2}-U^{\crit-2})\partial_{p_i}U+(U+B)_+^{\crit-2}\partial_{p_i}B).
\end{align*}
Using that $\crit>3$ together with \eqref{hyp:B:p} and \eqref{est:lem:1}, we obtain
$$\delta|\partial_{p_i}R|\leq C \tilde{U}^{\crit-2}|B|+C\tilde{U}|B|^{\crit-2}+C\delta|\partial_{p_i}B|\tilde{U}^{\crit-2}.$$
Since $B=k_n\delta^{\frac{n-2}{2}}\beta$, using the estimates of $\beta$ and its derivatives in Proposition~\ref{prop:est:beta:2} and the estimates \eqref{step:comput}, long but easy computations yield
$$\delta\Vert \partial_{p_i}R\Vert_{\frac{2n}{n+2}}\leq  C\left\{\begin{aligned}
&\delta&&\hbox{if }n=3\\
&\delta^2\ln(1/\delta)&&\hbox{if }n=4\\
&D_{h,\xi}\delta^2 \ln(1/\delta)+\delta^2&&\hbox{if }n=5.
\end{aligned}\right.$$
Therefore, we obtain \eqref{est:R}. This ends the proof of Proposition~\ref{prop:est:R}.
\endproof
 
With \eqref{est:R}, the estimates \eqref{est:J:1} and \eqref{est:derJ:1} become
$$J_{h}(W+\Phi)=J_h\left(W\right)+O\left(\begin{aligned}
&\delta^2&&\hbox{if }n=3\\
&\delta^4\left(\ln(1/\delta)\right)^{2}&&\hbox{if }n=4\\
&\delta^4+D_{h,\xi}^2\delta^4 (\ln(1/\delta))^2&&\hbox{if }n=5
\end{aligned}\right).$$
and
$$\partial_{p_i}J_h(W+\Phi)=\partial_{p_i}J_h(W)+\bigO\left(\begin{aligned}
&\delta&&\hbox{if }n=3\\
&\delta^3\left(\ln(1/\delta)\right)^{2}&&\hbox{if }n=4\\
&\delta^3+D_{h,\xi}^2\delta^3 (\ln(1/\delta))^2&&\hbox{if }n=5
\end{aligned}\right).$$
We now estimate $J_h(W+\Phi)$:

\begin{proposition}\label{prop:est:J:W} 
We have
\begin{multline}
J_{h}(W+\Phi)= \frac{1}{n}\int_{\rn} U_{1,0}^{\crit}\,dx+\frac{1}{2}\varphi_h(\xi)\left\{\begin{aligned}
&0&&\hbox{if }n=3\\
&8\omega_{n-1}\delta^2\ln(1/\delta)&&\hbox{if }n=4 \\
&\delta^2 \int_{\rn}U_{1,0}^2\, dx&&\hbox{if }n=5
\end{aligned}\right\}\\
-\frac{k_n^2}{2} m_{h_0}(\xi_0)\delta^{n-2}+\smallo(\delta^{n-2})\label{est:J:1:2}
\end{multline}
as $\delta\to0$, $\xi\to\xi_0$ and $h\to h_0$ in $C^2(M)$.
\end{proposition}

\proof[Proof of Proposition~\ref{prop:est:J:W}]
We have
\begin{align}\label{dev:W:1}
J_{h}(W)&= \frac{1}{2}\int_M\left(|\nabla W|^2+hW^2\right)dv_g-\frac{1}{\crit}\int_M W_+^{\crit} dv_g\\
&= \frac{1}{2}\int_MRW\, dv_g+\left(\frac{1}{2}-\frac{1}{\crit}\right)\int_M W_+^{\crit} dv_g.\nonumber
\end{align}
Using that $U\ge0$, we obtain
\begin{equation}\label{eq:W:2}
W_+^{\crit}=(U+B)_+^{\crit}=U^{\crit}+\crit B U^{\crit-1}+\bigO\left(B^2U^{\crit-2}+|B|^{\crit}\right).
\end{equation}
Plugging \eqref{est:R:1} and \eqref{eq:W:2} into \eqref{dev:W:1}, and using \eqref{hyp:B:p} and \eqref{est:lem:1}, we obtain
\begin{multline*}
J_{h}(W)= \frac{1}{n}\int_M U^{\crit} dv_g- \frac{1}{2}\int_M BU^{\crit-1} dv_g\\
+\bigO\left(\int_M (\tilde{U}^{\crit-2}B^2+\tilde{U}|B|^{\crit-1}+|B|^{\crit})\,dv_g\right).
\end{multline*}
Since $B=k_n\delta^{\frac{n-2}{2}}\beta$, the pointwise estimate \eqref{bnd:beta}, the definition \eqref{def:U:delta} and \eqref{eq:3} yield
\begin{equation}
J_{h}(W)= \frac{1}{n}\int_{\rn} U_{1,0}^{\crit}\,dx- \frac{1}{2}\int_M BU^{\crit-1} dv_g+O\left(\begin{aligned}
&\delta^2&&\hbox{if }n=3\\
&\delta^4\left(\ln(1/\delta)\right)^{3}&&\hbox{if }n=4\\
&\delta^4 &&\hbox{if }n=5
\end{aligned}\right).
\label{dev:W:2}
\end{equation}
The definitions \eqref{def:beta:h} and \eqref{def:B:beta} of $\beta$ and $B$ yield
\begin{equation}\label{eqB}
\Delta_gB+hB=U^{\crit-1}-(\Delta_g U+hU)\hbox{ in }M.
\end{equation}
Therefore, we obtain
\begin{align*}
\int_M BU^{\crit-1} dv_g&=\int_M B(U^{\crit-1}-(\Delta_g U+hU))\,dv_g+\int_M B(\Delta_g U+hU)\,dv_g\\
&= \int_M \left(|\nabla B|^2+hB^2\right)dv_g+\int_M (\Delta_g B+hB)U\, dv_g\\
&=\int_M \left(|\nabla B|^2+hB^2\right)dv_g-\delta^{\frac{n-2}{2}}\int_M F_\delta(\xi,\cdot)U\, dv_g-\int_M \hat{h}_\xi U^2 dv_g.
\end{align*}
Since $B=k_n\delta^{\frac{n-2}{2}}\beta$, using \eqref{cv:beta:h} and \eqref{lim:F} together with Lebesgue's convergence theorem, we obtain
\begin{multline}\label{est:B:R}
\int_M BU^{\crit-1} dv_g=\delta^{n-2}k_n^2\bigg(\int_M \left(|\nabla \beta_{h,0,\xi}|^2+h\beta_{h,0,\xi}^2\right)dv_g\\
-\frac{1}{k_n}\int_M F_0(\xi,\cdot)\Gamma_\xi\, dv_g\bigg)-\int_M \hat{h}_\xi U^2 dv_g+\smallo(\delta^{n-2}).
\end{multline}
Since $U(x)^2\leq C \delta^{n-2}d_g(\xi,x)^{4-2n}$, letting $\xi\to \xi_0$ and $h\to h_0$ in $C^2(M)$, integration theory yields
$$\int_M \hat{h}_\xi U^2 dv_g=\delta^{n-2}k_n^2 \int_M \hat{(h_0)}_{\xi_0}\Gamma_{\xi_0}^2\,dv_g+\smallo(\delta^{n-2})\hbox{ when }n=3.$$
We now assume that $n\in\left\{4,5\right\}$. We write
\begin{multline*}
\int_M \hat{h}_\xi U^2 dv_g=\hat{h}_\xi (\xi)\int_M U^2 dv_g+\partial_{\xi_i} \hat{h}_\xi (\xi)\int_M x^i U^2 dv_g\\
+\int_M (\hat{h}_\xi-\hat{h}_\xi (\xi)-\partial_{\xi_i}\hat{h}_\xi (\xi)x^i )U^2 dv_g,
\end{multline*}
where the coordinates are taken with respect to the exponential chart at $\xi$. As one checks, there exists $C>0$ such that 
$$|\hat{h}_\xi-\hat{h}_\xi (\xi)-\partial_{\xi_i}\hat{h}_\xi (\xi)x^i|U^2\leq C \delta^{n-2}d_g(\xi,x)^{6-2n}$$
for all $x,\xi\in M$, $x\ne\xi$. Since $n<6$ and $\xi$ remains in a neighborhood of $\xi$ (so that the exponential chart remains nicely bounded), integration theory then yields
\begin{multline*}
\int_M (\hat{h}_\xi-\hat{h}_\xi (\xi)-\partial_{\xi_i}\hat{h}_\xi (\xi)x^i )U^2 dv_g\\
=\delta^{n-2}k_n^2 \int_M (\hat{h}_\xi-\hat{h}_\xi (\xi)-\partial_{\xi_i}\hat{h}_\xi (\xi)x^i )\Gamma_\xi^2 dv_g+\smallo(\delta^{n-2}).
\end{multline*}
Furthermore, letting $\xi\to \xi_0$, $h\to h_0$ and using $\varphi_{h_0}(\xi_0)=|\nabla\varphi_{h_0}(\xi_0)|=0$, we obtain
\begin{equation}\label{est:mass:1}
\int_M (\hat{h}_\xi-\hat{h}_\xi (\xi)-\partial_{\xi_i}\hat{h}_\xi (\xi)x^i )U^2 dv_g=\delta^{n-2}k_n^2 \int_M \hat{(h_0)}_{\xi_0}\Gamma_{\xi_0}^2\,dv_g+\smallo(\delta^{n-2}).
\end{equation}
Via the exponential chart, using the radial symmetry of $U$, we obtain 
\begin{align*}
\int_M x^i U^2 dv_g&=\sqrt{n(n-2)}^{n-2}\int_{B_{r_0}(0)} x^i \left(\frac{\delta}{\delta^2+|x|^2}\right)^{n-2}(1+\bigO(|x|))\,dx\\
&= \bigO\left(\int_{B_{r_0}(0)} |x|^2 \left(\frac{\delta}{\delta^2+|x|^2}\right)^{n-2} dx\right)=\bigO(\delta^{n-2})
\end{align*}
since $n<6$. It then follows from \eqref{est:U:L2}, \eqref{est:U:L2:2}, \eqref{est:U:L2:eps} and the above estimates that
\begin{multline*}
\int_M \hat{h}_\xi U^2 dv_g=\hat{h}_\xi (\xi)\left\{\begin{aligned}
&0&&\hbox{if }n=3\\
&8\omega_{n-1}\delta^2\ln(1/\delta)&&\hbox{if }n=4 \\
&\delta^2 \int_{\rn}U_{1,0}^2\, dx&&\hbox{if }n=5
\end{aligned}\right\}\\
+\delta^{n-2}k_n^2 \int_M \hat{(h_0)}_{\xi_0}\Gamma_{\xi_0}^2\,dv_g+\smallo(\delta^{n-2}).
\end{multline*}
Combining this estimate with \eqref{est:B:R}, we obtain
$$\int_M B U^{\crit-1} dv_g=-\hat{h}_\xi (\xi)\left\{\begin{aligned}
&0&&\hbox{if }n=3\\
&8\omega_{n-1}\delta^2\ln(1/\delta)&&\hbox{if }n=4 \\
&\delta^2 \int_{\rn}U_{1,0}^2\, dx&&\hbox{if }n=5
\end{aligned}\right\}+\delta^{n-2}k_n^2 I_{h_0,\xi_0}+\smallo(\delta^{n-2}),$$
where
\begin{align}\label{def:I}
I_{h_0,\xi_0}&:=\int_M \left(|\nabla \beta_{h_0,0,\xi_0}|^2+h_0\beta_{h_0,0,\xi_0}^2\right)dv_g\\
&\qquad-\frac{1}{k_n}\int_M F_0(\xi,\cdot)\Gamma_{\xi_0}\, dv_g-\int_M \hat{(h_0)}_{\xi_0}\Gamma_{\xi_0}^2\,dv_g.\nonumber
\end{align}
Integrating by parts and using the definition \eqref{def:beta:h}, we obtain
\begin{align*}
I_{h_0,\xi_0}&=\int_M\beta_{h_0,0,\xi_0}\left(\Delta_g\beta_{h_0,0,\xi_0}+h_0\beta_{h_0,0,\xi_0}\right)dv_g\\
&\quad-\int_M\Gamma_{\xi_0}\left(\frac{1}{k_n}F_0(\xi,\cdot)\Gamma_{\xi_0}+ \hat{(h_0)}_{\xi_0}\Gamma_{\xi_0}\right)dv_g\\
&=\int_M\left(\beta_{h_0,0,\xi_0}+\Gamma_{\xi_0}\right)\left(\Delta_g\beta_{h_0,0,\xi_0}+h_0\beta_{h_0,0,\xi_0}\right)dv_g\\
&=\int_M G_{h_0,\xi_0}(\Delta_g\beta_{h_0,0,\xi_0}+h_0\beta_{h_0,0,\xi_0})\,dv_g.
\end{align*}
We now use \eqref{pointwise:beta} at the point $\xi_0$, which makes sense since $\beta_{h_0,0,\xi_0}$ is continuous on $M$. This yields
\begin{equation}\label{exp:I:m}
I_{h_0,\xi_0}=\beta_{h_0,0,\xi_0}(\xi_0)=m_{h_0}(\xi_0).
\end{equation}
Putting these results together yields Proposition~\ref{prop:est:J:W}.
\endproof
 
We now estimate the derivatives of $J_h(W+\Phi)$:

\begin{proposition}\label{prop:est:der:J}
We have
\begin{multline}\label{prop:est:der:J1}
\partial_\delta J_{h}(W+\Phi)=\varphi_h(\xi)\left\{\begin{aligned}
&0&&\hbox{if }n=3 \\
&8\omega_{n-1}\delta\ln(1/\delta)&&\hbox{if }n=4 \\
&\delta \int_{\rn}U_{1,0}^2\, dx&&\hbox{if }n=5\\
\end{aligned}\right\}\\
-\frac{n-2}{2}k_n^2 m_{h_0}(\xi_0)\delta^{n-3}+\smallo(\delta^{n-3})
\end{multline}
and
\begin{multline}\label{prop:est:der:J2}
\partial_{\xi_i}J_h(W+\Phi)= \frac{1}{2}\partial_{\xi_i}\varphi_{h}(\xi)\left\{\begin{aligned}
&0&&\hbox{if }n=3 \\
&8\omega_{n-1}\delta^2\ln(1/\delta)&&\hbox{if }n=4\\
&\delta^2 \int_{\rn}U_{1,0}^2\, dx &&\hbox{if }n=5
\end{aligned}\right\}\\
+\bigO\left(\begin{aligned}
&\delta&&\hbox{if }n=3 \\
&\delta^2+D_{h,\xi}\delta^2\ln(1/\delta)&&\hbox{if }n=4\\
&\delta^3+D_{h,\xi}\delta^2&&\hbox{if }n=5
\end{aligned}\right)
\end{multline}
for all $i=1,\dotsc,n$, as $\delta\to0$, $\xi\to\xi_0$ and $h\to h_0$ in $C^2(M)$.
\end{proposition}

\proof[Proof of Proposition~\ref{prop:est:der:J}]
We fix $i\in \{0,\dotsc,n\}$. With \eqref{est:R:1}, \eqref{est:lem:1} and \eqref{hyp:B:p}, we obtain 
\begin{align*}
&\partial_{p_i}J_h(W)= J_h^\prime(W)[\partial_{p_i}W]=\int_{M}(\Delta_g W+hW-W_+^{\crit-1})\partial_{p_i}W\, dv_g=\int_{M}R \partial_{p_i}Wdv_g\\
&= -(\crit-1)\int_M U^{\crit-2} B  \partial_{p_i}W\, dv_g+\bigO\left(\int_M\left(U^{\crit-3}B^2+|B|^{\crit-1}\right)|\partial_{p_i}W|\, dv_g\right)\\
&= -(\crit-1)\int_M U^{\crit-2} B  \partial_{p_i}W\, dv_g+\bigO\left(\delta^{-1}\int_M\left(\tilde{U}^{\crit-2}B^2+\tilde{U}|B|^{\crit-1}\right)dv_g\right).
\end{align*}
As in the proof of \eqref{dev:W:2}, it follows from \eqref{bnd:beta} and \eqref{def:U:delta} that
\begin{equation}
\partial_{p_i}J_h(W)=-(\crit-1)\int_M U^{\crit-2} B  \partial_{p_i}W\, dv_g+\bigO(\delta^{-1})\left\{\begin{aligned}
&\delta^2&&\hbox{if }n=3\\
&\delta^4\left(\ln(1/\delta)\right)^{3}&&\hbox{if }n=4\\
&\delta^4 &&\hbox{if }n=5.
\end{aligned}\right.\label{eq:78}
\end{equation}
The estimates \eqref{bnd:beta:xi} and \eqref{bnd:beta:d} and the definition $B=k_n\delta^{\frac{n-2}{2}}\beta$ yield
$$\int_M U^{\crit-2} B  \partial_{p_i}B\, dv_g=\bigO(\delta^{-1})\left\{\begin{aligned}
&\delta^2&&\hbox{if }n=3\\
&\delta^4\left(\ln(1/\delta)\right)^{3}&&\hbox{if }n=4\\
&\delta^4+\epsilon_{i0}D_{h,\xi}\delta^3\ln(1/\delta)&&\hbox{if }n=5,
\end{aligned}\right.$$
where $\epsilon_{i0}$ is the Kronecker symbol. Since $W=U+B$, \eqref{eq:78} then yields
\begin{multline*}
\partial_{p_i}J_h(W)=-(\crit-1)\int_M U^{\crit-2} B  \partial_{p_i}U\, dv_g\\
+\bigO(\delta^{-1})\left\{\begin{aligned}
&\delta^2&&\hbox{if }n=3\\
&\delta^4\left(\ln(1/\delta)\right)^{3}&&\hbox{if }n=4\\
&\delta^4+\epsilon_{i0}D_{h,\xi}\delta^3\ln(1/\delta)&&\hbox{if }n=5,
\end{aligned}\right.
\end{multline*}
Differentiating \eqref{eqB}, we obtain
$$(\Delta_g+h)\partial_{p_i}B=(\crit-1)U^{\crit-2}\partial_{p_i}U-(\Delta_g+h)\partial_{p_i}U.$$
Multiplying by $B$ and integrating by parts, we then obtain
\begin{equation}\label{estB}
\int_M \partial_{p_i} B(\Delta_g+h)B\, dv_g=(\crit-1)\int_MU^{\crit-2}B \partial_{p_i}U\, dv_g-\int_M\partial_{p_i}U(\Delta_g+h)B\, dv_g.
\end{equation}
We begin with estimating the term in the left-hand-side of \eqref{estB}. Using that $B=k_n\delta^{\frac{n-2}{2}}\beta$, we obtain
\begin{multline*}
\int_M \partial_{p_i} B(\Delta_g+h)B\, dv_g=k_n^2 \delta^{n-2}\int_M \partial_{p_i} \beta(\Delta_g+h)\beta\, dv_g\\
+\epsilon_{i0}\frac{n-2}{2}k_n^2\delta^{n-2-1}\int_M\beta(\Delta_g\beta+h\beta)\,dv_g.
\end{multline*}
With \eqref{def:beta:h} and the pointwise estimates \eqref{bnd:beta:xi} and \eqref{bnd:beta:d}, we obtain
$$\left|\int_M \partial_{p_i} \beta(\Delta_g+h)\beta\, dv_g\right|\leq C\left\{\begin{aligned}
&1&&\hbox{if }n=3,4\\
&1+ D_{h,\xi}^2\ln(1/\delta)&&\hbox{if }n=5.
\end{aligned}\right.
$$
Therefore, we obtain
\begin{multline}\label{est:B:1}
\int_M \partial_{p_i} B(\Delta_g+h)B\, dv_g=\epsilon_{i0}\frac{n-2}{2}k_n^2\delta^{n-2-1}\int_M\beta(\Delta_g\beta+h\beta)\,dv_g\\
+\bigO\left(\begin{aligned}
&\delta^{n-2} &&\hbox{if }n=3,4\\
&\delta^{3}+ D_{h,\xi}^2\delta^{3}\ln(1/\delta)&&\hbox{if }n=5
\end{aligned}\right).
\end{multline}
We now deal with the second term in the right-hand-side of \eqref{est:B:1}. We first consider the case where $i\geq 1$, so that $\partial_{p_i}=\partial_{\xi_i}$. In this case, it follows from \eqref{eq:30} that $\partial_{\xi_i}U=-\partial_{x_i}U+\bigO(\tilde{U})$. Then, using \eqref{def:beta:h}, we obtain
$$-\int_M\partial_{\xi_i}U(\Delta_g+h)B\, dv_g=\int_M\partial_{x_i}U(\Delta_g+h)B\, dv_g+\bigO\left(\int_M{\tilde{U}}(\delta^{\frac{n-2}{2}}+|\hat{h}_\xi|\tilde{U})\,dv_g\right).$$
With \eqref{esthath}, we obtain
$$\int_M{\tilde{U}}(\delta^{\frac{n-2}{2}}+|\hat{h}_\xi|\tilde{U})\,dv_g\leq C\left\{\begin{aligned}
&\delta &&\hbox{if }n=3\\
&\delta^2+D_{h,\xi}\delta^2\ln(1/\delta)&&\hbox{if }n=4\\
&\delta^3+D_{h,\xi}\delta^2&&\hbox{if }n=5.
\end{aligned}\right.$$
With \eqref{def:beta:h} and the estimate $\partial_{x_i}U=\bigO(\delta^{\frac{n-2}{2}}d_g(x,\xi)^{1-n})$ (see the definition \eqref{def:peak}), we obtain
$$\int_M\partial_{x_i}U(\Delta_g+h)B\, dv_g=-\int_M\hat{h}_\xi U \partial_{x_i}U\, dv_g+\bigO(\delta^{n-2}).$$
Putting together the above estimates yields
\begin{multline*}
-(\crit-1)\int_M U^{\crit-2} B  \partial_{\xi_i}U\, dv_g= -\int_M\hat{h}_\xi U \partial_{x_i}U\, dv_g\\
+\bigO\left(\begin{aligned}
&\delta &&\hbox{if }n=3\\
&\delta^2+D_{h,\xi}\delta^2\ln(1/\delta)&&\hbox{if }n=4\\
&\delta^3+D_{h,\xi}\delta^2&&\hbox{if }n=5
\end{aligned}\right).
\end{multline*}
Using the explicit expression \eqref{def:U:delta} of $U$ together with the facts that $\Lambda_\xi(\xi)=1$, $\nabla \Lambda_\xi(\xi)=0$ and $|x|\partial_{x_i}U=\bigO(\tilde{U})$, we obtain 
\begin{multline*}
\int_M\hat{h}_\xi U \partial_{x_i}U\, dv_g=\int_{B_{r_0}(0)}\hat{h}_\xi(\exp_{\xi}^{g_\xi}(x))U_{\delta,0}\partial_{x_i}U_{\delta,0} (1+\bigO(|x|^2))\,dx\\
+\bigO\left(\int_{B_{r_0}(0)} |\hat{h}_\xi(\exp_{\xi}^{g_\xi}(x))| |x|\tilde{U}_{\delta,0}^2\, dx\right) +\bigO(\delta^{n-2})
\end{multline*}
With a Taylor expansion of $\hat{h}_\xi$, using the radial symmetry of $U_{\delta,0}$ and the explicit expressions given in \eqref{est:der:U:i}, we then obtain that there exists $c'_4,c'_5>0$ such that
$$\int_M\hat{h}_\xi U \partial_{x_i}U\, dv_g= -\partial_{\xi_i}\varphi_{h}(\xi)\left\{\begin{aligned}
&0&&\hbox{if }n=3\\
&c_4'\delta^2\ln(1/\delta)&&\hbox{if }n=4\\
&c_5'\delta^2 &&\hbox{if }n=5
\end{aligned}\right\}+\bigO(\delta^{n-2})$$
and then
\begin{multline*}
\partial_{\xi_i}J_h(W)= \partial_{\xi_i}\varphi_{h}(\xi)\left\{\begin{aligned}
&0&&\hbox{if }n=3\\
&c'_4\delta^2\ln(1/\delta)&&\hbox{if }n=4\\
&c'_5\delta^2 &&\hbox{if }n=5
\end{aligned}\right\}\\
+\bigO\left(\begin{aligned}
&\delta&&\hbox{if }n=3\\
&\delta^2+D_{h,\xi}\delta^2\ln(1/\delta)&&\hbox{if }n=4\\
&\delta^3+D_{h,\xi}\delta^2&&\hbox{if }n=5.
\end{aligned}\right).
\end{multline*}
We now consider the case where $i=0$, so that $\partial_{p_i}=\partial_{p_0}=\partial_\delta$. In this case, we have
\begin{align*}
\int_M\partial_\delta U(\Delta_g+h)B\, dv_g&=-\int_{M}\hat{h}_\xi U\partial_\delta U\, dv_g-\delta^{\frac{n-2}{2}}\int_MF \partial_\delta U\, dv_g \\
&=-\int_{M}(\hat{h}_\xi(\xi)+x^i\partial_{\xi_i}\hat{h}_{\xi}(\xi)) U\partial_\delta U\, dv_g\\
&\quad-\int_M(\delta^{\frac{n-2}{2}}F+(\hat{h}_\xi-\hat{h}_\xi(\xi)-x^i\partial_{\xi i}\hat{h}_{\xi}(\xi))U) \partial_\delta U\, dv_g,
\end{align*}
where the coordinates are taken with respect to the exponential chart at $\xi$. With \eqref{der:U}, \eqref{def:Z:0} and \eqref{est:der:U:delta}, arguing as in the proof of \eqref{est:mass:1}, we obtain
\begin{multline*}
\delta^{\frac{n-2}{2}}\int_M(F+(\hat{h}_\xi-\hat{h}_\xi(\xi)-x^i\partial_i\hat{h}_{\xi}(\xi))\delta^{-\frac{n-2}{2}}U) \partial_\delta U\, dv_g\\
=\frac{n-2}{2}k_n^2\delta^{-1}\delta^{n-2}\int_M\left(\frac{F_{\xi,0}}{k_n}+\hat{h}_{\xi_0}\Gamma_{\xi_0}^{h}\right)\Gamma_{\xi_0}^h\, dv_g+\smallo(\delta^{-1}\delta^{n-2}).
\end{multline*}
Using \eqref{est:U:L2} and arguing as in the estimate of \eqref{est:A:U:L2}, we obtain that there exist $c''_4,c''_5>0$ such that
\begin{multline*}
\int_{M}(\hat{h}_\xi(\xi)+x^i\partial_{\xi_i}\hat{h}_{\xi}(\xi)) U\partial_\delta U\, dv_g=\frac{\hat{h}_\xi(\xi)}{\delta}\int_{B_0(r_0)}U_{\delta,0}Z_{\delta,0}\, dx+\smallo(\delta^{-1}\delta^{n-2})\\
= \frac{\hat{h}_\xi(\xi)}{\delta}\left\{\begin{aligned}
&0&&\hbox{if }n=3\\
&c''_4\delta^2\ln(1/\delta)&&\hbox{if }n=4\\
&c''_5\delta^2&&\hbox{if }n=5
\end{aligned}\right\}+\smallo(\delta^{-1}\delta^{n-2}).
\end{multline*}
Putting these estimates together yields
\begin{multline*}
-(\crit-1)\int_M U^{\crit-2} B  \partial_{p_i}U\, dv_g=\frac{\hat{h}_\xi(\xi)}{\delta}\left\{\begin{aligned}
&0&&\hbox{if }n=3\\
&c''_4\delta^2\ln(1/\delta)&&\hbox{if }n=4\\
&c''_5\delta^2&&\hbox{if }n=5
\end{aligned}\right\}\\
-\frac{n-2}{2}k_n^2I_{h_0,\xi_0}\delta^{-1}\delta^{n-2}+\smallo(\delta^{-1}\delta^{n-2}),
\end{multline*}
where $I_{h_0,\xi_0}$ is as in \eqref{def:I}. Since $I_{h_0,\xi_0}=m_{h_0}(\xi_0)$ (see \eqref{exp:I:m}), we obtain \eqref{prop:est:der:J1} and \eqref{prop:est:der:J2} up to the value of the constants. These values then follow from the above estimates together with Proposition~\ref{prop:est:J:W}.
\endproof

Theorem~\ref{th:main} for $n\in\left\{4,5\right\}$ will be proved in Section~\ref{sec:pf:1}.

\section{Proof of Theorem~\ref{th:main}}\label{sec:pf:1}

We let $h_0,f\in C^p(M)$, $p\geq 2$, and $\xi_0\in M$ satisfy the assumptions in  Theorem~\ref{th:main}. For small $\eps>0$ and $\tau\in\rn$, we define
\begin{equation}\label{hxieps}
h_\eps:=h_0+\eps f\hbox{ and }\xi_\eps(\tau):=\exp_{\xi_0}^{g_{\xi_0}}(\sqrt{\eps}\tau).
\end{equation}
We fix $R>0$ and $0<a<b$ to be chosen later. 
 
\subsection{Proof of Theorem~\ref{th:main} for $n\geq 6$} 

In this case, we let $(\delta_\eps)_{\eps>0}>0$ be such that $\delta_\eps\to 0$ as $\eps\to 0$. We define
\begin{equation}\label{deltaFeps}
\delta_\eps(t):=\delta_\eps t\hbox{ and }F_\eps(t,\tau):=J_{h_\eps}(U_{\delta_\eps(t),\xi_\eps(\tau)}+\Phi_{h_\eps,0,\delta_\eps(t),\xi_\eps(\tau)})
\end{equation}
for all $\tau\in\rn$ such that $|\tau|<R$ and $t>0$ such that $a<t<b$. Using the assumption $\varphi_{h_0}(\xi_0)=|\nabla\varphi_{h_0}(\xi_0)|=0$, we obtain
$$\varphi_{h_\eps}(\xi_{\eps}(\tau))=\frac{1}{2}\nabla^2\varphi_{h_0}(\xi_0)[\tau,\tau]\eps+f(\xi_0)\eps+\smallo(\eps)$$
and
$$\nabla\varphi_{h_\eps}(\xi_{\eps}(\tau))=\nabla^2\varphi_{h_0}(\xi_0)[\tau,\cdot]\sqrt{\eps}+\smallo(\sqrt{\eps})$$
as $\eps\to 0$ uniformly with respect to $|\tau|<R$. We distinguish two cases:

\smallskip\noindent
{\bf Case $n\geq 7$.} In this case, we set $\delta_\eps:=\sqrt{\eps}$. It follows from \eqref{est:J:W:2} that
\begin{equation}\label{lim:c0}
\lim_{\eps\to 0}\frac{F_\eps(t,\tau)-\frac{1}{n}\int_{\rn}U_{1,0}^{\crit}\,dx}{\eps^2}=E_0(t,\xi)\hbox{ in }C^0_{\loc}((0,\infty)\times\rn),
\end{equation}
where
$$E_0(t,\tau):=C_n\left(\frac{1}{2}\nabla^2\varphi_{h_0}(\xi_0)[\tau,\tau]+f(\xi_0)\right)t^2-D_n K_{h_0}(\xi_0)t^4,$$
with
\begin{equation}\label{def:C}
C_n:=\frac{1}{2}\int_{\rn}U_{1,0}^2\, dx\hbox{ and }D_n:=\frac{1}{4n}\int_{\rn}|x|^2U_{1,0}^2\, dx.
\end{equation}
Furthermore, we have
$$\partial_t F_\eps(t,\tau)=\sqrt{\eps}\left(\partial_\delta J_{h_\eps}(U_{(\delta_\eps(t),\xi_\eps(\tau)}+\Phi_{(\delta_\eps(t),\xi_\eps(\tau)})\right)$$
and
$$\partial_{\tau_i} F_\eps(t,\tau)=\sqrt{\eps}\left(\partial_{\xi_i} J_{h_\eps}(U_{(\delta_\eps(t),\xi_\eps(\tau)}+\Phi_{(\delta_\eps(t),\xi_\eps(\tau)})\right).$$
Therefore, it follows from \eqref{der:delta:2} and \eqref{der:xi:2} that the limit in \eqref{lim:c0} actually holds in $C^1_{\loc}((0,\infty)\times\rn)$. Assuming that $f(\xi_0)\times K_{h_0}(\xi_0)>0$, we can define
$$t_0:=\sqrt{\frac{C_nf(\xi_0)}{2D_n K_{h_0}(\xi_0)}}\,.$$
As one checks, $(t_0,0)$ is a critical point of $E_0$. In addition, the Hessian matrix at the critical point $(t_0,0)$ is
$$\nabla^2E_0(t_0,0)=\left(\begin{array}{cc}
-8t_0^2D_n K_{h_0}(\xi_0) & 0\\
0 & t_0^2C_n\nabla^2\varphi_{h_0}(\xi_0)
\end{array}\right).$$
Therefore, if $\xi_0$ is a nondegenerate critical point of $\varphi_{h_0}$, then $(t_0, 0)$ is a nondegenerate critical point of $E_0$. With the convergence in $C^1_{\loc}((0,\infty)\times\rn)$, we then obtain that there exists a critical point $(t_\eps,\tau_\eps)$ of $F_\eps$ such that $(t_\eps,\tau_\eps)\to(t_0,0)$ as $\eps\to 0$. It then follows from \eqref{cns} that 
$$u_\eps:=U_{\delta_\eps(t_\eps),\xi_\eps(\tau_\eps)}+\Phi_{h_\eps,0,\delta_\eps(t_\eps),\xi_\eps(\tau_\eps)}$$ 
is a solution to \eqref{eq:ue}. As one checks, $\ue\rightharpoonup 0$ weakly in $L^{\crit}(M)$ and $(\ue)_\eps$ blows up with one bubble at $\xi_0$. This proves Theorem~\ref{th:main} for $n\geq 7$.

\smallskip\noindent 
{\bf Case $n=6$.} In this case, we let $\delta_\eps>0$ be such that
\begin{equation}\label{delta6}
\delta_\eps^2\ln(1/\delta_\eps)=\eps.
\end{equation}
As one checks, $\delta_\eps\to0$ as $\eps\to 0$. As in the previous case, we obtain
$$\lim_{\eps\to 0}\frac{F_\eps(t,\tau)-\frac{1}{n}\int_{\rn}U_{1,0}^{\crit}\,dx}{\eps\delta_\eps^2}=E_0(t,\xi)\hbox{ in }C^1_{\loc}((0,\infty)\times\rn),$$
where
$$E_0(t,\tau):=C_6\left(\frac{1}{2}\nabla^2\varphi_{h_0}(\xi_0)[\tau,\tau]+f(\xi_0)\right) t^2-24^2\omega_5K_{h_0}(\xi_0) t^4$$
for all $t>0$ and $\tau\in\rn$. As in the previous case, $E_0$ has a nondegenerate critical point $(\tilde{t}_0,0)$, which yields the existence of a critical point of $F_\eps$ and, therefore, a blowing-up solution to \eqref{eq:ue} satisfying the desired conditions. This proves Theorem~\ref{th:main} for $n=6$.

\subsection{Proof of Theorem~\ref{th:main} for $n\in\left\{4,5\right\}$.} 

When $n\in\left\{4,5\right\}$, we define
$$F_\eps(t,\tau):=J_{h_\eps}(U_{\delta_\eps(t),\xi_\eps(\tau)}+B_{h_\eps,\delta_\eps(t),\xi_\eps(\tau)}+\Phi_{h_\eps,0,\delta_\eps(t),\xi_\eps(\tau)}),$$
where $\delta_\eps(t)$ will be chosen differently depending on the dimension.

\smallskip\noindent
{\bf Case $n=5$.} In this case, we set $\delta_\eps(t):=t\eps$. It follows from \eqref{est:J:1:2} that
$$\lim_{\eps\to 0}\frac{F_\eps(t,\tau)-\frac{1}{n}\int_{\rn}U_{1,0}^{\crit}\, dx}{\eps^3}=E_0(t,\xi)\hbox{ in }C^0_{\loc}((0,\infty)\times\rn),$$
where
$$E_0(t,\tau):=C_5\left(\frac{1}{2}\nabla^2\varphi_{h_0}(\xi_0)(\tau,\tau)+f(\xi_0)\right) t^2 -\frac{k_5^2}{2} m_{h_0}(\xi_0)t^3.$$
Furthermore, it follows from the $C^1-$estimates of Proposition~\ref{prop:est:der:J} that the convergence holds in $C^1_{\loc}((0,\infty)\times\rn)$. Assuming that $f(\xi_0)\times m_{h_0}(\xi_0)>0$, we then define
$$t_0:=\frac{4C_5f(\xi_0)}{(n-2)k_5^2m_{h_0}(\xi_0)}.$$
As in the previous cases, we obtain that $(t_0,0)$ is a nondegenerate critical point of $E_0$, which yields the existence of a critical point for $F_\eps$ and, therefore, a blowing-up solution to \eqref{eq:ue} satisfying the desired conditions. This proves Theorem~\ref{th:main} for $n=5$.

\smallskip\noindent
{\bf Case $n=4$.} In this case, we set $\delta_\eps(t):=e^{-t/\eps}$. It follows from the $C^1-$estimates of Proposition~\ref{prop:est:der:J} that
$$\lim_{\eps\to 0}\left(-\eps\delta_\eps(t)^{-2}\partial_t F_\eps(t,\tau),\delta_\eps(t)^{-2}\partial_{\tau} F_\eps(t,\tau)\right) =(\psi_0(t,\tau),\psi_1(t,\tau)),$$
in $C^0_{\loc}((0,\infty)\times\rn)$, where
$$\psi_0(t,\tau):=C_4\left(\frac{1}{2}\nabla^2\varphi_{h_0}(\xi_0)(\tau,\tau)+f(\xi_0)\right) t-\frac{n-2}{2}k_n^2m_{h_0}(\xi_0)$$
and
$$\psi_1(t,\tau):=\frac{1}{2}C_4\nabla^2\varphi_{h_0}(\xi_0)[\tau,\cdot]t.$$
As one checks, since $\xi_0$ is a nondegenerate critical point of $\varphi_{h_0}$, the function $\psi$ has a unique zero point in $(0,\infty)\times\rn$ which is of the form $(t_0, 0)$ for some $t_0>0$. Furthermore, the nondegeneracy implies that the Jacobian determinant of $\psi$ at $(t_0,0)$ is nonzero and, therefore, the degree of $\psi$ at $0$ is well-defined and nonzero. The invariance of the degree under uniform convergence then yields the existence of a critical point $(t_\eps,\tau_\eps)$ of $F_\eps$ such that $(t_\eps,\tau_\eps)\to(t_0,0)$ as $\eps\to 0$. It then follows from \eqref{cns} that 
$$u_\eps:=U_{\delta_\eps(t),\xi_\eps(\tau)}+B_{h_\eps,\delta_\eps(t),\xi_\eps(\tau)}+\Phi_{h_\eps,0,\delta_\eps(t),\xi_\eps(\tau)}$$
is a critical point of $J_{\he}$ that blows up at $\xi_0$ and converges weakly to 0 in $L^{\crit}(M)$. This proves Theorem~\ref{th:main} for $n=4$.\qed
\endproof

\section{Proof of Theorem~\ref{th:main:u0}}\label{sec:pf:2}

We let $h_0,f\in C^p(M)$, $p\geq 2$, $u_0\in C^2(M)$ and $\xi_0\in M$ satisfy the assumptions in Theorem~\ref{th:main:u0}. We let $h_\eps$ be as in \eqref{eq:ue}. We let $\xi_\epsilon(\tau)$ and $\delta_\epsilon(t)$ be as in \eqref{hxieps} and \eqref{deltaFeps}. Since $u_0$ is nondegenerate, the implicit function theorem yields the existence of $\epsilon'_0\in(0,\epsilon_0)$ and $(u_{0,\eps})_{0<\eps<\eps'_0}\in C^2(M)$ such that
\begin{equation}\label{eq:u0e}
\Delta_gu_{0,\eps}+h_\eps u_{0,\eps}=u_{0,\eps}^{\crit-1},\ u_{0,\eps}>0\hbox{ in }M
\end{equation}
and $(u_{0,\eps})_\eps$ is smooth with respect to $\eps$, which implies in particular that
$$\Vert u_{0,\eps}-u_0\Vert_{C^2}\leq C\eps.$$
We fix $0<a<b$ and $R>0$ to be chosen later. We define
$$F_\eps(t,\tau):=J_{h_\eps}(u_{0,\eps}+U_{\delta_\eps(t), \xi_\eps(t)}+\Phi_{h_\eps,u_{0,\eps},\delta_\eps(t), \xi_\eps(t)})$$
for all $\tau\in\rn$ such that $|\tau|<R$ and $t>0$ such that $a<t<b$.
With \eqref{est:c0phi}, we obtain that for $n\geq 7$,
\begin{multline*}
F_\eps(t,\tau)= J_{h_\eps}(u_{0,\eps})+\frac{1}{n}\int_{\rn}U_{1,0}^{\crit}\, dx+C_n\left(\frac{1}{2}\nabla^2\varphi_{h_0}(\xi_0)(\tau,\tau)+f(\xi_0)\right)t^2\eps \delta_\eps^2\\
 +\smallo(\eps\delta_\eps^2)-D_nK_{h_0}(\xi_0)t^4\delta_\eps^4+\smallo(\delta_\eps^4)-B_nu_0(\xi_0)t^{\frac{n-2}{2}}\delta_\eps^{\frac{n-2}{2}}+\smallo(\delta_\eps^{\frac{n-2}{2}})
\end{multline*}
as $\eps\to 0$ uniformly with respect to $a<t<b$ and $|\tau|<R$, where $C_n$ and $D_n$ are as in \eqref{def:C} and 
$$B_n:=\int_{\rn}U_{1,0}^{\crit-1} dx.$$ 
We distinguish three cases:

\smallskip\noindent
{\bf Case $7\leq n\leq 10$, that is $n\geq 7$ and $\frac{n-2}{2}\leq 4$.} In this case, we set $\delta_\eps:=\eps^{\frac{2}{n-6}}$, so that
$$\eps \delta_\eps^2=\delta_\eps^{\frac{n-2}{2}}.$$
We then obtain 
\begin{equation}\label{cv:F}
\lim_{\eps\to 0}\frac{F_\eps(t,\tau)-A_\eps}{\eps\delta_\eps^2}=E_0(t,\tau)
\end{equation}
uniformly with respect to $a<t<b$ and $|\tau|<R$, where
\begin{equation}\label{def:Ae}
A_\eps:=J_{h_\eps}(u_{0,\eps})+\frac{1}{n}\int_{\rn}U_{1,0}^{\crit}\, dx
\end{equation}
and
\begin{align*}
E_0(t,\tau)&:=C_n\left(\frac{1}{2}\nabla^2\varphi_{h_0}(\xi_0)(\tau,\tau)+f(\xi_0)\right)t^2\\
&\qquad-\left(B_nu_0(\xi_0)+{\bf 1}_{n=10}D_nK_{h_0}(\xi_0)\right)t^{\frac{n-2}{2}}.\nonumber
\end{align*}
Moreover, the estimates \eqref{der:J:u:0} and \eqref{der:J:u:i} yield the convergence \eqref{cv:F} in $C^1_{\loc}((0,\infty)\times\rn)$.
Straightforward changes of variable yield 
$$\frac{B_{10}}{D_{10}}=40\frac{\int_{\rr^{10}}U_{1,0}^{3/2}dx}{\int_{\rr^{10}}|x|^2U_{1,0}^{3/2}dx}=40\frac{\int_0^\infty\frac{r^9dr}{(1+r^2)^6}}{\int_0^\infty\frac{r^{11}dr}{(1+r^2)^8}}=40\frac{\int_0^\infty\frac{s^4ds}{(1+s)^6}}{\int_0^\infty\frac{s^5ds}{(1+s)^8}}\,.$$
Integrating by parts, we then obtain
\begin{multline*}
\frac{B_{10}}{D_{10}}=\frac{40\times\frac{4}{5}\times\frac{3}{4}\times\frac{2}{3}\times\frac{1}{2}\int_0^\infty\frac{ds}{(1+s)^2}}{\frac{5}{7}\times\frac{4}{6}\times\frac{3}{5}\times\frac{2}{4}\times\frac{1}{3}\int_0^\infty\frac{ds}{(1+s)^3}}=\frac{40\times6\times7\int_0^\infty\frac{ds}{(1+s)^2}}{5\int_0^\infty\frac{ds}{(1+s)^3}}=\frac{40\times6\times7\times2}{5}\\
=672.
\end{multline*}
The assumption $K_{h_0,u_0}(\xi_0)\ne0$ then gives $B_nu_0(\xi_0)+{\bf 1}_{n=10}D_nK_{h_0}(\xi_0)\neq 0$ with same sign as $f(\xi_0)$. As in the proof of Theorem~\ref{th:main} for $n\geq 7$, we obtain that $E_0$ has a unique critical point in $(0,\infty)\times\rn$, say $(t_0,0)$, and this critical point is nondegenerate. Mimicking again the proof of Theorem~\ref{th:main} for $n\geq 7$, we obtain the existence of a critical point $(t_\eps,\tau_\eps)$ of $F_\eps$ such that $(t_\eps,\tau_\eps)\to(t_0,0)$ as $\eps\to 0$. It then follows that
$$\ue:=u_{0,\eps}+U_{\delta_\eps(t_\eps), \xi_\eps(\tau)}+\Phi_{h_\eps,u_{0,\eps},\delta_\eps(t_\eps), \xi_\eps(\tau)}$$
is a solution to \eqref{eq:ue}. As one checks, $\ue\rightharpoonup 0$ weakly in $L^{\crit}(M)$ and $(\ue)_\eps$ blows up with one bubble at $\xi_0$. This proves Theorem~\ref{th:main:u0} for $7\le n\le10$.

\smallskip\noindent
{\bf Case $4<\frac{n-2}{4}$, that is $n\ge11$.} In this case, we set $\delta_\eps:=\sqrt{\eps}$, so that
$$\eps \delta_\eps^2=\delta_\eps^4\hbox{ and }\delta_\eps^{\frac{n-2}{2}}=\smallo(\delta_\eps^4)\hbox{ as }\eps\to 0.$$
We then obtain 
$$\lim_{\eps\to 0}\frac{F_\eps(t,\tau)-A_\eps}{\eps\delta_\eps^2}=E_0(t,\tau)\hbox{ in }C^0_{\loc}((0,\infty)\times\rn),$$
where $A_\eps$ is as in \eqref{def:Ae} and
$$E_0(t,\tau):=C_n\left(\frac{1}{2}\nabla^2\varphi_{h_0}(\xi_0)(\tau,\tau)+f(\xi_0)\right) t^2-D_nK_{h_0}(\xi_0)t^4.$$
As in the previous case, we obtain that the convergence holds in $C^1_{\loc}((0,\infty)\times\rn)$ and $E_0$ has a nondegenerate critical point in $(0,\infty)\times\rn$, which yields the existence of a blowing-up solution $(\ue)_\eps$ to \eqref{eq:ue} satisfying the desired conditions. This proves Theorem~\ref{th:main:u0} for $n\geq 11$.

\smallskip\noindent
{\bf Case n=6.} Note that in this case, we have $\crit-1=2$. Differentiating \eqref{eq:u0e} with respect to $\epsilon$ at 0, we obtain
$$(\Delta_g+h_0-2u_0)(\partial_\eps u_{0,\eps})_{|0}+fu_0=0\hbox{ in }M.$$
Using that $u_0$ is nondegenerate, we then obtain
$$(\partial_\eps u_{0,\eps})_{|0}=-(\Delta_g+h_0-2u_0)^{-1}(fu_0).$$ 
It follows that
$$\varphi_{h_\epsilon,u_\epsilon}=h_\epsilon-2u_{0,\epsilon}-c_n\Scal_g=\varphi_{h_0,u_0}+\tilde{f}\eps+\smallo(\eps)\hbox{ as }\eps\to0,$$
where 
$$\tilde{f}:=f+2(\Delta_g+h_0-2u_0)^{-1}(fu_0).$$
We let $\delta_\eps>0$ be as in \eqref{delta6}. With \eqref{co:6}, we then obtain
$$\lim_{\eps\to 0}\frac{F_\eps(t,\tau)-A_\eps}{\eps\delta_\eps^2}=E_0(t,\tau)\hbox{ in }C^0_{\loc}((0,\infty)\times\rn),$$
where $A_\eps$ is as in \eqref{def:Ae} and 
$$E_0(t,\tau):=C_6\left(\frac{1}{2}\nabla^2\varphi_{h_0}(\xi_0)(\tau,\tau)+\tilde{f}(\xi_0)\right) t^2-24^2\omega_5 K_{h_0,u_0}(\xi_0)t^4.$$
As in the previous case, using \eqref{der:0:6} and \eqref{der:i:6}, we obtain that the convergence holds in $C^1_{\loc}((0,\infty)\times\rn)$. Furthermore, using \eqref{condition6}, we obtain that $E_0$ has a nondegenerate critical point in $(0,\infty)\times\rn$ and, therefore, that there exists a blowing-up solution to \eqref{eq:ue} satisfying the desired conditions. This proves Theorem~\ref{th:main:u0} for $n=6$.\qed
\endproof

\section{Proof of Theorem~\ref{th:2}}\label{sec:pf:3}

We let $h_0\in C^p(M)$, $1\le p\le\infty$, and $\xi_0\in M$ be such that $\Delta_g+h_0$ is coercive and $\varphi_{h_0}(\xi_0)=|\nabla\varphi_{h_0}(\xi_0)|=0$. In the case where $p=1$, a standard regularization argument give the existence of $(\hat{h}_k)_{k\in\nn}\in C^2(M)$ such that $\hat{h}_k\to h_0$ in $C^1(M)$. In the case where $p\ge2$, we set $\hat{h}_k=h_0$. We then define 
$$\tilde{h}_k:=\hat{h}_k+f_k,\hbox{ where }f_k(x):=\chi(x)((h_0-\hat{h}_k)(\xi_0)+\langle\nabla(h_0-\hat{h}_k)(\xi_0),x\rangle+|x|^2/k),$$
where $\chi$ is a smooth cutoff function around 0 and the coordinates are taken with respect to the exponential chart at $\xi_0$. As one checks, we then have that $\tilde{h}_k\to h_0$ in $C^p(M)$, $\varphi_{\tilde{h}_k}(\xi_0)=\varphi_{h_0}(\xi_0)=0$, $|\nabla\varphi_{\tilde{h}_k}(\xi_0)|=|\nabla\varphi_{h_0}(\xi_0)|=0$ and for large $k$, $\xi_0$ is a nondegenerate critical point of $\varphi_{\tilde{h}_k}$.

\smallskip
Assume first that $n\in\left\{4,5\right\}$. Then the mass of $\tilde{h}_k$ is defined at $\xi_0$. If $m_{h_0}(\xi_0)\ne0$, then by continuity, $m_{\tilde{h}_k}(\xi_0)\neq 0$ for large $k$, with same sign as $m_{h_0}(\xi_0)$. If $m_{h_0}(\xi_0)=0$, then it follows from Proposition~\ref{prop:diff:mass} that $m_{\tilde{h}_k}(\xi_0)<0$ for large $k$. Therefore, in all cases, we have that $m_{\tilde{h}_k}(\xi_0)\neq 0$ for large $k$, with a sign independent of $k$.

\smallskip
Assume now that $n\geq 6$. With a similar argument, we obtain that, for large $k$, $K_{\tilde{h}_k}(\xi_0)\neq 0$ with a sign independent of $k$, where $K_{\tilde{h}_k}(\xi_0)$ is as in \eqref{Kh}.

\smallskip
In all cases, we can now fix $f_0\in C^\infty(M)$ such that $f_0(\xi_0)\times K_{\tilde{h}_k}(\xi_0)>0$ for large $k$. It then follows from Theorem~\ref{th:main} that  there exist $\eps_k>0$ and a family $(\tilde{u}_{k,\eps})_{0<\eps<\eps_k}$ of solutions to the equation
$$\Delta_g \tilde{u}_{k,\eps}+(\tilde{h}_k+\eps f_0)\tilde{u}_{k,\eps}=\tilde{u}_{k,\eps}^{\crit-1},\ \tilde{u}_{k,\eps}>0\hbox{ in }M$$
such that $\tilde{u}_{k,\eps}\rightharpoonup 0$ weakly in $L^{\crit}(M)$ and $(\tilde{u}_{k,\eps})_{\eps}$ blows up with one bubble at $\xi_0$ as $\eps\to 0$. Therefore, we obtain that for every $k\in\nn$, there exists $\eps'_k>0$ such that
$$0<\eps'_k<\min(1/k, \eps_k),\ \Vert \tilde{u}_{k,\eps'_k}\Vert_2<\frac{1}{k},\ \left|\int_{M}|\tilde{u}_{k,\eps'_k}|^{\crit} dv_g-\int_{\rn}U_{1,0}^{\crit}\, dx\right|<\frac{1}{k}$$
and
$$\int_{M\backslash B_{1/k}(\xi_0)}|\tilde{u}_{k,\eps'_k}|^{\crit} dv_g<\frac{1}{k}\,.$$
We then define $u_{k}:=\tilde{u}_{k,\eps'_k}$, so that 
$$\Delta_g u_k+h_k u_k=u_k^{\crit-1}\hbox{ in }M,\hbox{ where }h_k:=\tilde{h}_k+\eps'_kf_0=h_0+f_k+\eps'_kf_0.$$
As one checks, $u_k\rightharpoonup 0$ weakly in $L^{\crit}(M)$ and $(u_k)_k$ blows up with one bubble at $\xi_0$ as $k\to\infty$. This proves Theorem~\ref{th:2}.\qed
\endproof

\section{Proof of Theorem~\ref{th:3}}\label{sec:pf:4}

We let $h_0\in C^p(M)$, $1\le p\le\infty$, $u_0\in C^2(M)$ and $\xi_0\in M$ be such that $\Delta_g+h_0$ is coercive, $u_0$ is a solution of \eqref{eq} and the condition \eqref{phi0th3} is satisfied. We begin with proving the following:

\begin{lem}\label{lemu0}
There exists a neighborhood $\Omega_0$ of $\xi_0$ and sequences $(\tilde{h}_k)_{k\in\nn}\in C^p(M)$ and $(\tilde{u}_k)_{k\in\nn}\in  C^2(M)$ such that $\tilde{h}_k\to h_0$ in $C^p(M)$, $\tilde{u}_k\to u_0$ in $C^2(M)$, $\tilde{h}_k\equiv h_0$ and $\tilde{u}_k\equiv u_0$ in $\Omega_0$ and $\tilde{u}_k$ is a nondegenerate solution of 
\begin{equation}\label{eqtilduk}
\Delta_g\tilde{u}_k+\tilde{h}_k\tilde{u}_k=\tilde{u}_k^{\crit-1},\ \tilde{u}_k>0\hbox{ in }M\hbox{ for all }k\in\nn.
\end{equation}
\end{lem}

\proof[Proof of Lemma~\ref{lemu0}]
For all $v\in C^{p+2}(M)$ such that $v>-u_0$, we define 
$$u(v):=u_0+v\hbox{ and }h(v):=u(v)^{\crit-2}-\frac{u_0^{\crit-1}-h_0u_0+\Delta_gv}{u(v)}=u(v)^{\crit-2}-\frac{\Delta_gu(v)}{u(v)},$$
so that 
\begin{equation}\label{eqv}
\Delta_gu(v)+h(v)u(v)=u(v)^{\crit-1}\hbox{ in }M.
\end{equation}
By elliptic regularity, we have $u_0\in C^{p+1}(M)$. Since moreover $h_0\in C^p(M)$ and $v\in C^{p+2}(M)$, we obtain that $u(v)\in C^{p+1}(M)$ and $h(v)\in C^p(M)$. Furthermore, we have that $h(v)\to h_0$ in $C^p(M)$ and $u(v)\to u_0$ in $C^2(M)$ as $v\to 0$ in $C^{p+2}(M)$. As is easily seen, to prove the lemma, it suffices to show that there exists a neighborhood $\Omega_0$ of $\xi_0$ and a sequence $(v_k)_{k\in\nn}\in C^{p+2}(M)$  such that $v_k\to 0$ in $C^{p+2}(M)$, $v_k\equiv0$ in $\Omega_0$ and $u(v_k)$ is a nondegenerate solution of \eqref{eqv}. Assume by contradiction that this is not true, that is for every neighborhood $\Omega$ of $\xi_0$, there exists a small neighborhood $V_\Omega$ of 0 in $C^{p+2}(M)$ such that for every $v\in V_\Omega$, if $v\equiv0$ in $\Omega$, then $u(v)$ is degenerate i.e. there exists $\phi(v)\in K_v\backslash\left\{0\right\}$, where
$$K_v:=\{\phi\in H^2_1(M):\ \Delta_g\phi+h(v)\phi=(\crit-1)u(v)^{\crit-2}\phi\hbox{ in }M\}.$$
By renormalizing, we may assume that $\phi(v)\in\ss_{K_v}:=\{\phi\in K_v:\Vert \phi\Vert_{H^2_1}=1\}$. Since $h(tv),u(tv)\to h_0,u_0$ in $C^0(M)$ as $t\to0$, it then follows that there exists $\phi_v\in K_0$ and $(t_k)_{k\in\nn}>0$ such that $t_k\to0$ and $\phi(t_kv)\rightharpoonup\phi_v$ weakly in $H^2_1(M)$. By compactness of the embedding $H^2_1(M)\hookrightarrow L^2(M)$, we obtain that $\phi(t_kv)\to\phi_v$ strongly in $L^2(M)$. By elliptic theory that we apply to the linear equation satisfied by $\phi(t_kv)$, we then obtain that $\phi(t_kv)\to\phi_v$ strongly in $H^2_1(M)$, so that in particular $\phi_v\in\ss_{K_0}$. We then define 
$$\psi_k(v):=\frac{\phi(t_kv)-\phi_v}{t_k}\,.$$ 
It is easy to check that $\psi_k(v)$ satisfies the equation
\begin{equation}\label{eqpsik}
\Delta_g\psi_k(v)+h_0\psi_k(v)=(\crit-1)u_0^{\crit-2}\psi_k(v)+f_k(v)\phi(t_kv)\hbox{ in }M,
\end{equation}
where
\begin{align*}
f_k(v)&:=\frac{1}{t_k}((\crit-1)(u(t_kv)^{\crit-2}-u_0^{\crit-2})+h_0-h(t_kv))\\
&=\frac{1}{t_k}((\crit-2)(u(t_kv)^{\crit-2}-u_0^{\crit-2})+t_k\frac{u_0\Delta_gv-v\Delta_gu_0}{u_0u(t_kv)}).
\end{align*}
A straightforward Taylor expansion gives
\begin{equation}\label{fk}
f_k(v)=(\crit-2)^2u_0^{\crit-3}v+u_0^{-1}\Delta_gv-u_0^{-2}v\Delta_gu_0+\smallo(1)=u_0^{-1}L_0(v)+\smallo(1),
\end{equation}
as $k\to\infty$, uniformly in $v\in V_\Omega$, where 
\begin{equation}\label{f0}
L_0(v):=\Delta_gv+h_0v-(1-(\crit-2)^2)u_0^{\crit-2}v.
\end{equation}
It follows that 
$$\Vert\Pi_{K_0^\perp}(\Psi_k(v))\Vert_{H^2_1}\le C\Vert f_k(v)\phi(t_kv)\Vert_{\frac{2n}{n+2}} \le C\Vert\phi(t_kv)\Vert_{\frac{2n}{n+2}}\le C\Vert\phi(t_kv)\Vert_{H^2_1}\le C,$$  
where $\Pi_{K_0^\perp}$ is the orthogonal projection of $H^2_1$ onto $K_0^\perp$ and the letter $C$ stands for positive constants independent of $k\in\nn$ and  $v\in V_\Omega$. Since $(\Pi_{K_0^\perp}(\Psi_k(v)))_k$ is bounded in $H^2_1(M)$, up to a subsequence, we may assume that there exists $\psi_v\in K_0^\perp$ such that $\Pi_{K_0^\perp}(\Psi_k(v))\rightharpoonup\psi_v$ weakly in $H^2_1(M)$. Passing to the limit in \eqref{eqpsik} and using \eqref{fk}, we then obtain that $\psi_v$ satisfies the equation
$$\Delta_g\psi_v+h_0\psi_v=(\crit-1)u_0^{\crit-2}\psi_v+u_0^{-1}L_0(v)\phi_v\hbox{ in }M.$$
In particular, since $\phi_v\in K_0$, multiplying this equation by $\phi_v$ and integrating by parts yields
\begin{equation}\label{intv}
\int_Mu_0^{-1}L_0(v)\phi_v^2\,dv_g=0.
\end{equation}
We now construct $v$ contradicting \eqref{intv}. Let $\chi\in C^\infty(\rr)$ be such that $\chi(t)=0$ for $t\leq 1$ and $\chi(t)=1$ for $t\geq 2$. We define 
$$v_\eps(x):=\eps \chi(d_g(x,\xi_0)/\eps)u_0(x)\hbox{ for all }x\in M\hbox{ and }\eps>0.$$
As one checks, 
\begin{equation}\label{eq:33}
v_\eps\equiv 0\hbox{ in }B_\eps(\xi_0)\hbox{ and }\lim_{\eps\to 0}\eps^{-1}v_\eps= u_0\hbox{ in }L^p(M)\hbox{ for all }p\in [1,+\infty).
\end{equation}
Since $\Vert\phi_{v_\eps}\Vert_{H_1^2}=1$, $\phi_{v_\eps}\in K_0\subset C^2(M)$ and $K_0$ is of finite dimension, up to a subsequence, we can assume that there exists $\phi_0\in K_0$ such that
\begin{equation}\label{eq:34}
\lim_{\eps\to 0}\phi_{v_\eps}=\phi_{0}\neq 0\hbox{ in }C^2(M).
\end{equation}
Since $L_0$ is self-adjoint, it follows from \eqref{intv} that
$$\int_Mv_\eps L_0(u_0^{-1}\phi_{v_\eps}^2)\,dv_g=0\hbox{ for all }\eps>0.$$
Passing to the limit as $\eps\to 0$ in this equation and using \eqref{eq:33} and \eqref{eq:34}, we obtain
$$\int_M u_0 L_0(u_0^{-1}\phi_{0}^2)\,dv_g=0.$$
Integrating again by parts and noting that $L_0(u_0)=(\crit-2)^2u_0^{\crit-1}$, we then obtain
$$0=\int_M u_0^{-1}\phi_{0}^2L_0(u_0) \,dv_g=(\crit-2)^2\int_M u_0^{\crit-2}\phi_{0}^2\, dv_g,$$
which is a contradiction since $u_0>0$ and $\phi_0\not\equiv 0$. This ends the proof of Lemma~\ref{lemu0}.
\endproof

We can now end the proof of Theorem~\ref{th:3}. We let $\Omega_0$, $(\tilde{h}_k)_{k\in\nn}$ and $(\tilde{u}_k)_{k\in\nn}$ be given by Lemma~\ref{lemu0}. Since $\tilde{h}_k\equiv h_0$ and $\tilde{u}_k\equiv u_0$ in $\Omega_0$, we obtain that $\varphi_{\tilde{h}_k,\tilde{u}_k}\equiv\varphi_{h_0,u_0}$ in $\Omega_0$ and, therefore, $\varphi_{\tilde{h}_k,\tilde{u}_k}(\xi_0)=|\nabla\varphi_{\tilde{h}_k,\tilde{u}_k}(\xi_0)|=0$. For every $k\in\nn$, we can then mimick the first part of the proof of Theorem~\ref{th:2} to construct a sequence $(\tilde{h}_{k,j})_{j\in\nn}\in C^{\max(p,2)}(M)$ such that $\tilde{h}_{k,j}\to \tilde{h}_k$ in $C^p(M)$ as $j\to\infty$, $\varphi_{\tilde{h}_{k,j},\tilde{u}_k}(\xi_0)=0$, $\xi_0$ is a nondegenerate critical point of $\varphi_{\tilde{h}_{k,j},\tilde{u}_k}$ and $K_{\tilde{h}_{k,j},\tilde{u}_k}(\xi_0)\ne0$. We now distinguish two cases:

\smallskip\noindent
{\bf Case $n\ge7$.} Note that in this case, we have $\varphi_{\tilde{h}_{k,j},\tilde{u}_k}=\varphi_{\tilde{h}_{k,j}}$.  Since $\tilde{u}_k$ is nondegenerate and $\tilde{h}_{k,j}\to \tilde{h}_k$ in $C^1(M)$ as $j\to\infty$, the implicit function theorem gives that for large $j$, there exists a nondegenerate solution $\tilde{u}_{k,j}\in C^2(M)$ to the equation
$$\Delta_g\tilde{u}_{k,j}+\tilde{h}_{k,j}\tilde{u}_{k,j}=\tilde{u}_{k,j}^{\crit-1},\ \tilde{u}_{k,j}>0\hbox{ in }M$$
such that $\tilde{u}_{k,j}\to \tilde{u}_k$ in $C^2(M)$ as $j\to\infty$. By applying Theorem~\ref{th:main:u0}, we then obtain that there exist $\epsilon_{k,j}>0$, $(\tilde{h}_{k,j,\eps})_{0<\eps<\epsilon_{k,j}}\in C^{\max(p,2)}(M)$ and $(\tilde{u}_{k,j,\eps})_{0<\eps<\epsilon_{k,j}}\in C^2(M)$ satisfying
$$\Delta_g\tilde{u}_{k,j,\eps}+\tilde{h}_{k,j,\eps}\tilde{u}_{k,j,\eps}=\tilde{u}_{k,j,\eps}^{\crit-1}\hbox{ in }M,\ \tilde{u}_{k,j,\eps}>0\hbox{ for all }0<\eps<\epsilon_{k,j}$$
and such that $\tilde{h}_{k,j,\eps}\to\tilde{h}_{k,j}$ in $C^{\max(p,2)}(M)$, $\tilde{u}_{k,j,\eps}\rightharpoonup\tilde{u}_{k,j}$ weakly in $L^{\crit}(M)$ and $(\tilde{u}_{k,j,\eps})_\eps$ blows up with one bubble at $\xi_0$ as $\eps\to0$. Therefore, we obtain that for every $k\in\nn$, there exists $j_k\in\nn$ and $\eps'_k>0$ such that
$$\Vert \tilde{h}_{k,j_k}-\tilde{h}_k\Vert_{C^p}<\frac{1}{k},\ \Vert \tilde{u}_{k,j_k}-\tilde{u}_k\Vert_{C^2}<\frac{1}{k},\ 0<\eps'_{k}<\min(1/k, \eps_{k,j_k}),$$
$$\Vert \tilde{u}_{k,j_k,\eps'_k}-u_0\Vert_2<\frac{1}{k},\ \left|\int_{M}|\tilde{u}_{k,j_k,\eps'_k}-\tilde{u}_{k,j_k}|^{\crit} dv_g-\int_{\rn}U_{1,0}^{\crit}\, dx\right|<\frac{1}{k}$$
and
$$\int_{M\backslash B_{1/k}(\xi_0)}|\tilde{u}_{k,j_k,\eps'_k}-\tilde{u}_{k,j_k}|^{\crit} dv_g<\frac{1}{k}\,.$$
We then define $u_{k}:=\tilde{u}_{k,j_k,\eps'_k}$, so that 
$$\Delta_g u_k+h_k u_k=u_k^{\crit-1}\hbox{ in }M,\hbox{ where }h_k:=\tilde{h}_{k,j_k,\eps'_k}.$$
As one checks, $h_k\to h_0$ in $C^p(M)$, $u_k\rightharpoonup u_0$ weakly in $L^{\crit}(M)$ and $(u_k)_k$ blows up with one bubble at $\xi_0$ as $k\to\infty$. This proves Theorem~\ref{th:3} for $n\ge7$. 

\smallskip\noindent
{\bf Case $n=6$.} In this case, we have $\varphi_{\tilde{h}_{k,j},\tilde{u}_k}=\varphi_{\tilde{h}_{k,j}}-2\tilde{u}_k$. Furthermore, noting that $\crit-1=2$ when $n=6$, we can rewrite the equation \eqref{eqtilduk} as
$$\Delta_g\tilde{u}_k+(\tilde{h}_k-2\tilde{u}_k)\tilde{u}_k=-\tilde{u}_k^2\hbox{ in }M.$$
Since $\tilde{h}_{k,j}-2\tilde{u}_k\to \tilde{h}_k-2\tilde{u}_k$ in $C^0(M)$ as $j\to\infty$, a standard minimization method gives that for large $j$, there exists a unique nondegenerate solution $\tilde{u}_{k,j}$ to
$$\Delta_g\tilde{u}_{k,j}+(\tilde{h}_{k,j}-2\tilde{u}_k)\tilde{u}_{k,j}=-\tilde{u}_{k,j}^2,\ \tilde{u}_{k,j}>0\hbox{ in }M.$$
As is easily seen, this equation can be rewritten as 
\begin{equation}\label{eqtildukj}
\Delta_g\tilde{u}_{k,j}+\mathring{h}_{k,j}\tilde{u}_{k,j}=\tilde{u}_{k,j}^2,\ \tilde{u}_{k,j}>0\hbox{ in }M,\hbox{ where }\mathring{h}_{k,j}:=\tilde{h}_{k,j}-2\tilde{u}_k+2\tilde{u}_{k,j}.
\end{equation}
Since $\tilde{h}_{k,j}\to \tilde{h}_k$ in $C^p(M)$, we obtain that $\tilde{u}_{k,j}\to \tilde{u}_k$ and $\mathring{h}_{k,j}\to \tilde{h}_k$ in $C^p(M)$ as $j\to\infty$. Furthermore, since $\tilde{u}_k$ is nondegenerate, we have that $\tilde{u}_{k,j}$ is nondegenerate for large $j$. Similarly, since $K_{\tilde{h}_{k,j},\tilde{u}_k}(\xi_0)\ne0$, we obtain that $K_{\mathring{h}_{k,j},\tilde{u}_{k,j}}(\xi_0)\ne0$ for large $j$. Furthermore, we have
$$\varphi_{\mathring{h}_{k,j},\tilde{u}_{k,j}}=\mathring{h}_{k,j}-2\tilde{u}_{k,j}-c_n\Scal_g=\tilde{h}_{k,j}-2\tilde{u}_k-c_n\Scal_g=\varphi_{\tilde{h}_{k,j},\tilde{u}_k}.$$
In view of the properties satisfied by $\tilde{h}_{k,j}$, it follows that $\varphi_{\mathring{h}_{k,j},\tilde{u}_{k,j}}(\xi_0)=0$ and $\xi_0$ is a nondegenerate critical point of $\varphi_{\mathring{h}_{k,j},\tilde{u}_{k,j}}$. By applying Theorem~\ref{th:main:u0}, we then obtain that there exist $\epsilon_{k,j}>0$, $(\tilde{h}_{k,j,\eps})_{0<\eps<\epsilon_{k,j}}\in C^{\max(p,2)}(M)$ and $(\tilde{u}_{k,j,\eps})_{0<\eps<\epsilon_{k,j}}\in C^2(M)$ satisfying
$$\Delta_g\tilde{u}_{k,j,\eps}+\tilde{h}_{k,j,\eps}\tilde{u}_{k,j,\eps}=\tilde{u}_{k,j,\eps}^{\crit-1}\hbox{ in }M,\ \tilde{u}_{k,j,\eps}>0\hbox{ for all }0<\eps<\epsilon_{k,j}$$
and such that $\tilde{h}_{k,j,\eps}\to\mathring{h}_{k,j}$ in $C^{\max(p,2)}(M)$, $\tilde{u}_{k,j,\eps}\rightharpoonup\tilde{u}_{k,j}$ weakly in $L^{\crit}(M)$ and $(\tilde{u}_{k,j,\eps})_\eps$ blows up with one bubble at $\xi_0$ as $\eps\to0$. Finally, as in the previous case, we obtain the existence of $j_k\in\nn$ and $\epsilon'_k>0$ such that $u_{k}:=\tilde{u}_{k,j_k,\eps'_k}$ satisfies the desired conditions. This proves Theorem~\ref{th:3} for $n=6$.\qed
\endproof

\section{Necessity of the condition on the gradient}\label{sec:nec}

\begin{theorem}\label{th:cns}
Let $(M,g)$ be a compact Riemannian manifold of dimension $n\geq 4$. Let $h_0\in C^1(M)$ be such that $\Delta_g+h_0$ is coercive. Assume that there exist families $(h_\eps)_{\eps>0}\in C^p(M)$ and $(u_\eps)_{\eps>0}\in C^2(M)$ satisfying \eqref{eque} and such that $h_\eps\to h_0$ strongly in $C^1(M)$. Assume that $(M,g)$ is locally conformally flat. If $(u_\eps)_{\eps}$ blows up with one bubble at some point $\xi_0\in M$ and $u_\eps\rightharpoonup 0$ weakly as $\eps\to 0$, then  \eqref{phi0th2} holds true.
\end{theorem}

\smallskip\noindent{\it Proof of Theorem \ref{th:cns}.} Let $\varphi_{h_0}$ be as in \eqref{defphi0th3}. The identity $\varphi_{h_0}(\xi_0)=0$ is a consequence of the results of Druet~\cites{DruetJDG,DruetHDR}. Since $(M,g)$ is locally conformally flat, there exists $\Lambda\in C^\infty(M)$ positive such that $\hat{g}:=\Lambda^{\frac{4}{n-2}}g$ is flat around $\xi_0$. Define
$$\hue:=\Lambda^{-1}\ue\hbox{ and } \hhe:=\left(\he-c_n \Scal_g\right)\Lambda^{2-\crit}+c_n \Scal_{\hat{g}}.$$
The conformal law \eqref{conf:law} yields 
\begin{equation}\label{eq:hue}
\Delta_{\hat{g}}\hue +\hhe \hue =\hue^{\crit-1},\ \hue>0\hbox{ in }M.
\end{equation}
As one checks, on $(M,\hat{g})$, $\hue$ blows-up at $\xi_0$ in the sense that $\hue=U_{\delta_\eps,\xi_\eps}+\smallo(1)$ as $\eps\to0$ in $H^2_1(M)$, where $U_{\delta_\eps,\xi_\eps}$ is as in \eqref{def:peak} (with respect to the metric $\hat{g}$) and $(\delta_\eps,\xi_\eps)\to(0,\xi_0)$ as $\eps\to0$. It then follows from Druet--Hebey--Robert~\cite{DHR} that there exist $C,\eps_0>0$ such that for every $\eps\in(0,\eps_0)$,
\begin{equation}\label{est:c0:1}
\frac{1}{C} \left(\frac{\delta_\eps}{\delta_\eps^2+d_{\hat{g}}\left(x,\xi_\eps\right)^2}\right)^{\frac{n-2}{2}}\leq\hue(x)\leq C \left(\frac{\delta_\eps}{\delta_\eps^2+d_{\hat{g}}\left(x,\xi_\eps\right)^2}\right)^{\frac{n-2}{2}}
\end{equation} 
for all $x\in M$ and, defining 
$$U_\eps(x):=\delta_\eps^{\frac{n-2}{2}}\chi(x)\hue(\xi_\eps+\delta_\eps x)$$  
for all $x\in\rn$, where $\chi$ is a cutoff function on a small ball centered at $\xi_0$, we have
\begin{equation}\label{cv:Ue}
\lim_{\eps\to 0}U_\eps=U_{1,0}=\Bigg(\frac{\sqrt{n(n-2)}}{1+|\cdot|^2}\Bigg)^{\frac{n-2}{2}}\hbox{ in }C^2_{loc}(\rn).
\end{equation}
Without loss of generality, via a chart, we may assume that $\hat{g}$ is the Euclidean metric on $B_{2\nu}(\xi_0)$ for some $\nu>0$. We fix $i\in \{1,...,n\}$. By differentiating the Pohozaev identity for $\hue$ on $B_{\nu}(\xi_\eps)$ (see for instance Ghoussoub--Robert~\cite{GR}*{Proposition~7}) 
and integrating by parts, we obtain
\begin{multline}\label{eq:poho}
\frac{1}{2}\int_{B_{\nu}(\xi_\eps)}\partial_i\hhe \hue^2\,dx\\
=\int_{\partial B_{\nu}(\xi_\eps)}\left(\frac{x_i}{\left|x\right|}\left(\frac{|\nabla \hue|^2+\hhe\hue^2}{2}-\frac{\hue^{\crit}}{\crit}\right)-\left<\frac{x}{\left|x\right|},\nabla\hue\right>\partial_i\hue\right) d\sigma(x),
\end{multline} 
where $d\sigma$ is the volume element on $\partial B_{\nu}(\xi_\eps)$. It follows from \eqref{est:c0:1} that there exists $C(\nu)>0$ such that $\hue(x)\leq C(\nu)\delta_\eps^{\frac{n-2}{2}}$ for all $x\in M\setminus B_{\nu/4}(\xi_0)$ and $\eps\in(0,\eps_0)$. It then follows from \eqref{eq:hue} and elliptic theory that there exists $C_1>0$ such that $|\nabla \hue(x)|\leq C_1 \delta_\eps^{\frac{n-2}{2}}$ for all $x\in M\setminus B_{\nu/2}(\xi_0)$ and $\eps\in(0,\eps_0)$.  Plugging these inequalities into \eqref{eq:poho} yields
\begin{equation}\label{bnd:poho}
\int_{B_{\nu}(\xi_\eps)}\partial_i\hhe \hue^2\, dx=\bigO(\delta_\eps^{n-2})\hbox{ as }\eps\to 0.
\end{equation}
On the other hand, with a change of variable, we obtain
$$\int_{B_{\nu}(\xi_\eps)}\partial_i\hhe \hue^2\, dx=\delta_\eps^2\int_{B_{\nu/\delta_\eps}(0)}(\partial_i\hhe)(\xi_\eps+\delta_\eps x) U_\eps(x)^2\, dx.$$
The control \eqref{est:c0:1} gives $U_\eps\leq C U_{1,0}$. Therefore, when $n\geq 5$, Lebesgue's dominated convergence Theorem and \eqref{cv:Ue} yield
$$\int_{B_{\nu}(\xi_\eps)}\partial_i\hhe \hue^2\, dx=\delta_\eps^2\left( \partial_i\hhe(\xi_\eps)\int_{\rn}U_{1,0}^2\,dx +\smallo(1)\right)\text{ as }\eps\to0.$$
Combining this identity with \eqref{bnd:poho}, we obtain that $\partial_i(\varphi_{h_0}\Lambda^{2-\crit})(\xi_0)=0$ when $n\geq 5$. Since $\Lambda>0$ and $\varphi_{h_0}(\xi_0)=0$, it follows that $\partial_i\varphi_{h_0}(\xi_0)=0$ when $n\geq 5$.

\smallskip
We now assume that $n=4$. With \eqref{est:c0:1}, we obtain  
$$\int_{B_{\nu}(\xi_\eps)}|x-\xi_\eps| \hue^2\, dx=\bigO(\delta_\eps^{2}).$$ 
Therefore, with \eqref{bnd:poho}, we obtain
\begin{equation}\label{ineq1}
\partial_i\hhe(\xi_\eps)=\bigO\left(\delta_\eps^2\left(\int_{B_{\nu}(\xi_\eps)}\hue^2\, dx\right)^{-1}\right).
\end{equation}
With the lower bound in \eqref{est:c0:1}, we then obtain
\begin{equation}\label{ineq2}
\int_{B_{\nu}(\xi_\eps)}\hue^2\, dx\geq C\int_{B_{\nu}(\xi_\eps)}\left(\frac{\delta_\eps}{\delta_\eps^2+|x-\xi_\eps|^2}\right)^{n-2} dx\geq C \delta_\eps^{2}\ln(1/\delta_\eps).
\end{equation}
It follows from \eqref{ineq1} and \eqref{ineq2} that $\partial_i\hhe(\xi_\eps)=\smallo(1)$ as $\eps\to0$ and so again $\partial_i\varphi_{h_0}(\xi_0)=0$ when $n=4$. 

\smallskip 
In all cases, we thus obtain that $\nabla\varphi_{h_0}(\xi_0)=0$. 
This ends the proof of Theorem~\ref{th:cns}.\qed

\end{document}